\documentclass[11pt,a4paper,oneside]{amsart}
\pdfoutput=1

\usepackage[normalem]{ulem}
\usepackage{multirow}

\usepackage{tcolorbox}

\usepackage[arrow,curve,matrix,tips,2cell]{xy}

\usepackage[hmargin=2cm,vmargin=1.5cm]{geometry}

\usepackage{bm}
\usepackage{fancyhdr}

\usepackage{amsmath}
\usepackage{amscd}
\usepackage{amssymb}
\usepackage{latexsym}

\usepackage{graphicx}

\usepackage{mathrsfs}

\usepackage{bbm}

\usepackage[shortlabels]{enumitem}
\usepackage{hyperref}
\usepackage{comment}

\newtheorem{theorem}{Theorem}[section]
\newtheorem*{theorem*}{Theorem}
\newtheorem{lemma}[theorem]{Lemma}
\newtheorem*{lemma*}{Lemma}
\newtheorem{corollary}[theorem]{Corollary}
\newtheorem{conjecture}[theorem]{Conjecture}
\newtheorem{proposition}[theorem]{Proposition}

\newtheorem{remark}[theorem]{Remark}
\newtheorem{remarks}[theorem]{Remarks}
\newtheorem{definition}[theorem]{Definition}
\newtheorem*{definition*}{Definition}

\newtheorem{question}[theorem]{Question}
\newtheorem*{question*}{Question}

\newtheorem{example}[theorem]{Example}
\newtheorem{examples}[theorem]{Examples}

\makeatletter
\def\revddots{\mathinner{\mkern1mu\raise\p@
\vbox{\kern7\p@\hbox{.}}\mkern2mu
\raise4\p@\hbox{.}\mkern2mu\raise7\p@\hbox{.}\mkern1mu}}
\makeatother 
\newcommand{\bgl}{\begin{equation}} 
\newcommand{\egl}{\end{equation}}
\newcommand{\bgloz}{\begin{equation*}} 
\newcommand{\egloz}{\end{equation*}}
\newcommand{\bgln}{\begin{eqnarray}} 
\newcommand{\egln}{\end{eqnarray}}
\newcommand{\bglnoz}{\begin{eqnarray*}} 
\newcommand{\eglnoz}{\end{eqnarray*}}
\newcommand{\btheo}{\begin{theorem}}
\newcommand{\etheo}{\end{theorem}}
\newcommand{\btheooz}{\begin{theorem*}}
\newcommand{\etheooz}{\end{theorem*}}

\newcommand{\blemma}{\begin{lemma}}
\newcommand{\elemma}{\end{lemma}}
\newcommand{\blemmaoz}{\begin{lemma*}}
\newcommand{\elemmaoz}{\end{lemma*}}
\newcommand{\bproof}{\begin{proof}}
\newcommand{\eproof}{\end{proof}}
\newcommand{\bbew}{\begin{beweis}}
\newcommand{\ebew}{\end{beweis}}
\newcommand{\bremark}{\begin{remark}}
\newcommand{\eremark}{\end{remark}}
\newcommand{\bremarks}{\begin{remarks}}
	\newcommand{\eremarks}{\end{remarks}}
\newcommand{\bdefin}{\begin{definition}}
\newcommand{\edefin}{\end{definition}}
\newcommand{\bdefinoz}{\begin{definition*}}
\newcommand{\edefinoz}{\end{definition*}}
\newcommand{\bex}{\begin{example}}
\newcommand{\eex}{\end{example}}
\newcommand{\bexs}{\begin{examples}}
\newcommand{\eexs}{\end{examples}}

\newcommand{\bprop}{\begin{proposition}}
\newcommand{\eprop}{\end{proposition}}
\newcommand{\bcor}{\begin{corollary}}
\newcommand{\ecor}{\end{corollary}}
\newcommand{\bfa}{\begin{cases}} 
\newcommand{\efa}{\end{cases}}

\newcommand{\bquestion}{\begin{question}}
\newcommand{\equestion}{\end{question}}
\newcommand{\bquestionoz}{\begin{question*}}
\newcommand{\equestionoz}{\end{question*}}
\newcommand{\bconj}{\begin{conjecture}}
\newcommand{\econj}{\end{conjecture}}
\newcommand{\cB}{\mathcal B}
\newcommand{\cC}{\mathcal C}

\newcommand{\cE}{\mathcal E}
\newcommand{\cF}{\mathcal F}
\newcommand{\cG}{\mathcal G}

\newcommand{\cI}{\mathcal I}
\newcommand{\cJ}{\mathcal J}
\newcommand{\cK}{\mathcal K}

\newcommand{\cO}{\mathcal O}
\newcommand{\cP}{\mathcal P}

\newcommand{\cR}{\mathcal R}

\newcommand{\cU}{\mathcal U}

\def\Cz{\mathbb{C}}

\def\Kz{\mathbb{K}}

\def\Nz{\mathbb{N}}

\def\Qz{\mathbb{Q}}

\def\Tz{\mathbb{T}}
\def\Zz{\mathbb{Z}}

\def\1z{\mathbb{1}}
\newcommand{\fA}{\mathfrak A}

\newcommand{\fD}{\mathfrak D}
\newcommand{\fE}{\mathfrak E}
\newcommand{\fF}{\mathfrak F}

\newcommand{\fU}{\mathfrak U}

\newcommand{\fX}{\mathfrak X}

\newcommand{\fZ}{\mathfrak Z}
\newcommand{\an}[1]{``#1''} 
\newcommand{\ti}{\tilde}

\newcommand{\ma}{\mapsto} 
\newcommand{\onto}{\twoheadrightarrow} 
\newcommand{\into}{\hookrightarrow} 
\newcommand{\Rarr}{\Rightarrow} 
\newcommand{\Larr}{\Leftarrow} 
\newcommand{\LRarr}{\Leftrightarrow} 

\def\SEMI{\mbox{$\times\kern-2pt\vrule height5pt width.6pt \kern3pt $}}

\newcommand{\End}{{\rm End}\,}
\newcommand{\Aut}{{\rm Aut}}

\newcommand{\id}{{\rm id}}

\newcommand{\Ind}{\mathrm{ Ind}\,}

\renewcommand{\ker}{{\rm ker}\,}

\newcommand{\reg}{^\times} 
\newcommand{\Iso}{{\rm Isom}} 
\newcommand{\defeq}{\mathrel{:=}} 

\newcommand{\dop}{\text{: }} 
\newcommand{\ilim}{\varinjlim} 
\newcommand{\plim}{\varprojlim} 
\newcommand{\supp}{{\rm supp}} 
\newcommand{\Ell}{{\rm L}} 
\newcommand{\lge}{\left\{} 
\newcommand{\rge}{\right\}} 
\newcommand{\lsp}{\left\langle} 
\newcommand{\rsp}{\right\rangle} 
\newcommand{\gekl}[1]{\lge #1 \rge} 
\newcommand{\spkl}[1]{\lsp #1 \rsp} 
\newcommand{\menge}[2]{\gekl{ #1 \dop #2 }} 
\newcommand{\mfc}{\mathfrak{c}}
\newcommand{\mfg}{\mathfrak{g}}

\renewcommand{\c}{\mathfrak{c}}

\newcommand{\f}{\mathfrak{f}}
\newcommand{\p}{\mathfrak{p}}
\newcommand{\q}{\mathfrak{q}}
\newcommand{\m}{\mathfrak{m}}

\newcommand{\acts}{\curvearrowright} 

\newcommand{\spn}{\textup{span}}

\newcommand{\ann}{\textup{ann}}

\newcommand{\Isom}{\textup{Isom}}

\newcommand{\M}{\textup{M}}

\newcommand{\gp}[1]{\langle #1\rangle}
\newcommand{\mon}[1]{\langle #1\rangle^+}

\newcommand{\dom}{\textup{dom}}

\newcommand{\PIsom}{\textup{PIsom}}
\newcommand{\Reg}{\cR} 
\newcommand{\Asc}{\textup{Asc}}
\newcommand{\Ann}{\textup{Ann}}

\renewcommand{\Iso}{\textup{Iso}}
\newcommand{\im}{\textup{im}}
\newcommand{\fix}{\textup{fix}}

\newcommand{\E}{\widehat{\cE}}

\newcommand{\bbt}{\mathbbm{t}}

\newcommand{\bd}{\partial \widehat{\cE}}
\newcommand{\msc}[1]{\mathscr{#1}}
\newcommand{\es}{\textup{ess}}
\newcommand{\sing}{\textup{sing}}
\newcommand{\eps}{\varepsilon}
\newcommand{\core}{\textup{core}_G}

\begin{document}

\title{Algebraic actions I. C*-algebras and groupoids}

\thispagestyle{fancy}

\author{Chris Bruce}
\thanks{C. Bruce was supported by a Banting Fellowship administered by the Natural Sciences and Engineering Research Council of Canada (NSERC) and has received funding from the European Union’s Horizon 2020 research and innovation programme under the Marie Sklodowska-Curie grant agreement No 101022531. X. Li has received funding from the European Research Council (ERC) under the European Union’s Horizon 2020 research and innovation programme (grant agreement No. 817597).}
\address{School of Mathematics and Statistics, University of Glasgow, University Place, Glasgow G12 8QQ, United Kingdom}
\email{Chris.Bruce@glasgow.ac.uk}

\author{Xin Li}
\address{School of Mathematics and Statistics, University of Glasgow, University Place, Glasgow G12 8QQ, United Kingdom}
\email{Xin.Li@glasgow.ac.uk}


\subjclass[2020]{Primary 46L05, 46L55, 37A55; Secondary 20M18, 22A22, 37B99.}

\begin{abstract}
Given an algebraic action of a semigroup, we construct an inverse semigroup, and we characterize Hausdorffness, topological freeness, and minimality of the associated tight groupoid in terms of conditions on the initial algebraic action. We parameterize all closed invariant subspaces of the unit space of our groupoid, and characterize topological freeness of the associated reduction groupoids. We prove that our groupoids are purely infinite whenever they are minimal, which answers a general open question in the affirmative for our special class of groupoids.
In the topologically free case, we prove that the concrete C*-algebra associated with the algebraic action is always a (possibly exotic) groupoid C*-algebra in the sense that it sits between the full and essential C*-algebras of our groupoid. This provides a framework for studying such concrete C*-algebras, allowing us to obtain structural results that were only previously available for very special classes of algebraic actions. For instance, we obtain results on simplicity and pure infiniteness for C*-algebras associated with subshifts over semigroups, actions coming from commutative algebra, and non-commutative rings. These results were out of reach using existing techniques.
\end{abstract}

\maketitle

\parskip=0em


\setlength{\parindent}{0cm} \setlength{\parskip}{0.5cm}

\section{Introduction}

\subsection{Context}
\emph{Algebraic actions of groups} form an important class of dynamical systems with deep connections to commutative algebra and operator algebras; they have been studied intensely since their inception by Kitchens and Schmidt in the late 80s (see \cite{LindSchmidt}, \cite{Sch}, \cite[Chapters~13\&14]{KL} and the references therein).
The theory of their ``one-sided'' or ``irreversible'' counterparts, \emph{algebraic actions of semigroups}, is much less developed.
This work is the first in a series of papers in which we systematically study algebraic actions of semigroups from the perspective of operator algebras and groupoids. Here, we consider algebraic actions through the lens of C*-algebra and groupoid theory. Each such algebraic action gives rise to a concrete C*-algebra generated by the Koopman representation for the action together with the left regular representation of the group, and we develop a general framework for studying these C*-algebras using \'{e}tale groupoids, motivated by natural examples arising from a variety of sources, with the goal of developing a better understanding for their dynamical properties. The observation that algebraic actions of semigroups give rise to new constructions of C*-algebras and groupoids which turn out to be interesting in their own right provides further motivation for our work.

An algebraic action of a semigroup is an action of a semigroup on a group by injective endomorphisms. Until now, the algebraic actions that have been investigated from the point of view of C*-algebra theory fall into two classes.
The first class comes from ring theory: The multiplicative monoid of left regular elements of a ring acts by injective endomorphisms on the additive group of the ring, and the associated C*-algebra is the reduced ring C*-algebra of the ring. These were introduced and studied by Cuntz for the ring $\Zz$ \cite{Cuntz}, for integral domains with finite quotients by Cuntz and the second-named author in \cite{CuLi}, and for general rings by the second-named author in \cite{Li:Ring}. The ring C*-algebras of integral domains with finite quotients that are not fields were proven to be UCT Kirchberg algebras in \cite{CuLi} (following the approach in \cite{Cuntz} for the ring $\Zz$). For general commutative rings, conditions for pure infiniteness and simplicity were established in \cite{Li:Ring}, but the question of pure infiniteness and simplicity for non-commutative rings was left open.
Closely related to these are the C*-algebras associated with actions of congruence monoids on rings of algebraic integers, which have been studied by the authors in the context of boundary quotients of semigroup C*-algebras, and were proven to be UCT Kirchberg algebras, see \cite[\S~8]{Bru1} and \cite[\S~3]{BruLi}.
The second class comes from special actions of right LCM semigroups (i.e., semigroups in which all elements which have a right common multiple have a least right common multiple). Their study began with Hirshberg's work on the concrete C*-algebras associated with algebraic $\Nz$-actions on amenable groups \cite{Hirsh}. Cuntz and Vershik later proved that these C*-algebras are UCT Kirchberg algebras for exact, finite cokernel actions on Abelian groups \cite{CV} (the finite cokernel assumption was later weakened by Vieira \cite{Vieira1}). 
These works led Stammeier to introduce what are called \emph{irreversible algebraic dynamical systems} in \cite[Definition~1.5]{Sta}, which are special actions of free Abelian semigroups by endomorphisms on groups. Results on pure infiniteness and simplicity were obtained in this setting (see \cite[Theorem~3.26]{Sta} and \cite[Theorem~5.10]{Sta2}).
The actions studied in \cite{Sta} were generalized by the so-called \emph{algebraic dynamical systems} introduced by Brownlowe, Larsen, and Stammeier in \cite[Definition~2.1]{BLS2}. 
These are algebraic actions where the acting semigroup is required to be right LCM and the action is required to \emph{respect the order} in the sense of \cite[Definition~8.1]{BLS}. These conditions are quite restrictive, e.g., for the examples coming from rings, they are tantamount to requiring that the ring is a principal ideal domain. Yet, this class is already broad enough to include many interesting examples which have been studied also in, e.g., \cite{Sta2} and \cite{BS}.
These authors were primarily focused on universal C*-algebras attached to the actions and the relationship with (boundary quotients of) semigroup C*-algebras. In the case where these C*-algebras agree with the boundary quotients of the semi-direct product semigroup attached to the action, results on pure infiniteness and simplicity were obtained in  \cite[Theorem~8.12]{BLS} and \cite[Theorem~4.17]{BS}. One upshot of the present work is that we are able to obtain structural results even when tools from the theory of semigroup C*-algebras are not available. This allows us to treat many example classes of C*-algebras whose structure was previously intractable.

The aforementioned actions do not include several fundamental example classes: For instance, actions on solenoids, shifts over semigroups, and algebraic $\Nz^d$-actions on modules over multivariable polynomial rings, which are the one-sided analogues of the actions from Schmidt's classical book \cite{Sch}, are far from the ring-theoretic examples and typically fail to be algebraic dynamical systems in the sense of \cite[Definition~2.1]{BLS2}. 

\subsection{Motivation} Our motivation is twofold: First, it comes from the observation above that many of the most natural examples of algebraic actions do not fit into any of the existing frameworks, so that structural properties of the associated C*-algebras are not accessible with existing technologies; second, we wanted to find a good construction of an \'{e}tale groupoid from an algebraic action since the theory of \'{e}tale groupoids is of interest independent of C*-algebraic considerations.

The groupoid approach to studying a class of C*-algebras is by now a standard one: In order to analyze a class of C*-algebras constructed from some algebraic or combinatorial data, one finds groupoid models for the C*-algebras in question, i.e., one constructs a groupoid from the underlying data and then compares the C*-algebras of the groupoid to the initial C*-algebra of interest. Since there are powerful tools for the structure of groupoid C*-algebras, this often reduces many C*-algebraic questions to problems about groupoids, where ideas from dynamics can be employed. This is the general strategy used here to study the concrete C*-algebras associated with algebraic actions, however we stress that it differs significantly at the technical level from previous works for C*-algebras constructed from graphs, semigroups, or more generally left cancellative small categories.

\subsection{Novelty}
In this paper, we establish a general framework for studying C*-algebras and groupoids associated with general algebraic actions that includes these important classes, where previous approaches are not applicable, as well as the examples from ring theory and algebraic dynamical systems mentioned above. Our idea is to focus on actions of inverse semigroups which naturally arise from algebraic actions of semigroups. On the one hand, the corresponding groupoids give us access to various C*-algebras attached to algebraic actions. At this point, we encounter the interesting phenomenon that essential groupoid C*-algebras or even exotic groupoid C*-algebras naturally enter our discussion. On the other hand, we succeed in establishing structural properties of our groupoids for general classes of algebraic actions, which then have consequences for the C*-algebras we are interested in. In addition, our analysis paves the way for the discovery that the groupoids we associate with algebraic actions are interesting in their own right as they often exhibit surprising rigidity phenomena; this is explained in the subsequent paper \cite{BruLi:rigidity}. Compared to previous approaches, the key novelty of our work is that we do not insist on a strong connection to semigroup C*-algebras because universal C*-algebras arising from the semigroup C*-algebra approach do not provide good models for algebraic actions, as we demonstrate in this paper.

Our approach leads us naturally to a notion of \emph{exactness} for general algebraic actions (see Definition~\ref{def:exact}), which is a vast generalization of Rohlin's notion of exactness for a single endomorphism (see \cite{Rohlin}). This notion is important in our work from the perspective of C*-algebras and groupoids, but it also seems natural from the dynamical point of view.

\subsection{Overview of results}
Let us now explain the main construction of this paper (see \S~\ref{sec:C*andgpoids}) and our results on groupoids. Given any algebraic action $\sigma\colon S\acts G$ (see Definition~\ref{def:Saction}), we construct an inverse semigroup as follows: Let $I_\sigma$ be the inverse semigroup of partial bijections of $G$ generated by the endomorphisms implementing the action of $S$ together with the translations for the action of $G$ on itself (see Definition~\ref{def:invsemigp}). The action of $S$ on $G$ generates a distinguished family of subgroups of $G$, which we call \emph{$S$-constructible subgroups}.
The collection of non-zero idempotents of $I_\sigma$ is then equal to the collection of cosets for these subgroups. This collection of cosets defines a topology on $G$, and the resulting completion is a compact, totally disconnected space $\bd$, which carries a natural action of $I_{\sigma}$. The groupoid we build from the algebraic action $\sigma$ is then given by the associated transformation groupoid $\cG_\sigma:=I_\sigma\ltimes\bd$, which coincides with the tight groupoid of $I_\sigma$ in the sense of Exel \cite{Exel08,Exel09}. $\bd$ and its $I_{\sigma}$-action have natural interpretations from the perspective of C*-algebras: There is a canonical faithful representation of $I_{\sigma}$ by partial isometries on $\ell^2(G)$, and the C*-algebra generated by this representation is precisely the C*-algebra $\fA_\sigma$ generated by the Koopman representation of $S$ and the left regular representation of $G$. Now the Gelfand spectrum of the commutative C*-algebra $\fD_\sigma$ of $\fA_\sigma$ generated by the projections coming from idempotents in $I_\sigma$ is canonically homeomorphic to $\bd$. Furthermore, under Gelfand duality, the $I_{\sigma}$-action corresponds to the conjugation action by the partial isometries representing $I_{\sigma}$ on $\ell^2(G)$. As a next step, we carefully analyse our groupoids $\cG_{\sigma}$ and, among other things, establish the following general structural results (see \S~\ref{sec:properties}): 
\setlength{\parindent}{0cm} \setlength{\parskip}{0cm}

\begin{itemize}
    \item We completely characterize when $\cG_{\sigma}$ is Hausdorff, topologically free, or minimal in terms of the algebraic action $\sigma$.
    \item We parameterize all closed invariant subspaces of $\bd=\cG_\sigma^{(0)}$ and characterize topological freeness for the corresponding reduction groupoids, again in terms of $\sigma$.
    \item We prove that $\cG_\sigma$ is purely infinite as soon as it is minimal.
\end{itemize}

The third point can be reformulated by saying that as soon as $\cG_\sigma$ is minimal, it satisfies comparison in the sense of \cite{Ker,MW}, which is an important property in the classification programme for C*-algebras as well as in the study of topological full groups (see Question~\ref{ques:pi} and the discussion surrounding it). We then turn to the question of in what sense $\cG_{\sigma}$ is a groupoid model for our concrete C*-algebra $\fA_{\sigma}$. First, we prove that there is always a canonical, surjective *-homomorphism $C^*(\cG_\sigma)\to\fA_\sigma$ that is an isomorphism at the level of canonical commutative subalgebras. Thus $C^*(\cG_\sigma)$ can be viewed as a universal model for $\fA_\sigma$. At the level of concrete C*-algebras, it turns out that, interestingly, the essential C*-algebra $C_\es^*(\cG_\sigma)$ of $\cG_\sigma$ often provides the best approximation for $\fA_\sigma$. This contrasts with the case of C*-algebras associated to left cancellative small categories where the concrete C*-algebra is typically a quotient of the reduced C*-algebra of the groupoid model, see \cite{Li:GarI}. In \S~\ref{sec:C*-Comparison}, we use recent results by Christensen and Neshveyev on induced representations for groupoid C*-algebras \cite{CN} to prove that under mild assumptions, the concrete C*-algebra $\fA_\sigma$ surjects onto the essential C*-algebra of $\cG_\sigma$, and injectivity of this map is characterized by amenability of the homomorphic group image of the inverse semigroup $I^e$ generated by the endomorphisms coming from $S$. In the case where $\cG_{\sigma}$ is Hausdorff, reduced and essential groupoid C*-algebras coincide, so that we obtain criteria when $\fA_{\sigma}$ can be identified with $C^*_r(\cG_{\sigma})$. The essential C*-algebra of  an \'etale groupoid, as defined in \cite{ExelPitts,KM}, is an innovation at the heart of recent advances in the study of C*-algebras of non-Hausdorff groupoids (see, e.g., \cite{CEPSS,ExelPitts,KM,KKLRU,Li:GarI,CN,NS}). It captures---at the C*-algebraic level---minimality and pure infinitness of the groupoid even in the non-Hausdorff setting. Our identification of $\fA_{\sigma}$ with $C_\es^*(\cG_\sigma)$ allows us to deduce results on structural properties of $\fA_{\sigma}$ from our above-mentioned analysis of $\cG_{\sigma}$. For instance, under amenability assumptions, we characterize simplicity and pure infiniteness of $\fA_\sigma$ in terms of the initial action $\sigma \colon S\acts G$. Interestingly, we also identify criteria when $\fA_{\sigma}$ is an exotic groupoid C*-algebra lying properly between the full and reduced C*-algebras of $\cG_{\sigma}$ (see Corollary~\ref{cor:exotic}). Moreover, we compare $C_\es^*(\cG_\sigma)$ and $\fA_{\sigma}$ with C*-algebras from the context of semigroup C*-algebras (see \S~\ref{sec:semigpC*}) and show that the latter do not provide good models for $\fA_{\sigma}$ in general. 
\setlength{\parindent}{0cm} \setlength{\parskip}{0.5cm}

Having established our general framework, we turn to example classes. In \S~\ref{sec:FI}, we consider actions satisfying what we call the \emph{finite index property}. This condition means that every constructible subgroup is of finite index in $G$. Such actions include endomorphisms of (duals of) solenoids, examples from self-similar groups, and actions from torsion-free rings of finite rank. For such actions, our groupoid is always minimal and purely infinite. In particular, the reduced ring C*-algebra of a torsion-free ring $\cR$ of finite rank whose additive group does not contain a copy of $\Qz$ is simple if and only if the group of units in the $\Qz$-algebra $\Qz\otimes \cR$ is amenable. This is the first general structural result for ring C*-algebras of non-commutative rings, addressing a natural question left open in \cite{Li:Ring}. 
We then turn to actions by left reversible monoids in \S~\ref{sec:leftreversible}. This covers several large classes of actions, including actions from commutative algebra, which are the semigroup versions of the actions from \cite{Sch}. In this setting, too, our groupoids are always minimal and purely infinite. Even for the case of reduced ring C*-algebras of commutative rings, our results offer improvements on those from \cite{Li:Ring}.
We then consider shifts over semigroups in \S~\ref{sec:shifts}. The C*-algebras associated with the full $\Nz$-shift over $\Zz/n\Zz$ is the Cuntz algebra $\cO_n$, so the full $S$-shift over a group $\Sigma$ is a higher-dimensional analogue of $\cO_n$. For shifts over left reversible semigroups, we establish criteria when their C*-algebras are UCT Kirchberg algebras. We explicitly compute K-theory by applying the general K-theory formula from \cite{Li:Ktheory}. This allows to conclude that, surprisingly, many of these C*-algebras are actually abstractly isomorphic (see \S~\ref{sss:KTheoryClass}). For shifts over semigroups which are not left reversible, our C*-algebras are non-simple. We show that in many cases, there is a unique non-zero, proper ideal in our C*-algebras, which is given by the algebra of compact operators, and we prove that the minimal quotient is purely infinite and simple under mild assumptions on the action (see \S~\ref{sss:nonsimple}).

\textbf{Acknowledgments.}
C. Bruce would like to thank Bartosz Kosma Kwa\'{s}niewski for helpful comments on essential C*-algebras of non-Hausdorff groupoids, and Kevin Aguyar Brix and Gavin Goerke for inspiring discussions on non-Hausdorff groupoids and C*-algebras from actions of inverse semigroups. We also thank Julian Kranz for helpful comments which led to an improvement of Lemma~\ref{lem:PI}.

\section{Preliminaries}

\subsection{Algebraic actions of semigroups}
\label{sec:algactions}

This work is centered around the following class of actions, which are the ``one-sided'' analogues of the algebraic actions of groups from \cite{Sch}.
\bdefin
\label{def:Saction}
An \emph{algebraic action of a semigroup $S$} is an action of $S$ on a discrete group $G$ by injective group endomorphisms, i.e., a semigroup homomorphism $\sigma\colon S\to\End(G)$ such that $\sigma_s$ is an injective group homomorphism $G \to G$ for all $s\in S$.
\edefin
We shall write $\sigma\colon S\acts G$ to denote an algebraic action. When we want to emphasize the acting semigroup, we call $\sigma\colon S\acts G$ an \emph{algebraic $S$-action}. 

The case where $G$ is Abelian is the most interesting from the point of view of topological dynamics:

\bremark
When $G$ is Abelian, an algebraic $S$-action on $G$ is the same as having an action of $S$ on the compact Pontryagin dual $\widehat{G}$ of $G$ by surjective group endomorphisms.
More precisely, if $\sigma\colon S\acts G$ is an algebraic $S$-action, then we obtain a right action $\hat{\sigma}\colon S\acts \widehat{G}$ (i.e., a semigroup antihomomorphism from $S$ to the semigroup of continuous homomorphisms $\widehat{G} \to \widehat{G}$) characterized by $\spkl{\hat{\sigma}_s(\chi),g} = \spkl{\chi,\sigma_s(g)}$ for all $s \in S$, $\chi \in \widehat{G}$ and $g \in G$, where $\spkl{\cdot,\cdot} \colon \widehat{G} \times G \to \Tz$ is the duality pairing. Moreover, it is easy to see that, for all $s \in S$, $\hat{\sigma}_s$ is surjective if and only if $\sigma_s$ is injective.
\eremark

\textbf{Standing assumptions:} Assume that $\sigma \colon S\acts G$ is an algebraic $S$-action. 
The identity element of $G$ will be denoted by $e$. We shall always assume that $S$ is a monoid with identity element denoted by $1$ and that $G$ is non-trivial. Note that since $\sigma$ is injective, we have $\sigma_1=\id$. We will also always assume that the action $\sigma \colon S\acts G$ is \emph{faithful}, i.e., $\sigma$ is injective. This implies that $S$ is left cancellative.

We say $\sigma\colon S\acts G$ is \emph{automorphic} if $\sigma_s\in\Aut(G)$ for all $s\in S$; thus, $\sigma\colon S\acts G$ is non-automorphic if there exists $s\in S$ such that $\sigma_sG\lneq G$ (this implies that $S$ is non-trivial). 
We will primarily be interested in non-automorphic actions.

Let us now introduce the important concept of globalization.
\bdefin
\label{def:global}
We say that the algebraic action $\sigma \colon S\acts G$ has a \emph{globalization} if $S$ can be embedded into a group $\mathscr{S}$ and there is a group $\mathscr{G}$ containing $G$ together with an algebraic action $\tilde{\sigma} \colon \mathscr{S} \acts \mathscr{G}$ (which is necessarily automorphic) such that $\tilde{\sigma}_s\vert_G=\sigma_s$ for all $s\in S$.
\edefin

There are many large example classes of algebraic actions that admit a globalization.

\bex
Suppose $S$ is left Ore, i.e., cancellative and right reversible ($Ss \cap St \neq \emptyset$ for all $s, t \in S$). Then $S$ acts on the group $\mathscr{G}=S^{-1}G:=\ilim_s\{G\overset{\sigma_s}{\to}G\}$ by automorphisms, so we obtain a globalization $\ti{\sigma }\colon \mathscr{S}\acts S^{-1}G$, where $\mathscr{S} = S^{-1} S$ is the enveloping group of $S$. More precisely, consider the inductive system $\mathscr{G}_p \defeq G$ (for $p \in S$, where $p \leq q$ if $q = rp$ for some $r \in S$), with connecting map $\mathscr{G}_p \to \mathscr{G}_{rp}$ given by $\sigma_r \colon G = \mathscr{G}_p \to \mathscr{G}_{rp} = G$. Then the inductive limit $\mathscr{G} \defeq \ilim \mathscr{G}_p$ can be constructed as $\mathscr{G} \defeq \big( \bigsqcup_p \mathscr{G}_p \big) / { }_{\sim}$, with $\mathscr{G}_p \ni g \sim \sigma_r(g) \in \mathscr{G}_{rp}$ for all $p, r \in S$ and $g \in G$. Now define $\ti{\sigma}_s \colon \mathscr{G} \to \mathscr{G}$ as follows: Given $g \in \mathscr{G}_q$, find $q', s' \in S$ with $q'q = s's$. Such elements exist because $S$ is right reversible. Now define $\ti{\sigma}_s([g]) \defeq [\sigma_{q'}(g)]$, where we view $\sigma_{q'}(g)$ as an element of $\mathscr{G}_{s'}$. It is straightforward to check that $\ti{\sigma}_s$ is an automorphism of $\mathscr{G}$ which satisfies $\tilde{\sigma}_s\vert_G=\sigma_s$ for all $s \in S$, so that $\ti{\sigma}$ extends to the desired algebraic action $\tilde{\sigma} \colon \mathscr{S} \acts \mathscr{G}$. Note that when $G$ is Abelian, the dual action of $\ti{\sigma}$ can be constructed explicitly from $X = \widehat{G}$ and $\alpha = \hat{\sigma}$ (the underlying space will be given as a projective limit obtained from $(X,\alpha)$).
\eex

\bex
Suppose $G\subseteq \Qz^r$ is torsion-free and of finite rank $r\in\Zz_{>0}$. Then $S$ acts by automorphisms on $\Qz\otimes_\Zz G$, so we obtain a globalization by considering the action of the group generated by $S$ on $\Qz\otimes_\Zz G$.
\eex

\bex
\label{ex:genshift}
Let $\Sigma$ be a non-trivial group. The \emph{full $S$-shift over $\Sigma$} is the algebraic $S$-action
\[
\sigma\colon S\acts \textstyle{\bigoplus_S}\Sigma,\quad \sigma_s(x)_t:=\begin{cases}
x_{s^{-1}t} & \text{ if } t\in sS,\\
e & \text{ if } t\notin sS.
\end{cases}
\]
This action admits a globalization if and only if $S$ can be embedded into a group (say $\mathscr{S}$), in which case a globalization is given by the group shift
\[
\tilde{\sigma}\colon\mathscr{S} \acts \textstyle{\bigoplus_{\mathscr{S}}}\Sigma,\quad \tilde{\sigma}_g(x)_h:=x_{g^{-1}h}. 
\]
\eex

\bremark
Suppose that $S \subseteq \mathscr{S}$ is an embedding of $S$ into a group $\mathscr{S}$. Let $\sigma \colon S\acts G$ be an algebraic action with $G$ Abelian. Consider $\Zz \mathscr{S} \otimes_{\Zz S} G$ with the natural $\mathscr{S}$-action and the map $G \to \Zz \mathscr{S} \otimes_{\Zz S} G, \, g \ma 1 \otimes g$. This is the universal enveloping action of $\sigma$ with respect to $S \subseteq \mathscr{S}$, in the sense that any enveloping action of $\sigma$ with respect to $S \subseteq \mathscr{S}$ factors through it. As a consequence, $G \to \Zz \mathscr{S} \otimes_{\Zz S} G$ is injective if and only if there exists a globalization $\tilde{\sigma} \colon \mathscr{S} \acts \mathscr{G}$ of $\sigma$.
\eremark

\subsection{\'Etale groupoids and their C*-algebras}
In this and the next subsection, we introduce some notation for groupoids and their C*-algebras that will be used throughout this paper. For C*-algebras of non-Hausdorff groupoids, we refer the reader to \cite{Connes} or \cite{KS}. For the Hausdorff case, see, e.g., \cite{Renault} and \cite[Part~II]{SSW}.

Let $\cG$ be a (not necessarily Hausdorff) locally compact \'etale groupoid such that the unit space $\cG^{(0)}$ is Hausdorff in the relative topology. We let $r$ and $s$ denote the range and source maps on $\cG$. A subset $B\subseteq \cG$ is said to be a bisection if the restrictions $r\vert_B$ and $s\vert_B$ are injective on $B$.
Given an open bisection $U\subseteq\cG$, we let $C_c(U)$ denote the set of continuous compactly supported complex-valued functions on $U$. Extension-by-zero gives an embedding $C_c(U)\subseteq\ell^\infty(\cG)$, and we let $\cC(\cG)$ be the linear subspace of $\ell^\infty(\cG)$ spanned by the subspaces $C_c(U)$ as $U$ runs through the open bisections of $\cG$. Then $\cC(\cG)$ carries the natural structure of a *-algebra (see, e.g., \cite[\S~3]{Exel08}). The \emph{(full) C*-algebra} of $\cG$, which we denote by $C^*(\cG)$, is the enveloping C*-algebra of $\cC(\cG)$.
For each unit $x\in\cG^{(0)}$, there is a representation $\pi_x\colon \cC(\cG)\to\cB(\ell^2(\cG_x))$, where $\cG_x:=s^{-1}(x)$, such that 
\[
\pi_x(f)\delta_\gamma=\sum_{\alpha\in\cG_{r(\gamma)}}f(\alpha)\delta_{\alpha \gamma} \quad \text{ for all }f\in\cC(\cG).
\]
Here $\delta_{\gamma}$ is given by $\delta_{\gamma}(\gamma') = 1$ if $\gamma' = \gamma$ and $\delta_{\gamma}(\gamma') = 0$ if $\gamma' \neq \gamma$. The \emph{reduced C*-algebra} of $\cG$ is the completion of $\cC(\cG)$ with respect to the norm $||f||_r:=\sup_{x\in\cG^{(0)}}||\pi_x(f)||$. We shall view each $\pi_x$ as a representation of $C_r^*(\cG)$, and let $\pi_r\colon C^*(\cG)\to C_r^*(\cG)$ be the projection map.

The groupoids constructed in this paper will not always be Hausdorff, and we shall see that the essential C*-algebra, as defined in \cite{ExelPitts,KM}, will provide the best model for our concrete C*-algebras. \emph{Essential groupoid crossed products} are defined in \cite[Definition~7.12]{KM} for Fell bundles over $\cG$. Specializing to the case of the trivial Fell bundle whose fibres are all equal to $\Cz$, one arrives at the definition of the \emph{essential C*-algebra} of $\cG$, which is the quotient $C_\es^*(\cG):=C_r^*(\cG)/J_{\sing}$, where $J_{\sing}\subseteq C_r^*(\cG)$ is the ideal of so-called \emph{singular elements} (see \cite[\S~4]{KM}). Several characterizations of $J_{\sing}$ are given in \cite[Proposition~7.18]{KM}, but we shall not need them here. Basic properties of $C_\es^*(\cG)$ are established (in a more general setting) in \cite[\S~7]{KM}. If $\cG$ is Hausdorff, then $C_\es^*(\cG)$ agrees with $C_r^*(\cG)$.

\subsection{Induced representations}
We continue with the setup from the previous subsection. Fix $x\in\cG^{(0)}$. If $\pi$ is a representation of the group C*-algebra $C^*(\cG_x^x)$ on a Hilbert space $H_\pi$, then the associated induced representation $\Ind\pi$ of $C^*(\cG)$ can be explicitly defined as follows (cf. \cite[\S~1]{CN}): Let 
\[
H_{\Ind\pi}:=\Big\{\xi\colon \cG_x\to H_\pi :\xi(gh)=\pi(u_h)^*\xi(g)\text{ for all }g\in \cG_x,h\in\cG_x^x, \text{ and }\sum_{g\in\cG_x/\cG_x^x}||\xi(g)||^2<\infty\Big\},
\]
where $u_g$ is the canonical unitary in $C^*(\cG_x^x)$ corresponding to $g$. Then, $H_{\Ind\pi}$ is a Hilbert space with the obvious linear structure and inner product 
\[
\langle\xi,\eta\rangle_{H_{\Ind\pi}}:=\sum_{g\in \cG_x/\cG_x^x}\xi(g)\overline{\eta(g)}\quad \text{ for }\xi,\eta\in H_{\Ind\pi},
\]
and $\Ind\pi\colon C^*(\cG)\to\cB(\ell^2(\cG_x))$ is defined by
\[
((\Ind\pi)(f)\xi)(g)=\sum_{h\in\cG^{r(g)}}f(h)\xi(h^{-1}g)\quad \text{ for all }f\in\cC(\cG) \text{ and }\xi\in H_{\Ind\pi}.
\]
Let $N\subseteq \cG_x^x$ be a subgroup. The quasi-regular representation $\lambda_{\cG_x^x/N}$ of $C^*(\cG_x^x)$ on $\ell^2(\cG_x^x/N)$ is defined by $\lambda_{\cG_x^x/N}(u_g)\delta_{hN}=\delta_{ghN}$ for all $g,h\in\cG_x^x$. We shall make use of the following observation in several places below.

\bprop
\label{prop:indqr}
There is a unitary $W\colon H_{\Ind\lambda_{\cG_x^x/N}}\cong \ell^2(\cG_x/N)$ such that 
\begin{equation}
\label{e:WindW*}
\left(W(\Ind\lambda_{\cG_x^x/N})(f)W^*\xi\right)(gN)=\sum_{h\in\cG^{r(g)}}f(h)\xi(h^{-1}gN)
\end{equation}
for all $f\in\cC(\cG)$, $\xi\in \ell^2(\cG_x/N)$, and $gN\in \cG_x/N$.
\eprop
\setlength{\parindent}{0cm} \setlength{\parskip}{0cm}

\bproof
For $\xi\in H_{\Ind\pi}$, define $W\xi\colon \cG_x/N\to\Cz$ by $(W\xi)(\gamma N):=\xi(\gamma)(N)$. Let $\cR$ be a complete set of representatives for $\cG_x/\cG_x^x$, so that $\cG_x/N=\bigsqcup_{\gamma\in\cR}\{\gamma hN :hN\in\cG_x^x/N\}$. Then,
\begin{align*}
\sum_{\gamma N\in\cG_x/N}|(W\xi)(\gamma N)|^2&=\sum_{\gamma \in\cR}\sum_{hN\in\cG_x^x/N}|\xi(\gamma h)(N)|^2
=\sum_{\gamma \in\cR}\sum_{hN\in\cG_x^x/N}|(\lambda_{\cG_x^x/N}(h)^*\xi(\gamma))(N)|^2\\
&=\sum_{\gamma \in\cR}\sum_{hN\in\cG_x^x/N}|\xi(\gamma)(hN)|^2=\sum_{\gamma\in \cR}||\xi(\gamma)||_{\ell^2(\cG_x^x/N)}^2=||\xi||_{H_{\Ind\lambda_{\cG_x^x/N}}}^2,
\end{align*}
so $W$ defines an isometry $H_{\Ind\pi}\to \ell^2(\cG_x/N)$. For $\eta\in \ell^2(\cG_x/N)$ and $\gamma\in\cG_x$, define $(V\eta)(\gamma)\colon \cG_x^x/N\to\Cz$ by $[(V\eta)(\gamma)](gN):=\eta(\gamma gN)$. Then $\sum_{gN\in \cG_x^x/N}|\eta(\gamma gN)|^2\leq ||\eta||_{\ell^2(\cG_x/N)}^2$, so $(V\eta)(\gamma)\in\ell^2(\cG_x^x/N)$ and $V$ defines a map $\ell^2(\cG_x/N)\to H_{\Ind\pi}$. For $\eta\in \ell^2(\cG_x/N)$ and $\gamma N\in\cG_x/N$, we have 
$(WV\eta)(\gamma N)=(V\eta)(\gamma)(N)=\eta(\gamma N)$, so that $WV=I$, which shows that $W$ is surjective. It is now easy to verify \eqref{e:WindW*}.
\eproof

\bremark[cf. {\cite[\S~1]{CN}}]
The induced representation $\Ind\lambda_{\cG_x^x}$ coincides with $\pi_x\circ\pi_r$.
\eremark
\setlength{\parindent}{0cm} \setlength{\parskip}{0.5cm}

\section{C*-algebras and groupoids associated with algebraic actions}
\label{sec:C*andgpoids}

\subsection{The concrete C*-algebra associated with an algebraic action} 
Each algebraic action $\sigma \colon S\acts G$ naturally gives rise to a concrete C*-algebra acting on $\ell^2(G)$: Define an isometric representation $\kappa_\sigma\colon S\to\Isom(\ell^2(G))$ by $\kappa_{\sigma}(s) \delta_h =\delta_{\sigma_s(h)}$, where $\{\delta_h : h\in G\}$ is the canonical orthonormal basis for $\ell^2(G)$. We shall call $\kappa_\sigma$ the \emph{Koopman representation} associated with the action $\sigma$.

\bdefin
\label{def:fA}
Suppose $\sigma\colon S \acts G$ is an algebraic action. We let 
\[
  \fA_\sigma \defeq C^*(\{\kappa_\sigma(s): s\in S\}\cup \{\lambda(g) : g\in G\}),
\]
where $\lambda\colon G\to\cU(\ell^2(G))$ is the left regular representation of $G$.
\edefin

In the following, we will simply write $\kappa$ for $\kappa_{\sigma}$ if the algebraic action $\sigma$ is understood.

\bremark
The C*-algebra $\fA_\sigma$ depends only on the image of $\sigma$, so if we are only interested in $\fA_\sigma$, then there is no loss in generality in assuming faithfulness of $\sigma\colon S\acts G$.
\eremark

Let $P \defeq G\rtimes S$ denote the semi-direct product with respect to $\sigma$ taken in the category of monoids. We identify $G$ and $S$ via the embeddings $g \ma (g,1) \in P$ and $s \ma (e,s) \in P$ as submonoids of $P$.

\bremark
Define $\lambda\rtimes\kappa\colon P\to \Isom(\ell^2(G))$ by $\lambda\rtimes\kappa(g,s):=\lambda(g)\kappa(s)$. Then $\fA_\sigma=C_{\lambda\rtimes\kappa}^*(P)$ is the C*-algebra generated by the isometric representation $\lambda\rtimes\kappa$ of $P$. We will see that in general there will not exist a canonical *-homomorphism from the C*-algebra $C^*(P)$ to $C_{\lambda\rtimes\kappa}^*(P)$ (see \S~\ref{sec:semigpC*}), so we cannot appeal to the theory of semigroup C*-algebras to study $\fA_\sigma$. Here, $C^*(P)$ is the full semigroup C*-algebra of $P$ as defined in \cite[Definition~5.6.38.]{CELY}.
\eremark

The following explains our name for $\kappa$:
\bremark
Suppose $G$ is Abelian. Then in terms of the dual action $\alpha \defeq \hat{\sigma} \colon S \acts X = \widehat{G}$, $\kappa$ and $\fA_\sigma$ can be understood as follows: If $\mu$ denotes normalized Haar measure on $X$, then $\alpha$ is measure-preserving in the sense that $\mu(W) = \mu(\alpha_s^{-1}(W))$ for all Borel subsets $W \subseteq X$. Hence we obtain an isometric representation of $S$ on $\Ell^2(X,\mu)$ via $\Ell^2(X,\mu) \to \Ell^2(X,\mu), \, f \ma f \circ \hat{\sigma}_s$ (for $s \in S$), which is the analogue of the Koopman representation in the group case. Moreover, $C(X)$ acts on $\Ell^2(X,\mu)$ by multiplication operators. This yields a unitary representation of $G$ on $\Ell^2(X,\mu)$ via the canonical identification $C^*(G) \cong C(X)$. These two representations of $S$ and $G$ together give rise to a representation of $P$ which is unitarily equivalent to $\lambda\rtimes\kappa$ via the canonical unitary $\Ell^2(X,\mu) \cong \ell^2(G)$.
\eremark

A priori, we have the following description of $\fA_\sigma$:
\[
  \fA_\sigma =\overline{\spn}(\{\kappa(s_1)^* \lambda(g_1) \kappa(t_1) \cdots \kappa(s_m)^* \lambda(g_m) \kappa(t_m) : g_i\in G, s_i,t_i\in S, m\in\Zz_{>0}\}).
\]
We now set out to construct (a candidate for) a groupoid model for $\fA_\sigma$. For this, the language of inverse semigroups is very convenient.

\subsection{The inverse semigroup associated with an algebraic action}
For background on inverse semigroups, see \cite{Law}, and for background on their C*-algebras, see \cite{Pat} and \cite{Exel08,Exel09}.
Let $\cI_G$ denote the inverse semigroup of all partial bijections of $G$. We shall now view $\sigma_s$ as a partial bijection of $G$, so that $\sigma_s^{-1}$ makes sense as an element of $\cI_G$; namely, $\sigma_s^{-1}$ is the partial bijection $\sigma_sG\to G$ given by $\sigma_s(g)\mapsto g$. For each $g\in G$, let $\bbt_g\in\cI_G$ denote the bijection $G \to G$ given by $\bbt_g(h)=gh$ for all $h\in G$.

\bdefin
\label{def:invsemigp}
We let $I_{\sigma}$ denote the inverse sub-semigroup of $\cI_G$ generated by the endomorphisms $\sigma_s$ for $s\in S$ and the translations $\bbt_g$ for $g\in G$.
Explicitly, we have
\[
I_{\sigma}=\{\sigma_{s_1}^{-1}\bbt_{g_1}\sigma_{t_1}\cdots\sigma_{s_m}^{-1}\bbt_{g_m}\sigma_{t_m}: s_i,t_i\in S,\;g_i\in G,\, m\in\Zz_{>0}\}.
\]
\edefin

In the following, we will simply write $I$ for $I_{\sigma}$ if the algebraic action $\sigma$ is understood.

There is a canonical faithful representation by partial isometries $\Lambda \colon \cI_G\hookrightarrow \PIsom(\ell^2(G))$ such that, for $\phi\in\cI_G$ with domain $\dom(\phi)$ and $h\in G$,
\[
\Lambda_\phi\delta_h=\begin{cases}
\delta_{\phi(h)} & \text{ if } h\in\dom(\phi),\\
0 & \text{ if } h\notin\dom(\phi).
\end{cases}
\] 

From now on, we shall use $\Lambda$ to denote the restriction of the above representation to $I$. It is easy to see that this restriction extends the isometric representation $\kappa$, so that $\Lambda\colon I\hookrightarrow \PIsom(\fA_\sigma)$ is a representation of the inverse semigroup $I$ in $\fA_\sigma$. Now it follows immediately that $\fA_\sigma=\overline{\spn}(\{\Lambda_\phi : \phi\in I\})$, so it is reasonable to expect that the structure of $\fA_\sigma$ is closely related to properties of the inverse semigroup $I$. This motivates the analysis of the inverse semigroup $I$. In particular, we wish to compare the C*-algebras associated with $I$ with the C*-algebra $\fA_\sigma$.

In the following, let $\cE_\sigma$ be the idempotent semilattice of $I$. 

\bdefin 
\label{def:C}
We let $\cC_\sigma$ be the smallest family of subgroups of $G$ such that 
\setlength{\parindent}{0cm} \setlength{\parskip}{0cm}

\begin{enumerate}[\upshape(i)]
	\item $G\in\cC_\sigma$;
	\item if $C\in\cC_\sigma$, then $\sigma_sC\in\cC_\sigma$ and $\sigma_s^{-1}C\in\cC_\sigma$ for every $s\in S$.
\end{enumerate}
\edefin
\setlength{\parindent}{0cm} \setlength{\parskip}{0cm}

Here, for a subset $C\subseteq G$, we put $\sigma_s^{-1}C \defeq \{h\in G : \sigma_s(h)\in C\}$. We are writing $\sigma_s^{-1}C$ for the set-theoretic inverse image of $C$ under $\sigma_s$, rather than the  more cumbersome notation $\sigma_s^{-1}(C\cap\sigma_sG)$, which would be used when viewing $\sigma_s$ as a partial bijection.
\setlength{\parindent}{0cm} \setlength{\parskip}{0.5cm}

Members of $\cC_\sigma$ are called \emph{$S$-constructible subgroups}. In the following, we will simply write $\cC$ for $\cC_{\sigma}$ and $\cE$ for $\cE_\sigma$ if the algebraic action $\sigma$ is understood.
Note that
\[
\cC=\{\sigma_{t_1}^{-1}\sigma_{s_1}\cdots\sigma_{t_m}^{-1}\sigma_{s_m}G : s_i,t_i\in S, m\in\Zz_{>0}\}.
\]

\bremark
If $\sigma\colon S\acts G$ is non-automorphic, then $\cC$ is non-trivial, i.e., $\cC\supsetneq \{G\}$. 
\eremark

\bdefin
Let $I^e$ be the inverse sub-semigroup of $I$ generated by $\menge{\sigma_s}{s \in S}$, i.e., 
\[
I^e=\{\sigma_{s_1}^{-1}\sigma_{t_1}\cdots\sigma_{s_m}^{-1}\sigma_{t_m}: s_i,t_i\in S,\, m\in\Zz_{>0} \}.
\]
\edefin
Note that $\cC$ is the semilattice of idempotents of $I^e$.

\bremark
\label{rmk:partialgphom}
The elements in $I^e$ are partial group automorphisms, i.e., they are group isomorphisms from their domains onto their ranges. 
\eremark

\bprop 
\label{prop:C}
The family $\cC$ of $S$-constructible subgroups satisfies the following properties:
\setlength{\parindent}{0cm} \setlength{\parskip}{0cm}

\begin{enumerate}[\upshape(i)]
	\item $\cC$ is closed under taking finite intersections;
	\item if $C\in\cC$, then $\sigma_s^{-1}C$ and $\sigma_sC$ lie in $\cC$;
	\item if $s\in S$, $g\in G$, $C\in\cC$, then $\sigma_s^{-1}(gC)=\emptyset$ or $\sigma_s^{-1}(gC)=h\sigma_s^{-1}C$ for any $h\in \sigma_s^{-1}(gC)$.
\end{enumerate}
\eprop
\setlength{\parindent}{0cm} \setlength{\parskip}{0cm}

\bproof
(i) follows because $\cC$ is the idempotent semilattice of $I^e$, and (ii) follows immediately from the definition of $\cC$. Note that a direct proof of (i) can be given following the proof of \cite[Lemma~3.3]{Li:JFA}.  
\setlength{\parindent}{0.5cm} \setlength{\parskip}{0cm}

(iii): Let $s\in S$, $g\in G$, and $C\in\cC$. If $\sigma_s^{-1}(gC)\neq\emptyset$, then there exists $h\in G$ with $\sigma_s(h)\in gC$. We claim that $\sigma_s^{-1}(gC)=h\sigma_s^{-1}C$. The containment ``$\supseteq$'' is clear. For ``$\subseteq$'', let $k\in \sigma_s^{-1}(gC)$, so that $\sigma_s(k)\in gC$. Now $\sigma_s(h^{-1}k)=\sigma_s(h)^{-1}\sigma_s(k)\in C$, so that $h^{-1}k\in\sigma_s^{-1}C$, i.e., $k\in h\sigma_s^{-1}C$.
\eproof
\setlength{\parindent}{0cm} \setlength{\parskip}{0.5cm}

\bcor
\label{cor:cE=gC}
If $\sigma\colon S\acts G$ is non-automorphic, then $\cE=\{gC : C\in\cC, g\in G\}\cup\{\emptyset\}$.
\ecor
Members of $\cE$ shall be called \emph{$S$-constructible cosets}. We let $\cE^\times:=\cE\setminus\{\emptyset\}$.

Let us now develop a standard form for elements of $I$. 

\bprop
\label{prop:stform}
\begin{enumerate}[\upshape(i)]
\item For $\phi \in I$, we have $\phi \in I^e$ if and only if $\phi(e) = e$.
\item Every $\phi \in I^\times$ is of the form $\phi = \bbt_h \varphi \bbt_{g^{-1}}$ for some $\varphi \in I^e$ and $g,h \in G$. More precisely, if $\dom(\phi) = gC$ for some $g \in G$ and $C \in \cC$, and if $\phi(g) = h$, then $\varphi \defeq \bbt_{h^{-1}} \phi \bbt_g \in I^e$.
\end{enumerate}
\eprop
\setlength{\parindent}{0cm} \setlength{\parskip}{0cm}

\bproof
(i): By definition of $I$, we know that $\phi = \sigma_{s_1}^{-1}\bbt_{g_1}\sigma_{t_1}\cdots\sigma_{s_m}^{-1}\bbt_{g_m}\sigma_{t_m}$. We proceed inductively on $\# \menge{l}{g_l = e}$. Note that $\phi(e) = e$ implies that $e \in \dom(\phi)$, so that $g_m = \bbt_{g_m} \sigma_{t_m}(e) \in \dom(\sigma_{s_m}^{-1})$, i.e., $g_m = \sigma_{s_m}(h_m)$ for some $h_m\in G$. Hence
\begin{align*}
  \phi &= \sigma_{s_1}^{-1}\bbt_{g_1}\sigma_{t_1}\cdots\sigma_{s_m}^{-1}\bbt_{g_m}\sigma_{t_m} = \cdots\sigma_{s_m}^{-1}\bbt_{\sigma_{s_m}(h_m)}\sigma_{t_m}
  = \cdots \bbt_{h_m} \sigma_{s_m}^{-1} \sigma_{t_m} \\
  &= \cdots \sigma_{s_{m-1}}^{-1} \bbt_{g_{m-1}\sigma_{t_{m-1}}(h_m)} \sigma_{t_{m-1}} \sigma_{s_m}^{-1} \sigma_{t_m}
  = \cdots \sigma_{s_{m-1}}^{-1} \bbt_{g_{m-1}\sigma_{t_{m-1}}(h_m)} \sigma_{t_{m-1}} \sigma_{s_m}^{-1} \bbt_e \sigma_{t_m}.
\end{align*}
In this way, we reduce $\# \menge{l}{g_l = e}$, so that we arrive at $\phi = \bbt_g \varphi$ for some $\varphi \in I^e$. But then $\phi(e) = e$ implies that $g = \bbt_g \varphi(e) = \phi(e) = e$. Thus $\phi = \varphi \in I^e$, as desired.
\setlength{\parindent}{0.5cm} \setlength{\parskip}{0cm}

(ii): We know by Corollary~\ref{cor:cE=gC} that $\dom(\phi) = gC$ for some $g \in G$, $C \in \cC$. Let $h \in G$ be such that $\phi(g) = h$. Then $\bbt_{h^{-1}} \phi \bbt_g (e) = e$. (i) implies that $\varphi \defeq \bbt_{h^{-1}} \phi \bbt_g \in I^e$. Hence, $\phi = \bbt_h \varphi \bbt_{g^{-1}}$, as desired.
\eproof
\setlength{\parindent}{0cm} \setlength{\parskip}{0.5cm}

\bcor
\label{cor:0Eunitary}
The inverse semigroup $I$ is 0-E-unitary if and only if $I^e$ is E-unitary.
\ecor
\setlength{\parindent}{0cm} \setlength{\parskip}{0cm}

\bproof
It is easy to see that $I^e$ is E-unitary whenever $I$ is 0-E-unitary. Assume $I^e$ is E-unitary, and suppose $\phi\in I$ and $kD\in \cE^\times$ are such that $kD\subseteq \dom(\phi)$ and $\phi\vert_{kD}=\id_{kD}$. Let $gC=\dom(\phi)$. Note that $kD\subseteq gC$, so that $k\in gC$ and $D\subseteq C$.
By Proposition~\ref{prop:stform}, we can write $\phi=\bbt_h\varphi\bbt_{g^{-1}}$ for $\varphi\in I^e$ with $h=\phi(g)$ and $\dom(\varphi)=C$. Now $\phi\vert_{kD}=\id_{kD}$ is equivalent to
\begin{equation}
\label{eqn:fixxplusD}
    \varphi(g^{-1}kd)=h^{-1}kd\quad\text{ for all }d\in D.
\end{equation}
Taking $d=e$ in \eqref{eqn:fixxplusD} gives $\varphi(g^{-1}k)=h^{-1}k$. Since $g^{-1}k,d\in C=\dom(\varphi)$, we have $\varphi(g^{-1}kd)=\varphi(g^{-1}k)\varphi(d)$ (see Remark~\ref{rmk:partialgphom}), so \eqref{eqn:fixxplusD} reduces to $\varphi(d)=d$ for all $d\in D$. Since $I^e$ is E-unitary, it follows that $\varphi=\id_C$, so that \eqref{eqn:fixxplusD} implies $g=h$. Finally, we see that $\phi=\bbt_g\id_C\bbt_{g^{-1}}=\id_{gC}$.
\eproof
\setlength{\parindent}{0cm} \setlength{\parskip}{0.5cm}

\bremark
\label{rmk:zeros}
If $\sigma\colon S\acts G$ is non-automorphic, then $\emptyset\in\cE$, and we view $\emptyset$ as the distinguished zero element of $\cE$ (in the sense of \cite[Definition~5.5.2]{CELY}). We shall always regard $I^e$ as a semilattice without a distinguished zero element, even though $I^e$ may contain the trivial subgroup $\{e\}$.
\eremark

\subsection{The partial algebraic action associated with a globalization}
\label{ss:PartialGlobal}

For the basics of partial group actions, see, for instance, \cite{ExelPAMS}, \cite{KellLaw}, or \cite[\S~5.5]{CELY}.
If $\tilde{\sigma}\colon\msc{S}\acts \msc{G}$ is a globalization for $\sigma\colon S\acts G$, then the monoid $P=G\rtimes S$ embeds into the group $\Gamma:=\msc{G}\rtimes\msc{S}$. We get an affine action of $\Gamma$ on $\msc{G}$ by $(g,s).x:=g\tilde{\sigma}_s(x)$; we shall often identify $\msc{G}$ and $\msc{S}$ with their images in $\Gamma$. Elements of the group $\gp{P}\subseteq\Gamma$ are then of the form $s_1^{-1}g_1t_1\cdots s_m^{-1}g_mt_m$, where $s_i,t_i\in S$, $g_i\in G$, and $m\in\Zz_{>0}$.

By restricting to the subgroup $G\subseteq \msc{G}$, we obtain a partial affine action $\Gamma\acts G$, where for $\gamma=gs\in \Gamma$ with $g\in\msc{G}$ and $s\in\msc{S}$, $\gamma$ acts as follows: The domain of $\gamma$ is 
\[
G_{\gamma^{-1}}:=\tilde{\sigma}_s^{-1}(g^{-1}G)\cap G,
\]
and the action of $\gamma$ is given by $G_{\gamma^{-1}}\to G_\gamma, \, x \ma \gamma.x \defeq g\tilde{\sigma}_s(x)$. 

\bremark
\label{rmk:idpure}
Since the action of $g \in G$ is given by $\bbt_g$ and the action of $s \in S$ is given by $\sigma_s$, we see that $\bbt_h \sigma_{s_1}^{-1}\sigma_{t_1}\cdots \sigma_{s_m}^{-1} \sigma_{t_m} \bbt_{g^{-1}}\in I$ is a restriction of the partial bijection corresponding to $\gamma= h s_1^{-1}t_1\cdots s_m^{-1}t_m g^{-1} \in \gp{G\rtimes S}$, for $g, h \in G$, $s_i,t_i\in S$.
\setlength{\parindent}{0.5cm} \setlength{\parskip}{0cm}

For $\phi=\bbt_h\varphi\bbt_{g^{-1}}\in I$, with $\varphi=\sigma_{s_1}^{-1}\sigma_{t_1}\cdots \sigma_{s_m}^{-1}\sigma_{t_m}\in I^e$, we have $\dom(\phi)=g\sigma_{t_m}^{-1}\sigma_{s_m}\cdots \sigma_{t_1}^{-1}\sigma_{s_1}G$.
Moreover, for all $x\in\dom(\phi)$, we have $\phi(x)=g_\phi\tilde{\sigma}_{s_{\phi}}(x)$, where $s_\phi=s_1^{-1}t_1\cdots s_m^{-1}t_m\in\gp{S}$ and $g_\phi=h\tilde{\sigma}_{s_{\phi}}(g)^{-1}\in G$.
\eremark

In light of the above remark, it is natural to ask for conditions that will ensure there is a well-defined map $I^\times\to \Gamma$ given by ``$\phi\mapsto g_\phi s_\phi$''.

Suppose $\sigma\colon S\acts G$ has a globalization $\tilde{\sigma}\colon\msc{S}\acts\msc{G}$. Consider the following condition: 
\begin{equation}
	\label{eqn:(JF)}\tag{JF}
	C\subseteq \fix(\tilde{\sigma}_s)\implies s=e,\text{ for all }C\in\cC,\, s\in\gp{S},
\end{equation}
where $\fix(\tilde{\sigma}_s) \defeq \menge{g \in \msc{G}}{ \tilde{\sigma}_s(g) = g}$.

This condition is a kind of joint faithfulness for the partial action of $\gp{S}\subseteq\msc{S}$ on $G$.

\bprop
\label{prop:idpure}
Assume $\sigma\colon S\acts G$ has a globalization $\tilde{\sigma}\colon\msc{S}\acts\msc{G}$. Then $\tilde{\sigma}\colon\msc{S}\acts\msc{G}$ satisfies \eqref{eqn:(JF)} if and only if $I$ is strongly 0-E-unitary.
If these equivalent conditions hold, then there exists an idempotent pure partial homomorphism $\mfg \colon I^\times \to \Gamma$ such that $\mfg(\phi)=g_\phi s_\phi$, where $g_\phi s_\phi\in \Gamma$ is associated with $\phi$ as in Remark~\ref{rmk:idpure}.
\eprop  
\setlength{\parindent}{0cm} \setlength{\parskip}{0cm}

\bproof
Assume \eqref{eqn:(JF)} is satisfied. Let $\phi\in I^\times$. By Proposition~\ref{prop:stform}, we can write $\phi=\bbt_h\varphi\bbt_{g^{-1}}$ for some $g,h\in G$ and $\varphi=\sigma_{s_1}^{-1}\sigma_{t_1}\cdots\sigma_{s_m}^{-1}\sigma_{t_m}$, where $s_i,t_i\in S$ and $m\in\Zz_{>0}$. We have $\dom(\phi)=gC$, where $C=\sigma_{t_m}^{-1}\sigma_{s_m}\cdots \sigma_{t_1}^{-1}\sigma_{s_1}G$, and $\phi(x)=g_\phi\tilde{\sigma}_{s_\phi}(x)$ for all $x\in gC$, where $s_\phi=s_1^{-1}t_1\cdots s_m^{-1}t_m$ and $g_\phi=h\tilde{\sigma}_{s_{\phi}}(g)^{-1}$ (see Remark~\ref{rmk:idpure}). Suppose there exist $k\in G$ and $t\in\gp{S}$ such that 
\begin{equation}
	\label{eqn:sigmaWD}
	g_\phi\tilde{\sigma}_{s_\phi}(x)=k\tilde{\sigma}_t(x)
\end{equation}
for all $x\in gC$. Taking $x=g$ in \eqref{eqn:sigmaWD} gives $g_\phi\tilde{\sigma}_{s_\phi}(g)=k\tilde{\sigma}_t(g)$. Plugging this into \eqref{eqn:sigmaWD} gives $\tilde{\sigma}_{s_\phi}(c)=\tilde{\sigma}_t(c)$ for all $c\in C$, i.e., $\tilde{\sigma}_{t^{-1} s_\phi}(c)=c$ for all $c\in C$. Now $t^{-1}s_\phi=e$ by \eqref{eqn:(JF)}, i.e., $t=s_\phi$, and $k=g_\phi$ follows from \eqref{eqn:sigmaWD}.
This shows that $\mfg \colon I^\times \to \Gamma$ defined by $\mfg(\phi)=g_\phi s_\phi$ is well-defined. An argument similar to the one above shows that $\phi\in\cE^\times$ if and only if $(g_\phi,s_\phi)=(e,1)$, and it is easy to see that $\mfg$ is a partial homomorphism. Hence, $\mfg$ is an idempotent pure partial homomorphism, so that $I$ is strongly 0-E-unitary.
\setlength{\parindent}{0cm} \setlength{\parskip}{0.5cm}

Now assume that $I$ is strongly 0-E-unitary, so that there exists an idempotent pure partial homomorphism $\mfg \colon I^\times\to \Lambda$ for some discrete group $\Lambda$. Observe that the map $G \rtimes S \to I\reg, \, (g,s) \ma \bbt_g \sigma_s$ is injective, so that we may view $G \rtimes S$ as a submonoid in $I\reg$. The restriction of $\mfg$ to the copy of $G\rtimes S$ in $I^\times$ is injective by \cite[Lemma~5.5.7]{CELY}, which gives us an embedding $G\rtimes S\hookrightarrow \Lambda$.
Let $s=s_1^{-1}t_1\cdots s_m^{-1}t_m\in\gp{S}$, and suppose $\tilde{\sigma}_s(x)=x$ for all $x\in C$, where $C\in\cC$. This implies
\begin{equation}
	\label{eqn:phicircid_C}
	\sigma_{s_1}^{-1}\sigma_{t_1}\cdots\sigma_{s_m}^{-1}\sigma_{t_m}\id_C=\id_C
\end{equation}
in $I$ (cf. Remark~\ref{rmk:idpure}). Since both sides are non-zero, applying $\mfg$ to \eqref{eqn:phicircid_C} yields $s=s_1^{-1}t_1\cdots s_m^{-1}t_m=e$. Hence, \eqref{eqn:(JF)} is satisfied.
\eproof
\setlength{\parindent}{0cm} \setlength{\parskip}{0.5cm}

Let us now discuss the special case when $S$ is left Ore and $G$ is Abelian.
\bex
\label{ex:InvSgp-Ore}
Assume $S$ is left Ore and that $G$ is Abelian. We shall write $G$ additively. 
\setlength{\parindent}{0cm} \setlength{\parskip}{0cm}

\begin{enumerate}[\upshape(i)]
\item Let $t,u\in S$. Write $tu^{-1}=\alpha^{-1}\beta$ in $\gp{S}=S^{-1}S$, where $\alpha,\beta\in S$. Then $\sigma_t\sigma_u^{-1}=\sigma_\alpha^{-1}\id_{\sigma_{\alpha t}G}\sigma_\beta$.
\item $I=\{\sigma_s^{-1}\id_{h+D}\bbt_g\sigma_t : g,h\in G,\, s,t\in S,\, D\in\cC\}$.
\item $\sigma\colon S\acts G$ satisfies \eqref{eqn:(JF)} if and only if $C\subseteq \ker(\sigma_s-\sigma_t)\implies s=t$ for all $C\in\cC$ and $s,t\in S$.
\end{enumerate}
\setlength{\parindent}{0cm} \setlength{\parskip}{0.5cm}

For (i), observe that we have $\alpha t=\beta u$, so $\sigma_t\sigma_u^{-1}=\sigma_\alpha^{-1}\sigma_\alpha\sigma_t\sigma_u^{-1}\sigma_\beta^{-1}\sigma_\beta=\sigma_\alpha^{-1}\sigma_{\alpha t}\sigma_{\beta u}^{-1}\sigma_\beta=\sigma_\alpha^{-1}\id_{\sigma_{\alpha t}G}\sigma_\beta$.
\setlength{\parindent}{0.5cm} \setlength{\parskip}{0cm}

For (ii), it suffices to show that $\{\sigma_s^{-1}\id_{h+D}\bbt_g\sigma_t : g,h\in G,\, s,t\in S,\, D\in\cC\}$ is closed, up to $0$, under taking products. Let $s,t,u,v\in S, g,h,k,l\in G$, and $D,E\in \cC$. Since $S$ is left Ore, we can find $\alpha,\beta\in S$ such that $tu^{-1}=\alpha^{-1}\beta$ in $\gp{S}=S^{-1}S$. By (i), we have $\sigma_t\sigma_u^{-1}=\sigma_\alpha^{-1}\id_{\sigma_{\alpha t}G}\sigma_\beta$. 
Now we compute:
\begin{align*}
		\sigma_s^{-1}\id_{h+D}\bbt_g\sigma_t\sigma_u^{-1}\id_{l+E}\bbt_k\sigma_v&=\sigma_s^{-1}\id_{h+D}\bbt_g\sigma_\alpha^{-1}\id_{\sigma_{\alpha t}G}\sigma_\beta\id_{l+E}\bbt_k\sigma_v\\
		&=\sigma_s^{-1}\underbrace{\sigma_\alpha^{-1}\sigma_\alpha}_{=\id_G}\id_{h+D}\sigma_\alpha^{-1}\bbt_{\sigma_\alpha(g)}\id_{\sigma_{\alpha t}G}\sigma_\beta\id_{l+E}\underbrace{\sigma_\beta^{-1}\sigma_\beta}_{=\id_G}\bbt_k\sigma_v\\
		&=\sigma_{\alpha s}^{-1}\id_{\alpha(h+D)}\bbt_{\sigma_\alpha(g)}\id_{\sigma_{\alpha t}G}\id_{\beta(l+E)}\bbt_{\sigma_\beta(k)}\sigma_\beta\sigma_v\\
		&=\sigma_{\alpha s}^{-1}\id_{\alpha(h+D)}\bbt_{\sigma_\alpha(g)}\id_{\sigma_{\alpha t}G\cap \beta(l+E)}\bbt_{\sigma_\alpha(g)}^{-1}\bbt_{\sigma_\alpha(g)}\bbt_{\sigma_\beta(k)}\sigma_{\beta v}\\
		&=\sigma_{\alpha s}^{-1}\id_{\alpha(h+D)}\id_{\sigma_\alpha(g)+\sigma_{\alpha t}G\cap \beta(l+E)}\bbt_{\sigma_\alpha(g)+\sigma_\beta(k)}\sigma_{\beta v}\\
		&=\sigma_{\alpha s}^{-1}\id_{\alpha(h+D)\cap (\sigma_\alpha(g)+\sigma_{\alpha t}G)\cap \beta(l+E)}\bbt_{\sigma_\alpha(g)+\sigma_\beta(k)}\sigma_{\beta v}.
\end{align*}
It remains to observe that $\id_{\alpha(h+D)\cap (\sigma_\alpha(g)+\sigma_{\alpha t}G)\cap \beta(l+E)}$ is zero or of the from $\id_{x+C}$ for some $x \in G$. 

(iii) is true because, for $C\in\cC$ and $s,t\in S$, we have $C\subseteq \ker(\sigma_s-\sigma_t)$ if and only if $C\subseteq \ker(\id-\tilde{\sigma}_{s^{-1}t})$.
\eex

\bremark
\label{rmk:Oreform}
(i) in Example~\ref{ex:InvSgp-Ore} is equivalent to the statement that if $S$ is left Ore and $G$ is Abelian, then 
$\fA_\sigma=\overline{\spn}(\{\kappa(s)^*1_{hD}\lambda(g) \kappa(t) : s,t\in S, g,h\in G, D\in \cC\})$.
\eremark
\setlength{\parindent}{0cm} \setlength{\parskip}{0.5cm}

\subsection{The unit space of the groupoid model}
\label{ss:UnitSpace}

Let us assume in this subsection that our action is non-automorphic, so that $\cE=\{gC : C\in\cC, g\in G\}\cup\{\emptyset\}$ is the idempotent semilattice of our inverse semigroup $I$. As before, we write $\cE\reg \defeq \cE \setminus \gekl{\emptyset}$. Define $\widehat{\cE}$ as the space of characters of $\cE$, i.e., non-zero multiplicative maps $\cE \to \gekl{0,1}$ sending $\emptyset$ to $0$, equipped with the topology of point-wise convergence. A basis of open sets is given by 
\[
 \widehat{\cE}(gC; \gekl{g_iC_i}) \defeq \menge{\chi \in \widehat{\cE}}{\chi(gC) = 1; \chi(g_iC_i) = 0 \ \forall \ i},
\]
where $gC \in \cE\reg$ and $\gekl{g_iC_i} \subseteq \cE\reg$ is a finite subset. Without loss of generality we may assume $g_iC_i \subseteq gC$ for all $i$. There is a one-to-one correspondence between characters of $\cE$ and filters on $\cE$, i.e., subsets $\cF \subseteq \cE$ with the following properties: $\emptyset \notin \cF$; $G \in \cF$; if $gC \in \cF$ and $hD \in \cE$ with $gC \subseteq hD$, then $hD \in \cF$; and if $gC,hD \in \cF$, then $gC\cap hD \in \cF$. This one-to-one correspondence is implemented by the assignment $\widehat{\cE} \ni \chi \ma \cF(\chi) \defeq \menge{gC \in \cE}{\chi(gC) = 1}$.

\bdefin
Let $\widehat{\cE}_{\max}$ denote the characters $\chi$ of $\cE$ for which $\cF(\chi)$ is maximal with respect to inclusion.
\edefin
\setlength{\parindent}{0cm} \setlength{\parskip}{0cm}

In other words, $\chi \in \widehat{\cE}$ belongs to $\widehat{\cE}_{\max}$ if and only if we cannot find $\chi' \in \widehat{\cE}$ with $\cF(\chi) \subsetneq \cF(\chi')$.
\setlength{\parindent}{0cm} \setlength{\parskip}{0.5cm}

\bdefin
The boundary of $\widehat{\cE}$ is given by $\partial \widehat{\cE} \defeq \overline{\widehat{\cE}_{\max}}$.
\edefin
\setlength{\parindent}{0cm} \setlength{\parskip}{0cm}

Following \cite{Exel08,Exel09}, characters of $\cE$ which belong to $\partial \widehat{\cE}$ are called \emph{tight}, and we also call the corresponding filters tight. We briefly recall several notions from \cite{Exel08,Exel09,ExelPardo}. 
A \emph{cover} (resp. \emph{outer cover}) of a subset $\cF\subseteq \cE$ is a subset $\mfc \subseteq \cF$ (resp. $\mfc \subseteq \cE$) such that for each $gC\in\cF^\times:=\cF\cap \cE^\times$, there exists $hD\in \mfc$ with $hD \cap gC \neq \emptyset$. 
For $\phi\in I$, let $\fix(\phi):=\{g\in\dom(\phi) : \phi(g)=g\}$, and let $\cJ_\phi:=\{hD\in\cE : hD\subseteq\fix(\phi)\}\cup\{\emptyset\}$. 
A subset $\mfc\subseteq \cE$ is a \emph{cover} of the constructible coset $gC$ if $\mfc$ is a cover of $\cJ_{gC}$, where we view $gC$ as an element of $I$ using the identification of $\cE$ with the idempotent semilattice of $I$.
It is shown in \cite{Exel08,Exel09} that $\chi \in \widehat{\cE}$ belongs to $\partial \widehat{\cE}$ if and only if for every $gC \in \cE\reg$ with $\chi(gC) = 1$ and every cover $\mfc$ of $gC$, there exists $hD \in \mfc$ such that $\chi(hD) = 1$. To ease notation, let $\cup\cJ_\phi:=\bigcup_{gC\in\cJ_\phi}gC$.
\setlength{\parindent}{0cm} \setlength{\parskip}{0.5cm}

\blemma
\label{lem:covers=unions}
Let $\phi\in I^\times$. A finite collection $\c=\{h_iD_i\}\subseteq\cE^\times$ is a (finite) outer cover for $\cJ_\phi$ if and only if $\cup \cJ_\phi\subseteq \bigcup_ih_iD_i$. In particular, $\c$ is a cover of the constructible coset $gC$ if and only if $gC=\bigcup_ih_iD_i$.
\elemma
\setlength{\parindent}{0cm} \setlength{\parskip}{0cm}

\bproof
The ``if'' direction is clear, so suppose $\mfc = \gekl{h_iD_i}$ is a finite outer cover for $\cJ_\phi$. Set $D \defeq \bigcap_i D_i$. We have $D \in \cC$ because $\cC$ is closed under finite intersections. Let $gC\in\cJ_\phi$. For every $k \in gC$, we have $k(C\cap D) \subseteq gC$, so that, since $\mfc$ is a cover for $\cJ_\phi$, there exists $i$ such that $h_iD_i\cap k(C\cap D) \neq \emptyset$. Since $k(C\cap D) \subseteq kD_i$, this implies that $kD_i \cap h_iD_i \neq \emptyset$, i.e., $k \in h_iD_i$. Hence, $gC \subseteq \bigcup_i h_iD_i$, so that $\cup \cJ_\phi\subseteq\bigcup_ih_iD_i$. If $\c$ is a cover of $\cJ_\phi$, then $\mfc\subseteq\cJ_\phi$, so that $\bigcup_ih_iD_i\subseteq  \cup\cJ_\phi$. Since $\cup\cJ_{gC}=gC$, we see that $\c$ is a cover of the constructible coset $gC$ if and only if $gC=\bigcup_ih_iD_i$.
\eproof
\setlength{\parindent}{0cm} \setlength{\parskip}{0.5cm}

Let us now give a characterization of tight characters in our situation.
\blemma
\label{lem:CharTightChar}
Let $\chi \in \widehat{\cE}$. Then $\chi$ lies in $\partial \widehat{\cE}$ if and only if for all $gC \in \cE\reg$ with $\chi(gC) = 1$ and all $D \in \cC$ with $D \subseteq C$, $[C:D] < \infty$, and $C = \bigsqcup_i k_iD$ for some $k_i \in G$, there exists $i$ such that $\chi(gk_iD) = 1$.
\elemma
\setlength{\parindent}{0cm} \setlength{\parskip}{0cm}

\bproof
\an{$\Rarr$} follows from the characterization of tight characters in \cite{Exel08,Exel09} mentioned above because $\gekl{gk_iD}$ is a cover of $gC$.
For \an{$\Larr$}, suppose that $\gekl{h_jD_j}$ is a cover of $gC$. By \cite[(4.4)]{Neu}, we may without loss of generality assume that $[C:D_j] < \infty$ for all $j$. Set $D \defeq \bigcap_j D_j$. Then $[C:D] < \infty$, so that $C = \bigsqcup_i k_iD$ for some $k_i \in G$. By assumption, there exists $i$ such that $\chi(gk_iD) = 1$. Moreover, as $gC = \bigcup_j h_jD_j$ by Lemma~\ref{lem:covers=unions}, there exists $j$ such that $gk_iD \subseteq h_jD_j$. Hence $\chi(h_jD_j) = 1$, as desired.
\eproof
\setlength{\parindent}{0cm} \setlength{\parskip}{0.5cm}

Let us now relate characters and filters on $\cE$ with those on $\cC$. As above, a filter on $\cC$ is a subset $\fF \subseteq \cC$ with the following properties: $G \in \fF$, if $C \in \fF$ and $D \in \cC$ with $C \subseteq D$, then $D \in \fF$, and if $C,D \in \fF$, then $C\cap D \in \fF$. 

\bdefin
We call a filter $\fF$ on $\cC$ \emph{finitely hereditary} if whenever $C \in \fF$ and $D \in \cC$ satisfy $D \subseteq C$ with $[C:D] < \infty$, then $D \in \fF$. 
\edefin

Given $\chi\in\widehat{\cE}$, define $\Pi(\chi)\colon\cC\to\{0,1\}$ by
\[
\Pi(\chi)(C):=\begin{cases}
	1 & \text{if there exists $g\in G$ with }\chi(gC)=1,\\
	0 & \text{if }\chi(gC)=0\text{ for all }g\in G.
\end{cases}
\]
Moreover, set $\fF(\chi) \defeq \fF(\Pi(\chi)) \defeq \menge{C \in \cC}{\Pi(\chi)(C) = 1}$.
\blemma
\label{lem:filtersOnEvsOnC}
\begin{enumerate}[\upshape(i)]
\item For every $\chi \in \widehat{\cE}$, $\fF(\chi)$ is a filter on $\cC$.
\item $\Pi$ is a surjective map from $\widehat{\cE}$ onto the space of characters of $\cC$. More precisely, for each filter $\fF$ on $\cC$ there exists a character $\chi(\fF) \in \widehat{\cE}$ uniquely determined by $\cF(\chi(\fF)) = \fF$. Moreover, for every filter $\fF$ on $\cC$, $\chi \in \widehat{\cE}$ lies in $\Pi^{-1}(\fF)$ if and only if $\chi(g'C') = 0$ for all $C' \in \cC$ with $C' \notin \fF$, and for all $C \in \fF$, there exists $g_C \in G$ such that $\chi(g_CC) = 1$, and $(g_C)$ satisfies the compatibility condition that $g_C \in g_DD$ for all $C, D \in \fF_{\cC}$ with $C \subseteq D$. 
\item Given $\chi \in \widehat{\cE}$, we have $\chi \in \partial \widehat{\cE}$ if and only if $\fF(\chi)$ is a finitely hereditary filter.
\item The maximal filter on $\cC$ is given by $\fF_{\max} = \cC$.
\item $\Pi^{-1}(\fF_{\max})$ is a subset of $\widehat{\cE}_{\max}$ which is dense in $\partial \widehat{\cE}$.
\end{enumerate}
\elemma
\setlength{\parindent}{0cm} \setlength{\parskip}{0cm}

\bproof
(i) is true by construction. For surjectivity of $\Pi$ in (ii), just observe that, given a filter $\fF$ on $\cC$, $\chi(\fF)(C) \defeq 1$ if $C \in \fF$ and $\chi(\fF)(gC) \defeq 0$ for all $gC \in \cE$ with $gC \notin \fF$ defines a character of $\cE$ which satisfies $\cF(\chi(\fF)) = \fF$. The remaining claims in (ii) are easy to see.
\setlength{\parindent}{0.5cm} \setlength{\parskip}{0cm}

(iii) follows from Lemma~\ref{lem:CharTightChar}. (iv) is true because any two elements of $\cC$ always have non-empty intersection, as they all contain $e \in G$. It remains to show (v). $\Pi^{-1}(\fF_{\max}) \subseteq \widehat{\cE}_{\max}$ is clear. Given a non-empty basic open set $\partial \widehat{\cE}(gC;\gekl{h_iD_i})$, there exists $k \in gC \setminus \bigcup_i h_iD_i$. Hence $kD_i \cap h_iD_i = \emptyset$ for all $i$. Moreover, by (ii), there exists $\chi \in \Pi^{-1}(\fF_{\max})$ with the property that $\chi(kC) = 1$ for all $C \in \fF_{\max}$. It follows that $\chi \in \partial \widehat{\cE}(gC;\gekl{h_iD_i})$, as desired.
\eproof
\setlength{\parindent}{0cm} \setlength{\parskip}{0.5cm}

\bremark
The compatibility condition in Lemma~\ref{lem:filtersOnEvsOnC}~(ii) is equivalent to the condition that $(g_C)$ is an element of the projective limit $\plim_{C \in \fF} \gekl{G/C}$.
\eremark

The following notation will be convenient.
\bdefin
\label{def:chi_z}
Given $k \in G$, we denote by $\chi_k$ the character of $\cE$ which satisfies $\chi(kC) = 1$ for all $C \in \cC$.
\edefin
\setlength{\parindent}{0cm} \setlength{\parskip}{0cm}

Note that $\chi_k$ exists by Lemma~\ref{lem:filtersOnEvsOnC}~(ii). Indeed, we have $\chi_k = k. \chi(\fF_{\max})$.
\setlength{\parindent}{0cm} \setlength{\parskip}{0.5cm}

\bremark
\label{rmk:denseset}
It follows from Lemma~\ref{lem:filtersOnEvsOnC}~(v) that $\{\chi_k : k\in G\}$ is dense in $\partial\E$. Thus, $\bd$ is a completion of (a quotient of) $G$.
\eremark

\subsection{The groupoid associated with an algebraic action}
\label{ss:GroupoidModel}

We present (a candidate for) a groupoid model of $\fA_\sigma$. 
Let $I*\widehat{\cE}:=\{(\phi,\chi)\in I\times\widehat{\cE} : \chi(\dom(\phi))=1\}$, and define an equivalence relation on $I*\widehat{\cE}$ by
\[
(\phi,\chi)\sim (\psi,\omega) \text{ if }\chi=\omega \text{ and there exists }gC\in \cE^\times\text{ with } \phi\vert_{gC}=\psi\vert_{gC}\text{ and }\chi(gC)=1.
\]
We let $[\phi,\chi]$ denote the equivalence class of $(\phi,\chi)$. For each $gC\in \cE$, let $\widehat{\cE}(gC):=\{\chi\in\widehat{\cE} : \chi(gC)=1\}$. For $gC\in \cE$ and a finite subset $\f\subseteq\cE$, put
\[
\widehat{\cE}(gC;\f):=\widehat{\cE}(gC)\setminus\bigcup_{hD\in\f}\widehat{\cE}(hD).
\]

For $\chi\in \widehat{\cE}(\dom(\phi))$, define $\phi.\chi:=\chi(\phi^{-1}\sqcup \phi)\in \widehat{\cE}(\im(\phi))$.

\bdefin
Let 
\[
I\ltimes\widehat{\cE}:=I*\widehat{\cE}/ { }_{\sim}=\{[\phi,\chi] : (\phi,\chi)\in I*\cE\}
\] 
be the \emph{transformation groupoid} attached to $I\acts\widehat{\cE}$ with range and source maps given by $\mathrm{r}([\phi,\chi])=\phi.\chi$ and $\mathrm{s}([\phi,\chi])=\chi$, respectively. Multiplication and inversion are given by $[\phi,\psi.\chi][\psi,\chi]=[\phi\psi,\chi]$ and $[\phi,\chi]^{-1}=[\phi^{-1},\phi.\chi]$, respectively.
\edefin
\setlength{\parindent}{0cm} \setlength{\parskip}{0cm}

Now $I\ltimes\widehat{\cE}$ becomes an ample groupoid when equipped with the topology with basis of compact open bisections
\[
[\phi,\widehat{\cE}(gC;\f)]:=\{[\phi,\chi] : \chi\in \widehat{\cE}(gC;\f)\} \quad \text{for }gC\in\cE, \f\subseteq\cE \text{ finite with } \widehat{\cE}(gC;\f)\subseteq\widehat{\cE}(\dom(\phi)).
\]
\setlength{\parindent}{0cm} \setlength{\parskip}{0.5cm}

The subspace $\partial\E$ is $I$-invariant by \cite[Proposition~12.8.]{Exel08}, so that we get an action $I\acts \partial\widehat{\cE}$. The basic open subsets of $\bd$ are of the form 
\[
\partial \widehat{\cE}(gC;\gekl{h_iD_i}) \defeq \partial \widehat{\cE} \cap \widehat{\cE}(gC;\gekl{h_iD_i}).
\]

\bdefin 
We set $\cG_\sigma:=I\ltimes\partial\widehat{\cE} = (I \ltimes \widehat{\cE}) \vert_{\partial\widehat{\cE}}^{\partial\widehat{\cE}} = \menge{[\phi,\chi] \in I \ltimes \widehat{\cE}}{\chi, \, \phi.\chi \in \partial\widehat{\cE}}$.
\edefin

In the terminology from \cite{Exel08}, $\cG_\sigma$ is the \emph{tight groupoid} of $I$. When $I$ is countable, the C*-algebra $C^*(\cG_\sigma)$ is universal for \emph{tight representations} of the inverse semigroup $I$ by \cite[Theorem~13.3]{Exel08}. 

As $\Lambda\colon I\to \fA_\sigma$ is a representation, the universal property of $C^*(I)$ yields a (surjective) *-homomorphism $\tilde{\rho}\colon C^*(I\ltimes\widehat{\cE})\cong C^*(I)\to \fA_\sigma$ such that $\tilde{\rho}(1_{[\phi,\widehat{\cE}(\dom(\phi))]})=\Lambda_\phi$ for all $\phi\in I^\times$.

\bprop
\label{prop:Lambdatight}
Assume $I$ is countable and that $\sigma\colon S\acts G$ is non-automorphic. Then the representation $\Lambda\colon I\to \fA_\sigma$ is tight, so that $\tilde{\rho}\colon C^*(I\ltimes\widehat{\cE})\to \fA_\sigma$ factors through a representation $\rho\colon C^*(\cG_\sigma)\to \fA_\sigma$. Moreover, the restriction of $\rho$ to $C(\partial\E)$ is an isomorphism onto $\fD_\sigma:=\overline{\spn}(\{1_{gC} : g\in G,C\in\cC\})$.
\eprop
\setlength{\parindent}{0cm} \setlength{\parskip}{0cm}

\bproof
Since $\rho$ is unital, \cite[Corollary~4.3]{Exel21} implies that $\Lambda$ is tight if and only if it is cover-to-join in the sense of \cite[\S~3]{Exel21}, i.e., for every $gC\in \cE^\times$, we have
\begin{equation}
	\label{eqn:Lambdatight}
	\Lambda(gC)\leq \bigvee_{hD\in \mfc}\Lambda(hD)
\end{equation}
for all finite covers $\mfc$ of $gC$. Now suppose $\mfc=\{g_iC_i : i\in F\}$ is a finite cover of $gC$. By Lemma~\ref{lem:covers=unions}, we have $gC=\bigcup_{i\in F}g_iC_i$. In $\fD_\sigma$, we have 
\[
\bigvee_{i\in F}1_{g_iC_i}=\sum_{\emptyset\neq I\subseteq F}(-1)^{\# I-1}\prod_{i\in I}1_{g_iC_i}=1_{\cup_{i\in F}g_iC_i}=1_{gC}.
\]
Hence, $\bigvee_{i\in F}\Lambda(g_iC_i)=\Lambda(gC)$.
\setlength{\parindent}{0cm} \setlength{\parskip}{0.5cm}

We now turn to the second claim. It is clear that the restriction of $\rho$ to $C(\partial\E)$ has image equal to $\fD_\sigma$. It follows from \cite[Proposition 5.6.21.]{CELY} that there is an inverse *-homomorphism, which implies injectivity.
\eproof
\setlength{\parindent}{0cm} \setlength{\parskip}{0.5cm}

\bremark
The C*-algebra $C^*(\cG_\sigma)$ provides a universal model for $\fA_\sigma$.
\setlength{\parindent}{0.5cm} \setlength{\parskip}{0cm}

The C*-algebras $C^*(I)=C^*(I\ltimes\widehat{\cE})$ and $C_r^*(I)=C_r^*(I\ltimes\widehat{\cE})$ can be regarded as the full and reduced ``Toeplitz-type'' C*-algebras associated with $\sigma\colon S\acts G$.
\eremark
\setlength{\parindent}{0cm} \setlength{\parskip}{0.5cm}

\bremark
\label{rem:TrafoGPD}
Let $\sigma\colon S\acts G$ be an algebraic action that admits a globalization $\tilde{\sigma}\colon\msc{S}\acts\msc{G}$ as in \S~\ref{ss:PartialGlobal}. Assume \eqref{eqn:(JF)} is satisfied, and let $\mfg\colon I^\times\to \Gamma=\msc{G}\rtimes\msc{S}$ be the idempotent pure partial homomorphism from Proposition~\ref{prop:idpure}.
General results for inverse semigroups admitting idempotent pure partial homomorphisms to groups now give us a partial action of $\Gamma=\msc{G}\rtimes\msc{S}$ on $\widehat{\cE}$ (cf. \cite{Li:IMRN}, \cite[Chapter~5]{CELY}, or \cite{Li:Ktheory}): For each $\gamma\in \Gamma$, let 
\[
U_{\gamma^{-1}}:=\{\chi\in\widehat{\cE} : \chi(gC)=1,\text{ where }gC=\dom(\phi)\text{ for some }\phi\in I^\times\text{ with }\mfg(\phi)=\gamma\}.
\]
Then $\gamma \in \Gamma$ acts via the homeomorphism $U_{\gamma^{-1}}\to U_{\gamma}, \, \chi \ma \gamma.\chi$ defined by
\[
 (\gamma.\chi)(hD):=\chi(\gamma^{-1}.((hD) \cap \dom(\phi)))
\]
for any $\phi\in I^\times$ with $\mfg(\phi)=\gamma$. We have isomorphisms $I\ltimes\widehat{\cE}\cong \Gamma\ltimes\widehat{\cE}$ and $I\ltimes\partial\widehat{\cE}\cong \Gamma\ltimes \partial\widehat{\cE}$.
\eremark

\section{Properties of the groupoid model}
\label{sec:properties}

In this section, we study properties of the groupoid $\cG_\sigma$ attached to the non-automorphic algebraic action $\sigma\colon S\acts G$. These groupoid properties then translate into properties of the C*-algebras $C_r^*(\cG_\sigma)$ and $C_\es^*(\cG_\sigma)$.

\subsection{Hausdorffness}
\label{ss:Hd}

Combining Lemma~\ref{lem:covers=unions} with \cite[Theorem~3.16]{ExelPardo}, we obtain:
\bprop
\label{prop:(H)}
The groupoid $\cG_\sigma$ is Hausdorff if and only if the following condition is satisfied:
\begin{equation}
\label{eqn:(H)}\tag{H}
    \text{For all }\phi\in I, \text{ there exist constructible cosets }g_1C_1,...,g_nC_n\subseteq \fix(\phi)\text{ such that }\cup\cJ_\phi=\textstyle{\bigcup_i} g_iC_i.
\end{equation}
\eprop
\setlength{\parindent}{0cm} \setlength{\parskip}{0cm}

Recall that $\cJ_\phi=\{hD\in\cE : hD\subseteq\fix(\phi)\}\cup\{\emptyset\}$.
\setlength{\parindent}{0cm} \setlength{\parskip}{0.5cm}

\bremark
Clearly, if $\fix(\phi)$ is finite for all $\phi \in I\setminus \cE$, then \eqref{eqn:(H)} is satisfied and $\cG_\sigma$ is Hausdorff.

Assume that $\sigma\colon S\acts G$ has a globalization $\tilde{\sigma}\colon\mathscr{S}\acts\mathscr{G}$. 
\setlength{\parindent}{0.5cm} \setlength{\parskip}{0cm}

If \eqref{eqn:(JF)} is satisfied, then $\cG_\sigma \cong (\msc{G}\rtimes\mathscr{S}) \ltimes \bd$ is a partial transformation groupoid and hence Hausdorff. 

If $\mathscr{S}$ is torsion-free, $G$ is Abelian, and the dual action $\hat{\sigma}$ is mixing, then ${\fix}(\ti{\sigma}_s) = \gekl{e}$ for all $1 \neq s \in \mathscr{S}$ (see, e.g., \cite[Remark~4.1]{BruLi:rigidity}). Hence ${\fix}(\ti{\sigma}_s)$ is finite for all $1 \neq s \in \mathscr{S}$ and thus $\cG_\sigma$ is Hausdorff by the observation above. If we know in addition that $\gekl{e} \notin \cC$, then $I$ is strongly 0-E-unitary. 
\eremark
\setlength{\parindent}{0cm} \setlength{\parskip}{0.5cm}

\subsection{Closed invariant subspaces of the boundary} 

Recall that we introduced the notion of filters on $\cE$ in \S~\ref{ss:UnitSpace}, and that there is a one-to-one correspondence between characters of $\cE$ (i.e., elements of $\E$) and filters given by $\E \ni \chi \ma \cF(\chi) \defeq \menge{gC \in \cE}{\chi(gC) = 1}$. This bijection restricts to a one-to-one correspondence between tight characters and tight filters. Moreover, we defined the map $\Pi$ from characters of $\cE$ to characters of $\cC$ and set $\fF(\chi) \defeq \fF(\Pi(\chi)) \defeq \menge{C \in \cC}{\Pi(\chi)(C) = 1}$ for every $\chi \in \E$. We observed that for all $\chi \in \E$, $\chi$ is tight if and only if $\fF(\chi)$ is finitely hereditary. Moreover, for every filter $\fF$ on $\cC$ there exists a uniquely determined character $\chi(\fF) \in \E$ with $\cF(\chi(\fF)) = \fF$. In particular, for a finitely hereditary filter $\fF$ on $\cC$, $\chi(\fF)$ is a  tight character.

In order to describe closed invariant subspaces of $\bd$, we need some terminology.
\bdefin
Let $\bm{\fF}_{\cC}$ be the set of finitely hereditary filters on $\cC$. We define a partial action of $I^e$ on $\bm{\fF}_{\cC}$ as follows: Given $\fF \in \bm{\fF}_{\cC}$ and $\varphi \in I^e$, $\varphi.\fF$ is defined if there exists $C_{\varphi} \in \fF$ such that $C_{\varphi} \subseteq \dom(\varphi)$. In that case, $\varphi.\fF$ is defined as the smallest element of $\bm{\fF}_{\cC}$ containing $\menge{\varphi(C_{\varphi} \cap C)}{C \in \fF} = \menge{\varphi(D)}{D \in \fF, \, D \subseteq C_{\varphi}}$. A subset $\bm{\fF} \subseteq \bm{\fF}_{\cC}$ is called $I^e$-invariant if for all $\fF \in \bm{\fF}$ and $\varphi \in I^e$ such that $\varphi.\fF$ is defined, we have $\varphi.\fF \in \bm{\fF}$.
\setlength{\parindent}{0.5cm} \setlength{\parskip}{0cm}

A subset $\bm{\fF} \subseteq \bm{\fF}_{\cC}$ is called $\subseteq$-closed if for all $\fE \in \bm{\fF}_\cC$, $\fE \subseteq \bigcup_{\fF \in \bm{\fF}} \fF$ implies $\fE \in \bm{\fF}$.
\edefin
\setlength{\parindent}{0cm} \setlength{\parskip}{0.5cm}

Our main result concerning closed invariant subspaces of $\bd$ reads as follows:
\btheo
\label{thm:ClosedInvariantSubsp}
There is a one-to-one correspondence between closed invariant subspaces of $\bd$ and $I^e$-invariant, $\subseteq$-closed subsets of $\bm{\fF}_{\cC}$ sending $\fX \subseteq \bd$ to $\fF(\fX) \defeq \menge{\fF(\chi)}{\chi \in \fX}$. The inverse map sends $\bm{\fF} \subseteq \bm{\fF}_{\cC}$ to $\fF^{-1}(\bm{\fF}) \defeq \menge{\chi \in \bd}{\fF(\chi) \in \bm{\fF}}$.
\etheo
For the proof, we first show that closed, $G$-invariant subspaces of $\bd$ are in one-to-one correspondence with $\subseteq$-closed subsets of $\bm{\fF}_{\cC}$. 
\blemma
\label{lem:Pi-1Pi}
Given $\fX \subseteq \bd$, we have $\Pi^{-1}(\Pi(\fX)) \subseteq \overline{\menge{\bbt_k.\chi}{k \in G, \, \chi \in \fX}}$.
\elemma
\setlength{\parindent}{0cm} \setlength{\parskip}{0cm}

\bproof
Suppose that $\chi, \omega \in \bd$ satisfy $\Pi(\omega) = \Pi(\chi)$. Then we claim that $\omega \in \overline{\menge{\bbt_k.\chi}{k \in G}}$. Indeed, take a basic open neighbourhood $\bd(gC;\gekl{h_iD_i})$ of $\omega$. Then $\omega(gC) = 1$ and $\omega(h_iD_i) = 0$ for all $i$. Without loss of generality, we can assume that $D_i \notin \fF(\omega)$ (otherwise replace $C$ by $C \cap D_i$). As $\Pi(\omega) = \Pi(\chi)$, this implies that there exists $k \in G$ with $\chi(kC) = 1$ while $\bbt_h.\chi(h_iD_i) = 0$ for all $i$ and all $h \in G$. Therefore, $\bbt_{gk^{-1}}.\chi(gC) = \chi(kC) = 1$ whereas $\bbt_{gk^{-1}}.\chi(h_iD_i) = 0$ for all $i$. Hence $\bbt_{gk^{-1}}.\chi(gC) \in \bd(gC;\gekl{h_iD_i})$, as desired.
\eproof
\setlength{\parindent}{0cm} \setlength{\parskip}{0.5cm}

\bprop
\label{prop:A-ClosureVSunion-closed}
Given $\fX \subseteq \bd$ and $\omega \in \bd$, we have $\omega \in \overline{\menge{\bbt_k.\chi}{k \in G, \, \chi \in \fX}}$ if and only if $\fF(\omega) \subseteq \bigcup_{\chi \in \fX} \fF(\chi)$.
\eprop
\setlength{\parindent}{0cm} \setlength{\parskip}{0cm}

\bproof
\an{$\Rarr$}: We have $C \in \fF(\omega)$ if and only if $\omega(gC) = 1$ for some $g \in G$. If $\omega \in \overline{\menge{\bbt_k.\chi}{k \in G, \, \chi \in \fX}}$, then there exists $\chi \in \fX$ and $k \in G$ such that $\bbt_k.\chi(gC) = 1$. Hence $C \in \fF(\bbt_k.\chi) = \fF(\chi)$.
\setlength{\parindent}{0cm} \setlength{\parskip}{0.5cm}

\an{$\Larr$}: Without loss of generality, we may assume that $\fF(\omega) = \cF(\omega)$ because Lemma~\ref{lem:Pi-1Pi} allows us to replace $\omega$ by $\chi(\fF(\omega))$ if necessary. Take a basic open neighbourhood $\bd(gC;\gekl{h_iD_i})$ of $\omega$, with $h_iD_i \subseteq gC$ for all $i$. We may assume $g=e$. Hence $\omega(C) = 1$ and $\omega(h_iD_i) = 0$ for all $i$. We may also assume that $D_i \notin \fF(\omega)$ for all $i$ (otherwise replace $C$ by $C \cap D_i$). Thus $[C:D_i] = \infty$ for all $i$. Since $\fF(\omega) \subseteq \bigcup_{\chi \in \fX} \fF(\chi)$, we can find $\chi \in \fX$ with $\chi(lC) = 1$ for some $l\in G$. Without loss of generality assume that $l=e$, i.e., $\chi(C) = 1$. For each $i$ with $D_i \in \fF(\chi)$ choose $k_i \in G$ so that $\chi(k_iD_i) = 1$. 
For each $i$, we have $h_iD_ik_i^{-1}=D_i'h_ik_i^{-1}$ where $D_i':=h_iD_ih_i^{-1}$. Since $h_iD_i\subseteq C$, we have $D_i'\subseteq C$. Moreover, $[C:D_i']=[C:D_i]=\infty$ for all $i$, so by \cite[Lemma~4.1]{Neu} we have that $C \not\subseteq \bigcup_i D_i'h_ik_i^{-1}$. Thus there exists $h \in C$ with $h \notin h_iD_ik_i^{-1}$ for all $i$.
Therefore $\bbt_h.\chi(C) = 1$ and $\bbt_h.\chi(hk_iD_i) = 1$ for all $i$. But $h \notin h_iD_ik_i^{-1}$ implies $hk_iD_i\neq h_iD_i$ for all $i$. Hence $(hk_iD_i)\cap (h_iD_i)  = \emptyset$ for all $i$, and thus $\bbt_h.\chi(h_iD_i) = 0$ for all $i$. We conclude that $\bbt_h.\chi \in \bd(C;\gekl{h_iD_i})$, as desired.
\eproof
\setlength{\parindent}{0cm} \setlength{\parskip}{0.5cm}

\bcor
\label{cor:ClAinv=union-cl}
The map in Theorem~\ref{thm:ClosedInvariantSubsp} implements a one-to-one correspondence between closed, $G$-invariant subspaces of $\bd$ and $\subseteq$-closed subsets of $\bm{\fF}_{\cC}$.
\ecor
\setlength{\parindent}{0cm} \setlength{\parskip}{0cm}

\bproof
Proposition~\ref{prop:A-ClosureVSunion-closed} implies that for any closed, $G$-invariant subspace $\fX$ of $\bd$, $\fF(\fX)$ is $\subseteq$-closed. Moreover, given a $\subseteq$-closed subset $\bm{\fF}$ of $\bm{\fF}_{\cC}$, $\fF^{-1}(\bm{\fF})$ is clearly $G$-invariant. To see that it is also closed, take $\omega \in \overline{\fF^{-1}(\bm{\fF})}$. Proposition~\ref{prop:A-ClosureVSunion-closed} implies that $\fF(\omega) \subseteq \bigcup_{\fF \in \bm{\fF}} \fF$. But since $\bm{\fF}$ is $\subseteq$-closed, this implies $\fF(\omega) \in \bm{\fF}$ and thus $\omega \in \fF^{-1}(\bm{\fF})$, as desired.
\setlength{\parindent}{0cm} \setlength{\parskip}{0.5cm}

To see that these maps are inverse to each other, first note that $\fF(\fF^{-1}(\bm{\fF})) = \bm{\fF}$ because $\fF$ is surjective. Moreover, to show $\fF^{-1}(\fF(\fX)) = \fX$, it suffices to show \an{$\subseteq$} as \an{$\supseteq$} is clear. So take $\omega \in \bd$ with $\fF(\omega) \in \fF(\fX)$. Then $\fF(\omega) = \fF(\chi)$ for some $\chi \in \fX$. It follows that $\Pi(\omega) = \Pi(\chi)$ and thus $\omega \in \overline{\menge{\bbt_k.\chi}{k \in G}}$ by Lemma~\ref{lem:Pi-1Pi}. This shows \an{$\subseteq$}.
\eproof
\setlength{\parindent}{0cm} \setlength{\parskip}{0.5cm}

For the proof of Theorem~\ref{thm:ClosedInvariantSubsp}, it remains to show that the one-to-one correspondence in Corollary~\ref{cor:ClAinv=union-cl} restricts to a one-to-one correspondence between closed invariant subspaces of $\bd$ and $I^e$-invariant, $\subseteq$-closed subsets of $\bm{\fF}_{\cC}$. In other words, we have to show that the $I^e$-action is preserved. This is a consequence of the following observations:
\blemma
\begin{enumerate}[\upshape(i)]
\item For all $\varphi \in I^e$ and $\chi \in \bd$, if $\varphi.\chi$ is defined, then $\varphi.\fF(\chi)$ is defined.
\item For all $\varphi \in I^e$ and $\chi \in \bd$, if $\varphi.\fF(\chi)$ is defined, then there exists $g \in G$ such that $\varphi.(\bbt_g.\chi)$ is defined.
\item For all $\varphi \in I^e$ and $\chi \in \bd$ such that $\varphi.\chi$ and $\varphi.\fF(\chi)$ are defined, we have $\fF(\varphi.\chi) = \varphi.\fF(\chi)$.
\end{enumerate}
\elemma
\setlength{\parindent}{0cm} \setlength{\parskip}{0cm}

\bproof
(i) and (ii) follow from the observations that $\varphi.\chi$ is defined if and only if there exist $g \in G$ and $C \in \cC$ with $\chi(gC) = 1$ and $gC \subseteq \dom(\varphi)$, while the latter is equivalent to $g \in \dom(\varphi)$ and $C \subseteq \dom(\varphi)$ because $\varphi \in I^e$, and that $\varphi.\fF(\chi)$ is defined if and only if there exist $C \in \fF(\chi)$ with $C \subseteq \dom(\varphi)$.
\setlength{\parindent}{0.5cm} \setlength{\parskip}{0cm}

For (iii), assume that $g_\varphi C_{\varphi} \in \cE^\times$ satisfies $g_\varphi C_{\varphi} \subseteq \dom(\varphi)$ and $\chi(g_\varphi C_\varphi)=1$. Observe that $D \in \fF(\varphi.\chi)$ if and only if there exists $h \in G$ such that $\varphi.\chi(hD) = 1$, which is equivalent to existence of $h \in G$ such that $\chi(\varphi^{-1}((\varphi(g_{\varphi}) \varphi(C_{\varphi})) \cap (hD))) = 1$. Furthermore, observe that $D \in \varphi.\fF(\chi)$ if and only if there exists $C \in \fF(\chi)$ such that $\varphi(C_{\varphi} \cap C) \subseteq D$. Therefore, we need to show that there exists $h \in G$ such that $\chi(\varphi^{-1}((\varphi(g_{\varphi}) \varphi(C_{\varphi})) \cap (hD))) = 1$ if and only if there exists $C \in \fF(\chi)$ such that $\varphi(C_{\varphi} \cap C) \subseteq D$. For \an{$\Rarr$}, set $C \defeq \varphi^{-1}(\varphi(C_{\varphi}) \cap D)$. Then $C \in \fF(\chi)$, and $\varphi(C_{\varphi} \cap C) = \varphi(C_{\varphi}) \cap D \subseteq D$. For \an{$\Larr$}, $C \in \fF(\chi)$ implies that there exists $g \in G$ such that $\chi(gC) = 1$. Hence there exists $k \in G$ such that $k (C_{\varphi} \cap C) = (g_{\varphi} C_{\varphi}) \cap (gC) \in \cF(\chi)$. So $h \defeq \varphi(k)$ satisfies $h \varphi(C_{\varphi} \cap C) \subseteq (\varphi(g_{\varphi})\varphi(C_{\varphi})) \cap (hD)$ and $\chi(\varphi^{-1}(h\varphi(C_{\varphi} \cap C))) = \chi(kC_{\varphi} \cap C) = 1$.
\eproof
\setlength{\parindent}{0cm} \setlength{\parskip}{0.5cm}

\bremark
Theorem~\ref{thm:ClosedInvariantSubsp} reduces the problem of computing all closed invariant subspaces of $\bd$ to the study of certain subsets of $\bm{\fF}_{\cC}$, and the relevant subsets are singled out by two conditions involving the $I^e$-action and set-theoretical properties, but no topology. The question remains whether it is possible to compute closed invariant subspaces more concretely. In principle, $I^e$-invariant and $\subseteq$-closed subsets $\bm{\fF}$ of $\bm{\fF}_{\cC}$ are completely determined by the subset $\bigcup_{\fF \in \bm{\fF}} \fF$ of $\cC$. However, in general, it seems to be a challenge to characterize which subsets of $\cC$ arise in this way.
\eremark

\subsection{Minimality}

Let us now give several characterizations of minimality of our groupoid model. 

\btheo
\label{thm:minimal}
Let $\sigma\colon S\acts G$ be a non-automorphic algebraic action. The following conditions are each equivalent to minimality of $\cG_\sigma$:
\setlength{\parindent}{0cm} \setlength{\parskip}{0cm}

\begin{enumerate}[label=\upshape (M\theenumi),ref= M\theenumi]
    \item\label{eqn:(M)} For all $C,D\in\cC$, there exist $\phi_1,...,\phi_n\in I$ such that $C\subseteq\textstyle{\bigcup_{i=1}^n}\phi_i(\dom(\phi_i) \cap D)$.
    \item\label{eqn:(M2)} For all $D\in\cC$, there exist $\phi_1, \dotsc, \phi_n \in I$ with $G= \bigcup_{i=1}^n \phi_i(\dom(\phi_i) \cap D)$.
    \item\label{eqn:(M3)} For all $C,D\in \cC$, there exist $\varphi\in I^e$ and  $D'\in\cC$ with $D'\subseteq D\cap \dom(\varphi)$ such that $\varphi(D')\subseteq C$ and $[C:\varphi(D')]<\infty$.
    \item\label{eqn:(M4)} For all $D\in \cC$, there exist $\varphi\in I^e$ and  $D'\in\cC$ with $D'\subseteq D\cap \dom(\varphi)$ such that $[G : \varphi(D')]<\infty$.   
    \item\label{eqn:(M5)} $\bm{\fF}_{\cC}$ contains no proper, $I^e$-invariant, $\subseteq$-closed subsets.
\end{enumerate}
\etheo
\setlength{\parindent}{0cm} \setlength{\parskip}{0cm}

\bproof
By Lemma~\ref{lem:covers=unions} and \cite[Theorem~5.5]{ExelPardo}, minimality of $\cG_\sigma$ implies condition \eqref{eqn:(M)}. The implications  \eqref{eqn:(M)}$\Rightarrow$\eqref{eqn:(M2)} and \eqref{eqn:(M3)}$\Rightarrow$\eqref{eqn:(M4)} are clear, and that \eqref{eqn:(M5)} implies minimality of $\cG_\sigma$ follows from Theorem~\ref{thm:ClosedInvariantSubsp}.
\setlength{\parindent}{0cm} \setlength{\parskip}{0.5cm}

\eqref{eqn:(M2)}$\Rightarrow$\eqref{eqn:(M3)}: Let $C,D\in\cC$. By assumption, there exist $\phi_1,\dotsc,\phi_n \in I$ with $G=\bigcup_{i=1}^n \phi_i(\dom(\phi_i) \cap D)$, where $\dom(\phi_i) \cap D\neq\emptyset$ for all $i$. We have 
\[
C=\textstyle{\bigcup_{i=1}^n} C\cap \phi_i(\dom(\phi_i) \cap D)=\textstyle{\bigcup_{i=1}^n}\phi_i\left(\phi_i^{-1}(C\cap\im(\phi_i))\cap D\right).
\]
Without loss in generality, we may assume that, for every $i$, $\phi_i\left(\phi_i^{-1}(C\cap\im(\phi_i))\cap D\right)$ is non-empty, so that there exists a constructible coset $k_iD_i\subseteq\dom(\phi_i)\cap D$ with $\phi_i^{-1}(C\cap\im(\phi_i))\cap D=k_iD_i$.
By Proposition~\ref{prop:stform}, we can write $\phi_i=\bbt_{h_i}\varphi_i\bbt_{g_i^{-1}}$ for some $g_i,h_i\in G$ and $\varphi_i\in I^e$ such that $\dom(\varphi_i)\in\cC$ and $\dom(\phi_i)=g_i\dom(\varphi_i)$. Since $k_iD_i\subseteq \dom(\phi_i)$, we have $k_i=g_ic_i$ for some $c_i\in \dom(\varphi_i)$ and $D_i\subseteq \dom(\varphi_i)$, so that $c_iD_i\subseteq \dom(\varphi_i)$. By Remark~\ref{rmk:partialgphom}, $\varphi_i$ is a homomorphism on $\dom(\varphi_i)$, so that
\[
\phi_i(k_iD_i)=h_i\varphi_i(g_i^{-1}k_iD_i)=h_i\varphi_i(c_iD_i)=h_i\varphi_i(c_i)\varphi_i(D_i).
\]
Now $C=\bigcup_{i=1}^nh_i\varphi_i(c_i)\varphi_i(D_i)$, so by \cite[(4.3)]{Neu}, we may assume that $\varphi_i(D_i)$ has finite index in $C$ for all $i$. Since $e\in C$, there exists $i$ such that $e\in h_i\varphi_i(c_i)\varphi_i(D_i)$, i.e., $h_i\varphi_i(c_i)\varphi_i(D_i)=\varphi_i(D_i)$. It remains to observe that $D_i\subseteq D$ because $k_iD_i\subseteq D$.

\eqref{eqn:(M4)}$\Rightarrow$\eqref{eqn:(M5)}: Suppose $\bm{\fF} \subseteq \bm{\fF}_{\cC}$ is a non-empty, $I^e$-invariant, $\subseteq$-closed subset. Let $\fF\in \bm{\fF}$. Given $D\in\cC$, there exists a constructible subgroup $D'\subseteq D$ and $\varphi\in I^e$ such that $D'\subseteq\dom(\varphi)$ and $\varphi(D')\subseteq G$ has finite index. Since $G\in \fF$ and $\fF$ is finitely hereditary, it follows that $\varphi(D')\in\fF$. Thus, $\varphi^{-1}.\fF$ is defined and contains $D'$, so that $D\in \varphi^{-1}.\fF$ because $\varphi^{-1}.\fF$ is a filter. Moreover, $\varphi^{-1}.\fF\in \bm{\fF}$ by $I^e$-invariance. Since $D$ was arbitrary, it follows that $\cC\subseteq \bigcup_{\fF\in\bm{\fF}}\fF$. Since $\bm{\fF}$ is $\subseteq$-closed, it follows that $\fF\in \bm{\fF}$ for all $\fF\in \bm{\fF}_{\cC}$.
\eproof
\setlength{\parindent}{0cm} \setlength{\parskip}{0.5cm}

\subsection{Topological freeness}

Following \cite[\S~2.4]{KM}, we call an \'etale groupoid $\cG$ \emph{topologically free} if for every open bisection $U$ of $\cG$, $U\subseteq \Iso(\cG)\setminus\cG^{(0)}$ implies that $\{x\in\cG^{(0)} : \cG_x^x\cap U\neq\emptyset\}$ has empty interior.
In general, $\cG$ is topologically free whenever it is effective (i.e., $\Iso(\cG)^\circ=\cG^{(0)}$, where $\Iso(\cG):=\bigcup_{x\in\cG^{(0)}}\cG_x^x$ is the isotropy bundle of $\cG$). The converse holds if $\cG$ is Hausdorff. Topological freeness is important because of its relationship to the intersection property, see \cite[\S~7]{KM} and \cite[\S~7.1]{KKLRU}.
 
We now turn to topological freeness for the groupoids $I\ltimes \fX$, where $\fX\subseteq\bd$ is a closed, invariant subset.

\bdefin
\label{def:exact}
The algebraic action $\sigma \colon S \acts G$ is called \emph{exact} if $\bigcap_{C \in \cC} C = \gekl{e}$.
\edefin

\bremark
\label{rem:Intersection}
The group $\bigcap_{C \in \cC} C$ is the biggest subgroup $G_c$ of $G$ which is invariant under $\sigma_s$ for all $s \in S$ such that $\sigma_s \vert_{G_c}$ is surjective for all $s \in S$. Thus, $\sigma\colon S\acts G_c$ is an automorphic $S$-action, and exactness of $\sigma \colon S \acts G$ is equivalent to saying that there are no $S$-invariant subgroups $H\subseteq G$ such that the associated $S$-action $S\acts H$ is automorphic.
\eremark

\bremark
Our definition of exactness for an algebraic action is a vast generalization of the notion of exactness for a single endomorphism given by Rohlin in \cite{Rohlin}.

When $\sigma\colon S\acts G$ is an algebraic dynamical system in the sense of \cite{BLS2}, then it is exact if and only if $\bigcap_{s\in S}\sigma_sG=\{e\}$. This stronger condition is called ``minimal'' in \cite{Sta}. Our condition \eqref{eqn:(M)} is automatically satisfied for the actions in \cite{Sta,BLS2} (see \S~\ref{sec:leftreversible} below), which could explain their choice of terminology (cf. \cite[Remark~1.7]{Sta}). 
\eremark

For $\fF\in\bm{\fF}$, put $\cap\fF:=\bigcap_{C\in\fF}C$, and for $H\subseteq G$, let $\core(H):=\bigcap_{g\in G}gHg^{-1}$. Our main result on topological freeness reads as follows:
\btheo
\label{thm:TopFree}
Suppose that $\fX \subseteq \bd$ is a closed invariant subspace, and let $\bm{\fF} = \fF(\fX)$.
\setlength{\parindent}{0cm} \setlength{\parskip}{0cm}

\begin{enumerate}[\upshape(i)]
\item If the groupoid $I \ltimes \fX$ obtained by restricting $\cG_\sigma$ to $\fX$ is topologically free, then for all $D \in \bigcup_{\fF \in \bm{\fF}} \fF$, we have $\bigcap_{\fF \in \bm{\fF}, \, D \in \fF} \core(\cap\fF) = \gekl{e}$.
\item If for all $D \in \bigcup_{\fF \in \bm{\fF}} \fF$, we have $\bigcap_{\fF \in \bm{\fF}, \, D \in \fF}\cap\fF = \gekl{e}$, then $I \ltimes \fX$ is topologically free.
\end{enumerate}
In particular, if $G$ is Abelian, then $\cG_\sigma$ is topologically free if and only if $\sigma\colon S\acts G$ is exact.
\etheo
\setlength{\parindent}{0cm} \setlength{\parskip}{0.5cm}

For the proof, we need some preparations.

\blemma
\label{lem:bd()Inbd(.,.)}
Let $gC,g_iC_i\in\cE^\times$ for $1\leq i\leq n$ with $gC\setminus \bigcup_ig_iC_i\neq\emptyset$. Then there exists $hD\in\cE^\times$ such that $hD\subseteq gC\setminus\bigcup_ig_iC_i$. Moreover, if $\bigcup_ig_iC_i\subsetneq gC$, then for any constructible coset $hD\subseteq gC\setminus\bigcup_ig_iC_i$, we have  $\partial\widehat{\cE}(hD)\subseteq\partial\widehat{\cE}(gC;\{g_iC_i\})$.
\elemma
\setlength{\parindent}{0cm} \setlength{\parskip}{0cm}

\bproof
Put $D:=C\cap \bigcap_iC_i$. Since $D\subseteq C_i$ for all $i$ and $D\subseteq C$, each of the cosets $gC$ and $g_iC_i$ can be written as a non-trivial (possibly infinite) disjoint union of $D$-cosets. Hence, $gC\setminus\bigcup_ig_iC_i$ contains a $D$-coset, $hD$ say.
\setlength{\parindent}{0cm} \setlength{\parskip}{0.5cm}

Now suppose $\bigcup_ig_iC_i\subsetneq gC$. If $\chi\in \bd(hD)$, i.e., $\chi(hD)=1$, then $\chi(gC)=1$ because $hD\subseteq gC$. If we had $\chi(g_iC_i)=1$ for some $i$, then we would have $\chi(hD\cap g_iC_i)=1$; this is impossible since $hD$ and $g_iC_i$ are disjoint. Thus,  $\partial\widehat{\cE}(hD)\subseteq\partial\widehat{\cE}(gC;\{g_iC_i\})$.
\eproof
\setlength{\parindent}{0cm} \setlength{\parskip}{0.5cm}

\blemma
\label{lem:TopFree}
Let $\fX$ and $\bm{\fF}$ be as in Theorem~\ref{thm:TopFree}. The following are true:
\setlength{\parindent}{0cm} \setlength{\parskip}{0cm}

\begin{enumerate}[\upshape(i)]
\item $\bigcup_{\fF \in \bm{\fF}} \menge{\bbt_k.\chi(\fF)}{k \in G}$ is dense in $\fX$. In particular, $\menge{\chi_k}{k \in G}$ is dense in $\bd$.
\item Given $\fF \in \bm{\fF}$ and $\phi \in I$, if $\phi.\chi(\fF) = \chi(\fF)$, then $\phi(\cap\fF) = \cap\fF$.
\end{enumerate}
\elemma
\setlength{\parindent}{0cm} \setlength{\parskip}{0cm}

\bproof
(i) follows from Proposition~\ref{prop:A-ClosureVSunion-closed}. To see (ii), observe that for $D \subseteq \dom(\phi)$, we have $\chi(\fF)(D) = 1$ if and only if $\chi(\fF)(\phi(D)) = 1$. This implies $\phi(\bigcap_{C \in \fF} C) \supseteq \bigcap_{C \in \fF} C$. The reverse inclusion follows by replacing $\phi$ by $\phi^{-1}$.
\eproof
\setlength{\parindent}{0cm} \setlength{\parskip}{0.5cm}

\bproof[Proof of Theorem~\ref{thm:TopFree}]
(i): Suppose that there exists $D \in \bigcup_{\fF \in \bm{\fF}} \fF$ with $\bigcap_{\fF \in \bm{\fF}, \, D \in \fF} \core(\cap\fF) \neq \gekl{e}$. Take $e \neq k \in \bigcap_{\fF \in \bm{\fF}, \, D \in \fF}\core(\cap\fF)$. For all $\chi \in \fX$ with $D \in \fF(\chi)$, we have $k^{-1}gC=gC$ for every $g\in G$ and $C\in \fF(\chi)$, which implies $\bbt_k.\chi = \chi$. Hence, $[\gekl{\bbt_k} \times (\fX \cap \bd(D))] \subseteq \Iso(I \ltimes \fX) \setminus \fX$. Since $\fX \cap \bd(D)\neq\emptyset$, this implies that $I \ltimes \fX$ is not topologically free.

(ii): Suppose that $[\phi,U] \subseteq \Iso(I \ltimes \fX)$ for some $\phi \in I$ and non-empty open set $U \subseteq \fX \cap \partial\widehat{\cE}(\dom(\phi))$. Then $U$ contains a basic open set of the form $\fX \cap \partial\widehat{\cE}(gC;\gekl{h_iD_i})$, and by Lemma~\ref{lem:bd()Inbd(.,.)}, there exists $hD\in\cE^\times$ such that $\partial\widehat{\cE}(hD)\subseteq\partial\widehat{\cE}(gC;\gekl{h_iD_i})$. Hence $\phi.\chi = \chi$ for all $\chi \in \fX$ with $\chi(hD) = 1$. In particular, for all $\fF \in \bm{\fF}$ with $D \in \fF$ and $k \in hD$, we have $\phi.(\bbt_k(\chi(\fF))) = \bbt_k(\chi(\fF))$. It follows that $(\bbt_{k^{-1}} \phi \bbt_k).\chi(\fF) = \chi(\fF)$ and thus $(\bbt_{k^{-1}} \phi \bbt_k)(\cap\fF) = \cap\fF$ by Lemma~\ref{lem:TopFree}, i.e., $\phi(k(\cap\fF)) =  k (\cap\fF)$. 
By assumption, $\bigcap_{\fF \in \bm{\fF}, \, D \in \fF} \cap\fF = \gekl{e}$, which implies $\phi(k) = k$ for all $k \in hD$, i.e., $\phi \vert_{hD} = \id_{hD}$. Hence $\fX \cap \partial\widehat{\cE}(hD) = [\gekl{\phi} \times (\fX \cap \partial\widehat{\cE}(hD))] \subseteq [\phi,U] \cap \fX$ and thus $[\phi,U] \cap \fX \neq \emptyset$, as desired.

The last claim follows from parts (i) and (ii).
\eproof

Since $\cG_\sigma$ is ample, we derive the following immediate consequence from \cite[Corollary~3.15]{BoLi}.
\bcor
If the conditions in part (ii) of Theorem~\ref{thm:TopFree} are satisfied for all closed, invariant subspaces $\fX\subseteq\bd$ and $\cG_\sigma$ is Hausdorff and inner exact, then $C^*_r(\cG_\sigma)$ has the ideal property.
\ecor

We also deduce the following consequence, which is special to our situation and in general only holds for amenable groupoids. As in Remark~\ref{rem:Intersection}, set $G_c:=\bigcap_{C \in \cC}C$.

\bcor
\label{cor:separatesideals}
Consider the following statements:
\setlength{\parindent}{0cm} \setlength{\parskip}{0cm}

\begin{enumerate}[\upshape(i)]
\item $\sigma\colon S\acts G$ is exact.
\item $\cG_\sigma$ is topologically free.
\item $C(\bd) \subseteq C_{\es}^*(\cG_\sigma)$ has the ideal intersection property.
\end{enumerate}
Then $(i)\Rarr (ii)\Rarr(iii)$. If $G_c$ is amenable and $G_c=\bigcap_{C \in \cC} \core(C)$ (e.g., if $G$ is Abelian), then $(iii)\Rarr(i)$. 
\ecor
\setlength{\parindent}{0cm} \setlength{\parskip}{0cm}

\bproof
The implication ``$(i)\Rarr (ii)$'' follows from part (ii) of Theorem~\ref{thm:TopFree}, and ``$(ii)\Rarr(iii)$'' follows from \cite[Theorem~7.29]{KM}.
\setlength{\parindent}{0.5cm} \setlength{\parskip}{0cm}

``$(iii)\Rarr(i)$'': Assume $G_c$ is amenable and $G_c=\bigcap_{C \in \cC} \core(C)$, and 
suppose that $G_c \neq \gekl{e}$. Arguing as in the proof of part (i) of Theorem~\ref{thm:TopFree}, we see that $[G_c,\chi]:=\{[\bbt_g,\chi] : g\in G_c]\}$ is contained in $\Iso(\cG_\sigma)$ and that $[G_c,\chi]\cap\bd =\{[e,\chi]\}$. Remark~\ref{rem:Intersection} implies that $\menge{(\chi,[G_c,\chi])}{\chi \in \bd}$ is an essentially confined amenable section of isotropy groups of $\cG_\sigma$, in the sense of \cite[Definition~7.1]{KKLRU}. Hence \cite[Theorem~7.2]{KKLRU} implies that $C(\bd) \subseteq C_{\es}^*(\cG_\sigma)$ does not have the ideal intersection property.
\eproof
\setlength{\parindent}{0cm} \setlength{\parskip}{0.5cm}

\subsection{Amenability}

We now consider amenability for our groupoids.

\btheo
\label{thm:Amenability}
Suppose $\sigma\colon S\acts A$ has a globalization $\tilde{\sigma}\colon\mathscr{S}\acts\mathscr{G}$ and that \eqref{eqn:(JF)} is satisfied. Assume (without loss of generality) that $\mathscr{S}$ is generated by $S$. If $\cG_\sigma$ is amenable, then $\mathscr{S}$ is amenable. The converse holds if $\msc{G}$ is amenable.
\etheo
\bproof
As explained in \S~\ref{ss:GroupoidModel}, our assumptions mean that we have the identification $\cG_\sigma \cong (\msc{G}\rtimes\mathscr{S}) \ltimes \bd$. 
Our assumptions also give $\mathscr{S} \subseteq (\cG_\sigma)^{\chi_e}_{\chi_e}$, so amenability of $\mathscr{S}$  follows from amenability of $\cG_\sigma$ by \cite[Proposition~5.1.1]{AR}. If $\msc{G}$ and $\mathscr{S}$ are amenable, then $\msc{G}\rtimes\mathscr{S}$ is amenable, so that $(\msc{G}\rtimes\mathscr{S}) \ltimes \bd$ is amenable.
\eproof

\bremark
In general, $C_r^*(\cG_\sigma)$ is nuclear if and only if $\cG_\sigma$ is amenable (see \cite{AR}).
Assume we are in the setting of Theorem~\ref{thm:Amenability}, so that $\cG_\sigma\cong (\msc{G}\rtimes\mathscr{S}) \ltimes \bd$. If the group $\msc{G}\rtimes\mathscr{S}$ is exact, then $C_r^*(\cG_\sigma)$ is nuclear if and only if the canonical map $C^*(\cG_\sigma)\to C_r^*(\cG_\sigma)$ is an isomorphism by \cite[Theorem~4.12]{BFS}. 
\eremark

\subsection{Pure infiniteness}

We now turn to pure infiniteness following \cite[\S~4]{Mat15}. Let $\cG$ be an ample \'{e}tale groupoid with compact unit space. A subset $\fX\subseteq \cG^{(0)}$ is said to be \emph{properly infinite} if there exist compact open bisections $U,V\subseteq \cG$ such that $s(U)=s(V)=\fX$ and $r(U)\sqcup r(V)\subseteq \fX$. For many classes of \'{e}tale groupoids, one has a dichotomy (see, e.g., \cite{BoLi,RS,Ma}): Either the groupoid is purely infinite, or the groupoid admits an invariant measure. Indeed, we have the following open question: 
\begin{question}
\label{ques:pi}
Let $\cG$ be a second countable, topologically free, minimal, ample \'{e}tale groupoid with compact unit space. If there is no invariant measure on $\cG^{(0)}$, then must $\cG$ be purely infinite?
\end{question}

For our groupoid $\cG_\sigma$, if it is minimal, then the action $\sigma\colon S\acts G$ must be non-automorphic, from which it is easy to see that $\bd$ is properly infinite, so that $\bd$ has no invariant measures. Thus, Theorem~\ref{thm:PI} below answers Question~\ref{ques:pi} in the affirmative for the class of groupoids arising from algebraic actions.

\btheo
\label{thm:PI}
If $\cG_\sigma$ is minimal, then $\cG_\sigma$ is purely infinite.
\etheo

Combining this with our results above and general results for groupoid C*-algebras, we obtain:

\bcor
\label{cor:PI}
Suppose $\sigma\colon S\acts G$ is a non-automorphic algebraic action. If $\sigma\colon S\acts G$ is exact and satisfies one of the conditions in Theorem~\ref{thm:minimal}, then $C_{\es}^*(\cG_\sigma)$ is simple and purely infinite. 
\ecor
\setlength{\parindent}{0cm} \setlength{\parskip}{0cm}

\bproof
By Theorem~\ref{thm:TopFree}, $\cG_\sigma$ is topologically free, and by Theorem~\ref{thm:minimal}, $\cG_\sigma$ is minimal. The groupoid $\cG_\sigma$ is purely infinite by Theorem~\ref{thm:PI}, which implies $\cG_\sigma$ is locally contacting. Now the result follows from \cite[Theorem~7.26]{KM} (see \cite[Remark~7.27]{KM}).
\eproof
\setlength{\parindent}{0cm} \setlength{\parskip}{0.5cm}

In addition to the application to C*-algebras above, pure infiniteness has implications for topological full groups and homology: It was recently proven in \cite{Li:TFG} that purely infinite groupoids satisfy Matui's AH Conjecture from \cite{Mat12}, and general results for topological full groups of purely infinite groupoids have been established in \cite{GT}.

\bremark
Condition \eqref{eqn:(M)} implies that $\gekl{e} \notin \cC$. If $\{e\}\in\cC$, then $\fA_\sigma\cong C_{\es}^*(\cG_\sigma)$ (see Corollary~\ref{cor:rho<Indchi} below), and $\fA_\sigma$ contains the compact operators, so $C_{\es}^*(\cG_\sigma)$ cannot be purely infinite in this case.
\eremark

Before proceeding to the proof of Theorem~\ref{thm:PI}, we need two lemmas. 
\blemma
\label{lem:fibasics}
Let $B\in\cC$, and let $\{k_iB_i\}\subseteq\cE^\times$ be a finite (possibly empty) collection. Then for all $gC\in\cE^\times$ with $gC\subseteq B$ and $1<[B:C]<\infty$, we have 
\[
\bd(B;\{k_iB_i\}\cup\{gC\})=\bigsqcup_{hC\in B/C,hC\neq gC}\bd(B\cap hC;\{k_iB_i\}).
\]
\elemma
\setlength{\parindent}{0cm} \setlength{\parskip}{0cm}

\bproof
``$\subseteq$'': Suppose $\chi\in \bd(B;\{k_iB_i\}\cup\{gC\})$. Then, $\chi(B)=1$ and $\chi(gC)=\chi(k_iB_i)=0$ for all $i$. We have the finite decomposition $B=gC\sqcup\bigsqcup_{hC\in B/C,hC\neq gC}hC$, so $\chi(hC)=1$ for some $hC\in B/C$ by Lemma~\ref{lem:CharTightChar}. Since $\chi(gC)=0$, $hC\neq gC$, so that $\chi\in \bd(B\cap hC;\{k_iB_i\})$.
\setlength{\parindent}{0.5cm} \setlength{\parskip}{0cm}

``$\supseteq$'': Suppose $\chi\in \bd(B\cap hC;\{k_iB_i\})$ for some $hC\in B/C$ with $hC\neq gC$. Then $\chi(k_iB_i)=0$ for all $i$ and $1=\chi(B\cap hC)=\chi(B)\chi(hC)$, so that $\chi(B)=\chi(hC)=1$. Since $gC\cap hC=\emptyset$, $\chi(gC)=0$, so $\chi \in\bd(B;\{k_iB_i\}\cup\{gC\})$.
\eproof
\setlength{\parindent}{0cm} \setlength{\parskip}{0.5cm}

\blemma
\label{lem:PI}
Assume $\cG_\sigma$ is minimal. If $B\in\cC$ satisfies $[B:C]<\infty$ for all constructible subgroups $C\subseteq B$, then $\sigma\colon S\acts G$ satisfies \eqref{eqn:(FI)}.
\elemma
\setlength{\parindent}{0cm} \setlength{\parskip}{0cm}

\bproof
By Theorem~\ref{thm:minimal}, there exist a constructible subgroup $C\subseteq B$ and $\varphi\in I^e$ such that $C\subseteq\dom(\varphi)$ and $\varphi(C)\subseteq G$ has finite index. Since the partial isomorphism $\varphi^{-1}$ maps $\varphi(C)\cap B$ onto $C':=C\cap \varphi^{-1}(B\cap\im(\varphi))\subseteq C$, we have an isomorphism $\varphi(C)/(\varphi(C)\cap B)\cong C/C'$. By assumption, $[B:C']<\infty$, which implies $[C:C']<\infty$, so that $[\varphi(C):\varphi(C)\cap B]<\infty$. Now $[G:\varphi(C)\cap B]=[G:\varphi(C)][\varphi(C):\varphi(C)\cap B]<\infty$, so that $[G:B]<\infty$. Now let $D\in \cC$. We have $[G:D]\leq [G:D\cap B]=[G:B][B:D\cap B]<\infty$, so \eqref{eqn:(FI)} holds.
\eproof
\setlength{\parindent}{0cm} \setlength{\parskip}{0.5cm}

We are now ready for the proof of our theorem.
\bproof[Proof of Theorem~\ref{thm:PI}]
By \cite[Lemma~4.1]{Li:IMRN}, it suffices to prove that every non-empty basic open set $\bd(kB;\{k_iB_i\})$ is properly infinite.
Suppose we are given $\bd(kB;\{k_iB_i\})\neq\emptyset$, where $kB\in \cE^\times$ and $\{k_iB_i\}$ is a finite (possibly empty) collection of constructible cosets. By conjugating by the homeomorphism $\bbt_k$ if necessary, we may assume $k=e$. We may also assume $k_iB_i\subsetneq B$ for all $i$.

First, let us suppose $[B:C]<\infty$ for all $C\in\cC$ with $C\subseteq B$. By Lemma~\ref{lem:fibasics}, we may assume $\{k_iB_i\}=\emptyset$ in this case. By Lemma~\ref{lem:PI}, $[G:B]<\infty$. Since $\sigma\colon S\acts G$ is non-automorphic, there exists $s\in S$ such that $\sigma_sG\lneq G$. Put $C:=B\cap\sigma_s^{-1}B$. Then $C$ is a constructible subgroup of $B$ satisfying $\sigma_sC\subseteq B$, and we have that $C,\sigma_sC\subseteq G$ are of finite index. Since $\sigma_sG/\sigma_sC=\sigma_s(G/C)$, we have $[G:\sigma_sC]=[G:\sigma_sG][\sigma_sG:\sigma_sC]=[G:C][G:\sigma_sG]$ and all indices are finite.
Thus, 
\[
[G:B][B:\sigma_sC]=[G:\sigma_sC]=[G:C][G:\sigma_sG]=[G:B][B:C][G:\sigma_sG],
\]
where all indices are finite, so that $[B:\sigma_sC]=[B:C][G : \sigma_sG]$. Let $g_1,...,g_m\in B$ be a complete set of representatives for $B/C$. Since $[G:\sigma_sG]\geq 2$, we can find $h_1,...,h_{2m}\in B$ such that the cosets $h_j\sigma_sC$ are pairwise disjoint for $1\leq j\leq 2m$.
Then
\[
U:=\bigsqcup_{j=1}^m\{\bbt_{h_j}\sigma_s\bbt_{g_j^{-1}}\}\times\bd(g_jC)\quad\text{ and }\quad V:=\bigsqcup_{j=1}^m\{\bbt_{h_{m+j}}\sigma_s\bbt_{g_j^{-1}}\}\times\bd(g_jC)
\]
are compact open bisections. By Lemma~\ref{lem:CharTightChar}, 
\[
s(U)=s(V)=\bigsqcup_{j=1}^m\bd(g_jC)=\bd(B).
\]
Moreover,
\[
r(U)=\bigsqcup_{j=1}^m\bd(h_j\sigma_sC)\quad\text{ and }\quad r(V)=\bigsqcup_{j=1}^m\bd(h_{m+j}\sigma_sC)
\]
are disjoint subsets of $\bd(B)$ by our choice of $h_j$'s. Thus, $\bd(B)$ is properly infinite.

Now let us suppose there exists a constructible subgroup $C'\subseteq B$ such that $[B:C']=\infty$. By Lemma~\ref{lem:bd()Inbd(.,.)}, there exists a constructible coset $hD\subseteq B\setminus\bigcup_ik_iB_i$, so that $\bd(hD)\subseteq \bd(B;\{k_iB_i\})$. 
By replacing $D$ with $D\cap C'$, we may assume that $D\subseteq C'$, so that $[B:D]=\infty$. If $[D:D']<\infty$ for all constructible subgroups $D'\subseteq D$, then $[G:D]<\infty$ by Lemma~\ref{lem:PI}; since the inclusion map $B/D\to G/D$ is injective, this would then imply that $[B:D]<\infty$, a contradiction. Thus, there exists a constructible subgroup $D'\subseteq D$ such that $[D:D']=\infty$.
By Theorem~\ref{thm:minimal}, there exist a constructible subgroup $C\subseteq D'$ and $\varphi\in I^e$ with $C\subseteq\dom(\varphi)$ such that $\varphi(C)\subseteq B$ has finite index. Let $g_1,...,g_m$ be a complete set of representatives for $B/\varphi(C)$. Choose $h_1,...,h_{2m}\in D$ such that the cosets $h_iC$ are pairwise disjoint for $1\leq j\leq 2m$ (here we are using that $[D:C]=\infty$). 
Consider the compact open bisections
\[
U:=\bigsqcup_{j=1}^m\{\bbt_{hh_j}\varphi^{-1}\bbt_{g_j^{-1}}\}\times\bd(g_j\varphi(C))\quad\text{ and }\quad V:=\bigsqcup_{j=1}^m\{\bbt_{hh_{m+j}}\varphi^{-1}\bbt_{g_j^{-1}}\}\times\bd(g_j\varphi(C)).
\]
Using Lemma~\ref{lem:CharTightChar}, we have 
\[
s(U)=s(V)=\bigsqcup_{j=1}^m\bd(g_j\varphi(C))=\bd(B)\supseteq \bd(B;\{k_iB_i\}).
\] 
Moreover, 
\[
r(U)=\bigsqcup_{j=1}^m\bd(hh_jC) \quad\text{ and } \quad r(V)=\bigsqcup_{j=1}^m\bd(hh_{m+j}C)
\]
are disjoint by our choice of $h_j$'s and are contained in $\bd(hD)$ because $hh_kC\subseteq hD$ for all $k$. Now put $U':=(s\vert_U)^{-1}(\bd(B;\{k_iB_i\}))$ and $V':=(s\vert_V)^{-1}(\bd(B;\{k_iB_i\}))$. Then $U'$ and $V'$ are compact open bisections satisfying $s(U')=s(V')=\bd(B;\{k_iB_i\})$ and $r(U')\sqcup r(V')\subseteq \bd(hD)\subseteq\bd(B;\{k_iB_i\})$. Thus, $\bd(B;\{k_iB_i\})$ is properly infinite.
\eproof

\section{Comparison of the groupoid model and the concrete C*-algebra}
\label{sec:C*-Comparison}
We now compare the C*-algebras of our groupoid $\cG_{\sigma}$ with $\fA_\sigma$. Throughout this section, $\sigma\colon S\acts G$ will be a non-automorphic algebraic action with $S$ and $G$ countable. 

\subsection{Comparison with the essential groupoid C*-algebra}

We shall first describe $C_{\es}^*(\cG_\sigma)$ using induced representations. For this, we need some preliminary results.
We let $\msc{S}_e:=(\cG_\sigma)_{\chi_e}^{\chi_e}$, where $\chi_e$ is the character from Definition~\ref{def:chi_z}.

\blemma
\label{lem:phichig=chig}
Suppose we have $g,h\in G$ and $\phi\in I^\times$ with $g\in\dom(\phi)$. Then $\phi.\chi_g=\chi_{\phi(g)}$, and $\chi_h=\chi_g$ if and only if $hG_c=gG_c$.
\elemma
\setlength{\parindent}{0cm} \setlength{\parskip}{0cm}

\bproof
The first claim is easy to see. 
We have $\chi_h=\chi_g$ if and only if $\chi_h(kC)=\chi_g(kC)$ for all $kC\in\cE^\times$, i.e., $h\in kC$ if and only if $g\in kC$ for all $kC\in\cE^\times$. This in turn is equivalent to having $hC=gC$ for all $C\in\cC$, i.e., $g^{-1}h\in C$ for all $C\in\cC$.  
\eproof
\setlength{\parindent}{0cm} \setlength{\parskip}{0.5cm}

\blemma
\label{lem:isogroup}
We have $\msc{S}_e=\{[\bbt_h\varphi,\chi_e] : h\in G_c, \varphi\in I^e\}$. 
\elemma
\setlength{\parindent}{0cm} \setlength{\parskip}{0cm}

\bproof
It follows from Lemma~\ref{lem:phichig=chig} that $[\bbt_h\varphi,\chi_e] \in\msc{S}_e$ for all $h\in G_c$ and $\varphi\in I^e$. Suppose $[\phi,\chi_e]\in\msc{S}_e$. Then $\phi.\chi_e=\chi_e$, so in particular $e\in\dom(\phi)$ which by Proposition~\ref{prop:stform} implies that $\dom(\phi)$ is a constructible subgroup and $\phi=\bbt_h\varphi$ for some $h\in G$ and $\varphi\in I^e$. Since $\chi_e=\phi.\chi_e=\chi_{\phi(e)}=\chi_h$, we have $h\in G_c$ by Lemma~\ref{lem:phichig=chig}.
\eproof
\setlength{\parindent}{0cm} \setlength{\parskip}{0.5cm}

Let us introduce some notation following \cite{ExelPardo}. For $\phi\in I^\times$, let $F_\phi:=\{\chi\in\bd : \chi(\dom(\phi))=1,\phi.\chi=\chi\}$ be the set of \emph{fixed characters} of $\phi$, and we let
\[
TF_\phi:=\bigcup_{gC\in\cE^\times,\, gC\subseteq\fix(\phi)}\bd(gC)
\]
be the set of \emph{trivially fixed characters} of $\phi$. It is straightforward to see that $TF_\phi\subseteq F_\phi$.

\blemma
\label{lem:notdangerous}
If $h\in G_c$ and $\varphi\in I^e$ with $[\bbt_h\varphi,\chi_e]\in \msc{S}_e\setminus \{\chi_e\}$, then $[\bbt_h\varphi,\bd(\dom(\varphi))]\cap\bd=\emptyset$.
\elemma 
\setlength{\parindent}{0cm} \setlength{\parskip}{0cm}

\bproof 
By \cite[Proposition~3.14]{ExelPardo}, we have $[\bbt_h\varphi,\bd(\dom(\bbt_h\varphi))]\cap\bd=[\bbt_h\varphi,TF_{\bbt_h\varphi}]$, so it suffices to show that $TF_{\bbt_h\varphi}=\emptyset$. 
Since $[\bbt_h\varphi,\chi_e]\neq \chi_e$, we have $\bbt_h\varphi\vert_D\neq \id_D$ for all $D\in\cC$. 
Suppose there exists $gC\in\cE^\times\setminus\cC$ with $gC\subseteq\fix(\bbt_h\varphi)$.
Then $g\in \dom(\varphi)=\dom(\bbt_h\varphi)$, and since $\varphi\in I^e$, $\dom(\varphi)$ is a (constructible) subgroup, so we then get $C\subseteq \dom(\varphi)$. Since $\varphi$ is a group homomorphism on its domain (see Remark~\ref{rmk:partialgphom}), we have
$gc=\bbt_h\varphi(gc)=h\varphi(g)\varphi(c)=hg\varphi(c)$
for every $c\in C$. Taking $c=e$ gives $g=hg$, so that $h=e$. But now we have $gc=g\varphi(c)$ for all $c\in C$, which implies $\varphi\vert_C=\id_C$, a contradiction. Hence, $\fix(\bbt_h\varphi)$ contains no constructible cosets, so that $TF_{\bbt_h\varphi}=\emptyset$.
\eproof
\setlength{\parindent}{0cm} \setlength{\parskip}{0.5cm}

Let $\lambda_{\msc{S}_e}$ be the left regular representation of the group C*-algebra $C^*(\msc{S}_e)$, and let $\Ind \lambda_{\msc{S}_e}$ be the representation of $C^*(\cG_\sigma)$ on $\ell^2((\cG_\sigma)_{\chi_e})$ induced from $\lambda_{\msc{S}_e}$.
Inspired by a recent result in the setting of semigroup C*-algebras from \cite{NS}, we now prove the following:

\bprop[cf. {\cite[Proposition~2.5]{NS}}]
\label{prop:Indchi}
We have $C_{\es}^*(\cG_\sigma)=(\Ind\lambda_{\msc{S}_e})(C^*(\cG))$.
\eprop
\setlength{\parindent}{0cm} \setlength{\parskip}{0cm}

\bproof
First, we claim that if $(\chi_i)_i$ is a net in $\bd$ that converges to $\chi_k$ for some $k\in G$, then $(\chi_i)_i$ does not converge to any point of $(\cG_\sigma)_{\chi_k}^{\chi_k}\setminus\{\chi_k\}$. 
In the terminology from \cite[\S~7]{KM}, this says that none of the points in $\{\chi_k : k\in G\}$ are dangerous.
Since $\bbt_k^{-1}$ is a homeomorphism of $\bd$ taking $\chi_k$ to $\chi_e$, it suffices to consider the case $k=e$.
Suppose $[\phi,\chi_e]\in \msc{S}_e\setminus \{\chi_e\}$. By Lemma~\ref{lem:isogroup}, there exists $h\in G_c$ and $\varphi\in I^e$ such that $\phi=\bbt_h\varphi$. We need to show that $(\chi_i)_i$ does not converge to $[\phi,\chi_e]$. But this follows immediately from Lemma~\ref{lem:notdangerous}.
\setlength{\parindent}{0cm} \setlength{\parskip}{0.5cm}

Since none of the points in $\{\chi_k : k\in G\}$ are dangerous, \cite[Proposition~1.12]{NS} implies that $J_\sing=\bigcap_{k\in G}\ker(\pi_{\chi_k})$ (note that the set $D_0$ in the statement of \cite[Proposition~1.12]{NS} is contained in the set of dangerous points). Here we used that $\cG_\sigma$ can be covered by countably many open bisections since $S$ and $G$ are assumed to be countable.
Since every $\pi_{\chi_k}$ is unitarily equivalent to $\pi_{\chi_e}$, it follows that $J_\sing=\ker(\pi_{\chi_e})$. Since $\Ind\lambda_{\msc{S}_e}=\pi_{\chi_e}\circ\pi_r$, we are done.
\eproof
\setlength{\parindent}{0cm} \setlength{\parskip}{0.5cm}

Let $\pi_\es\colon C^*(\cG_\sigma)\to C_{\es}^*(\cG_\sigma)$ and $\pi_r \colon C^*(\cG_\sigma)\to C^*_r(\cG_\sigma)$ be the canonical projection maps.
When $\cG_\sigma$ is topologically free, the C*-algebra $C_{\es}^*(\cG_\sigma)$ enjoys the following co-universal property:
\bprop[{\cite[Proposition~5.8~\&~Theorem~7.29]{KM}}]
\label{prop:essential}
Suppose $\cG_\sigma$ is topologically free (e.g., if $\sigma\colon S\acts G$ is exact). If $\pi\colon C^*(\cG_\sigma)\to B$ is a *-homomorphism to a C*-algebra $B$ such that $\pi\vert_{C(\bd)}$ is injective, then there exists a *-homomorphism $\pi(C^*(\cG_\sigma))\to C_{\es}^*(\cG_\sigma)$ such that $\pi(a)\mapsto \pi_\es(a)$ for all $a\in C^*(\cG_\sigma)$. In particular, there exists a *-homomorphism $\fA_\sigma\to C_{\es}^*(\cG_\sigma)$ such that $\Lambda_\phi\mapsto \pi(v_\phi)$ for all $\phi\in I$.
\eprop
\setlength{\parindent}{0cm} \setlength{\parskip}{0cm}

\bproof
This is a direct application of \cite[Proposition~5.8]{KM} together with \cite[~Theorem~7.29]{KM}. 
\eproof
\setlength{\parindent}{0cm} \setlength{\parskip}{0.5cm}

In order to compare $C_{\es}^*(\cG_\sigma)$ with $\fA_\sigma$, we shall now also describe $\fA_\sigma$ using induced representations.
Let $\rho: \: C^*(\cG_\sigma) \to \fA_\sigma$ be the *-homomorphism from Proposition~\ref{prop:Lambdatight}. Then $\rho$ is determined by $\rho(v_{\phi})=\Lambda_{\phi}$, where for $\phi \in I$, we let $v_{\phi} \defeq 1_{[\phi,\partial \widehat{\cE}(\dom(\phi))]} \in \cC(\cG_\sigma)$. 

Recall that if $\pi_1$ and $\pi_2$ are representations of $C^*(\cG_\sigma)$, then $\pi_1$ is said to be weakly contained in $\pi_2$ (written $\pi_1 \preceq \pi_2$) if $\ker \pi_2\subseteq\ker \pi_1$. Clearly, $\pi_1 \preceq \pi_2$ if and only if $\pi_1$ factors through $\pi_2$. For instance, $\rho\preceq \pi_\es$ if and only if there exists a *-homomorphism $C_\es^*(\cG_\sigma)\to\fA_\sigma$ such that $v_\phi\mapsto \Lambda_\phi$.

Consider the subgroup $\check{\msc{S}}_e:=\{[\phi,\chi_e] : \phi\in I^e\}\subseteq \msc{S}_e$. 
By Proposition~\ref{prop:indqr}, $\Ind\lambda_{\msc{S}_e/\check{\msc{S}}_e}$ is (unitarily equivalent to) the representation on $\ell^2((\cG_\sigma)_{\chi_e}/\check{\msc{S}}_e)$ given by 
\begin{equation*}
	\left((\Ind \lambda_{\msc{S}_e/\check{\msc{S}}_e})(f)\xi\right)([\bbt_k,\chi_e]\check{\msc{S}}_e)=\sum_{[\psi,\chi]\in\cG_\sigma^{\chi_k}}f([\psi,\chi])\xi([\psi^{-1}\bbt_k,\chi_e]\check{\msc{S}}_e)
\end{equation*}
for all $f\in \cC(\cG_\sigma)$ and all $\xi\in\ell^2((\cG_\sigma)_{\chi_e}/\check{\msc{S}}_e)$.

The essential observation is the following.

\bprop
\label{prop:rhoisinduced}
The representation $\rho$ is unitarily equivalent to $\Ind \lambda_{\msc{S}_e/\check{\msc{S}}_e}$. 
\eprop
\setlength{\parindent}{0cm} \setlength{\parskip}{0cm}

\bproof
First, we will show that the map $G\to (\cG_\sigma)_{\chi_e}/\check{\msc{S}}_e$ given by $g\mapsto [\bbt_g,\chi_e]\check{\msc{S}}_e$ is bijective. 
Let $[\phi,\chi_e]\in (\cG_\sigma)_{\chi_e}$. 
Since $\chi_e(\dom(\phi))=1$ (i.e., $e\in\dom(\phi)$), it follows from Proposition~\ref{prop:stform} that we can write $\phi=\bbt_{h}\varphi$ for some $h\in G$ and $\varphi\in I^e$. We have $[\phi,\chi_e]\check{\msc{S}}_e=[\bbt_h,\chi_e]\check{\msc{S}}_e$, so our map is surjective.
If $[\bbt_g,\chi_e]\check{\msc{S}}_e=[\bbt_h,\chi_e]\check{\msc{S}}_e$ for $g,h\in G$, then there exists $\varphi\in I^e$ such that $[\bbt_g,\chi_e]=[\bbt_h\varphi,\chi_e]$, so that there exists $C\in \cC$ with $\bbt_g\vert_C=(\bbt_h\varphi)\vert_C$. Evaluating at $e\in C$ gives $g=h$. Thus, our map is injective.
\setlength{\parindent}{0cm} \setlength{\parskip}{0.5cm}

Second, we will show that the unitary $\ell^2(G)\cong \ell^2((\cG_\sigma)_{\chi_e}/\check{\msc{S}}_e)$ associated with the above bijection intertwines $\rho$ and $\Ind \lambda_{\msc{S}_e/\check{\msc{S}}_e}$.
For $\phi\in I$, we have
\[
\left((\Ind \lambda_{\msc{S}_e/\check{\msc{S}}_e})(v_\phi)\delta_{[\bbt_g,\chi_e]\check{\msc{S}}_e}\right)([\bbt_k,\chi_e]\check{\msc{S}}_e)=\sum_{[\psi,\chi]\in\cG_\sigma^{\chi_k}}v_\phi([\psi,\chi])\delta_{[\bbt_g,\chi_e]\check{\msc{S}}_e}([\psi^{-1}\bbt_k,\chi_e]\check{\msc{S}}_e).
\] 
In order for the sum to be non-zero, there must exist $[\psi,\chi]\in\cG_\sigma$ such that $\psi.\chi=\chi_k$, $[\psi,\chi]\in [\phi,\bd(\dom(\phi))]$, and $[\bbt_g,\chi_e]\check{\msc{S}}_e=[\psi^{-1}\bbt_k,\chi_e]\check{\msc{S}}_e$. The condition $[\psi,\chi]\in [\phi,\bd(\dom(\phi))]$ means that $\chi(\dom(\phi))=1$ and there exists $lC\in\cE^\times$ with $\psi\vert_{lC}=\phi\vert_{lC}$ and $\chi(lC)=1$. Since $\psi.\chi=\chi_k$ is equivalent to $\chi=\psi^{-1}.\chi_k=\chi_{\psi^{-1}(k)}$, $\chi(\dom(\phi))=1$ implies $\psi^{-1}(k)\in\dom(\phi)$.  
The condition $[\bbt_g,\chi_e]\check{\msc{S}}_e=[\psi^{-1}\bbt_k,\chi_e]\check{\msc{S}}_e$ means that there exist $\varphi\in I^e$ and $D\in\cC$ such that $\bbt_g\vert_D=\psi^{-1}\bbt_k\varphi\vert_D$. Evaluating this equation at $e$ gives $g=\psi^{-1}(k)$. Moreover, $\psi^{-1}.\chi_k(lC)=\chi(lC)=1$ implies that $g=\psi^{-1}(k)$ lies in $lC\cap\dom(\psi)$, so that $\psi\vert_{lC}=\phi\vert_{lC}$ implies $\psi(g)=\phi(g)$.
Therefore, the sum is zero unless $g\in\dom(\phi)$ and $k=\phi(g)$, in which case there is a single non-zero summand corresponding to $[\phi,\chi_g]$. Thus, $(\Ind \lambda_{\msc{S}_e/\check{\msc{S}}_e})(v_\phi)\delta_{[\bbt_g,\chi_e]\check{\msc{S}}_e}=\delta_{[\bbt_{\phi(g)},\chi_e]\check{\msc{S}}_e}$ when $g\in\dom(\phi)$ and $(\Ind \lambda_{\msc{S}_e/\check{\msc{S}}_e})(v_\phi)\delta_{[\bbt_g,\chi_e]\check{\msc{S}}_e}=0$ otherwise.
\eproof
\setlength{\parindent}{0cm} \setlength{\parskip}{0.5cm}

\bcor
\label{cor:rho<Indchi}
If $\check{\msc{S}}_e$ is amenable, then $\rho\preceq\Ind\lambda_{\msc{S}_e}$. If $\rho\preceq\Ind\lambda_{\msc{S}_e}$ and $\cG_\sigma$ is topologically free (e.g., if $\sigma\colon S\acts G$ is exact), then 
\begin{equation}
 \label{eqn:esstofA}
C_{\es}^*(\cG_\sigma)=(\Ind\lambda_{\msc{S}_e})(C^*(\cG_\sigma))\to\fA_\sigma,\quad (\Ind\lambda_{\msc{S}_e})(v_\phi)\mapsto\Lambda_\phi,
\end{equation}
is an isomorphism. In particular, if $\cG_{\sigma}$ is Hausdorff, $\rho \preceq \pi_r$, and $\cG_\sigma$ is topologically free (e.g., if $\sigma\colon S\acts G$ is exact), then the *-homomorphism $\bar{\rho}\colon C_r^*(\cG_\sigma)\to\fA_\sigma, \, v_\phi \ma \Lambda_\phi$ is an isomorphism.
\ecor 
\setlength{\parindent}{0cm} \setlength{\parskip}{0cm}

\bproof
If $\check{\msc{S}}_e$ is amenable, then $\lambda_{\msc{S}_e/\check{\msc{S}}_e}\preceq\lambda_{\msc{S}_e}$, which, because weak containment is preserved under induction (cf. \cite[Lemma~2.1]{KS}), implies $\Ind\lambda_{\msc{S}_e/\check{\msc{S}}_e}\preceq\Ind\lambda_{\msc{S}_e}$. Since $\rho\sim_u \Ind\lambda_{\msc{S}_e/\check{\msc{S}}_e}$ by Proposition~\ref{prop:rhoisinduced}, it follows that $\rho\preceq \Ind\lambda_{\msc{S}_e}$. If $\cG_\sigma$ is topologically free (e.g., if $\sigma\colon S\acts G$ is exact), then the map in \eqref{eqn:esstofA} is invertible by Proposition~\ref{prop:essential}.
\eproof
\setlength{\parindent}{0cm} \setlength{\parskip}{0.5cm}

Combining Corollary~\ref{cor:rho<Indchi} with Corollary~\ref{cor:PI}, we obtain:

\bcor
\label{cor:essentialtofA}
Let $\sigma\colon S\acts G$ be a non-automorphic algebraic action with $S$ and $G$ countable. If $\sigma\colon S\acts G$ is exact, satisfies \eqref{eqn:(M)}, and $\check{\msc{S}}_e$ is amenable, then $\fA_\sigma$ is simple and purely infinite.
\ecor 

\bremark
Surprisingly, amenability of $\check{\msc{S}}_e$ is often also necessary for $\rho\preceq\Ind\lambda_{\msc{S}_e}$, see Proposition~\ref{prop:reducedtofA}.
\eremark

Let us now make a couple of observations about the group $\check{\msc{S}}_e$. Lemma~\ref{lem:isogroup} gives the following:

\bcor
\label{cor:S=checkS}
We have $\msc{S}_e=\check{\msc{S}}_e$ if and only if $\sigma\colon S\acts G$ is exact. 
\ecor

\bremark
\label{rmk:rho<Indpir} 
If we have a globalization $\tilde{\sigma}\colon\msc{S}\acts\msc{G}$ that satisfies \eqref{eqn:(JF)}, so that $\cG_\sigma\cong (\msc{G}\rtimes\msc{S})\ltimes\bd$ by Remark~\ref{rem:TrafoGPD}, then $\msc{S}_e\cong G_c\rtimes\gp{S}$ and $\check{\msc{S}}_e\cong\gp{S}$ by Lemma~\ref{lem:isogroup}, where $\gp{S}$ is the subgroup of $\msc{S}$ generated by $S$.

Any semigroup $P$ that can be embedded into a group admits a universal group embedding, i.e., there exists a group $G_{\textup{univ}}$, called the \emph{universal group of $P$}, and an embedding $P\to G_{\textup{univ}}$ such that for any homomorphism from $P$ to a group $H$ extends (uniquely) to a homomorphism $G_{\textup{univ}}\to H$ (see \cite[\S~12]{CP67} and the discussion in \cite[\S~5.4.1]{CELY}). Thus, if the universal group of the monoid $\{[\sigma_s,\chi_e] :s\in S\}$ is amenable, then $\check{\msc{S}}_e$ is amenable. If $S$ can be embedded into a group and the universal group of $S$ is amenable, then $\check{\msc{S}}_e$ is amenable.

The maximal group image of $I^e$ is defined to be the quotient of $I^e$ by the congruence $\phi\sim\psi$ if there exists $C\in\cC$ such that $\phi\vert_C=\psi\vert_C$. Since $[\phi,\chi_e]=[\psi,\chi_e]$ if and only if $\phi\sim\psi$, we see that the group $\check{\msc{S}}_e$ is canonically isomorphic to the maximal group image of $I^e$.
\eremark

When $I$ is strongly 0-E-unitary---which implies, in particular, that $C_r^*(\cG_\sigma)=C_\es^*(\cG_\sigma)$---we shall use recent results from \cite{CN} to characterize $\rho\preceq\Ind\lambda_{\msc{S}_e}$. 
Let $C_e^*(\msc{S}_e)$ denote the completion of the complex group algebra $\Cz\msc{S}_e$ with respect to the norm $||\cdot||_e$ from \cite[Definition~2.1]{CN}, i.e., $||\cdot||_e$ is given by
\[
||f||_e:=\sup\{||\pi(f)|| : \pi \text{ a representation of $C^*(\msc{S}_e)$ such that }\Ind\pi\prec \lambda_{\msc{S}_e}\}
\]
for all $f\in \Cz\msc{S}_e$.
Denote by $\lambda_{\msc{S}_e}^e$ the canonical projection $C^*(\msc{S}_e)\to C_e^*(\msc{S}_e)$. The following is a consequence of Proposition~\ref{prop:rhoisinduced}.

\bprop
\label{prop:reducedtofA}
Consider the following statements:
\setlength{\parindent}{0cm} \setlength{\parskip}{0cm}

\begin{enumerate}
    \item[\upshape(i)] $\check{\msc{S}}_e$ is amenable,
    \item[\upshape(ii)] $\lambda_{\msc{S}_e/\check{\msc{S}}_e}\preceq\lambda_{\msc{S}_e}^e$,
    \item[\upshape(iii)] $\rho\preceq \pi_r$.
\end{enumerate}
We always have (i) $\Rarr$ (ii) $\LRarr$ (iii). If $\lambda_{\msc{S}_e}^e=\lambda_{\msc{S}_e}$, then (ii) $\Rarr$ (i). In particular, if $I$ is strongly 0-E-unitary (e.g., if $\sigma\colon S\acts G$ has a globalization $\tilde{\sigma}\colon\msc{S}\acts\msc{G}$ which satisfies \eqref{eqn:(JF)}), then (i) $\LRarr$ (ii) $\LRarr$ (iii).
\eprop
\bproof
We always have $\lambda_{\msc{S}_e}\preceq\lambda_{\msc{S}_e}^e$, and amenability of $\check{\msc{S}}_e$ implies $\lambda_{\msc{S}_e/\check{\msc{S}}_e}\preceq\lambda_{\msc{S}_e}$ and thus  $\lambda_{\msc{S}_e/\check{\msc{S}}_e}\preceq\lambda_{\msc{S}_e}^e$.
\setlength{\parindent}{0.5cm} \setlength{\parskip}{0cm}

We have (ii) $\LRarr$ (iii), because Proposition~\ref{prop:rhoisinduced} implies that $\rho$ is unitarily equivalent to $\Ind\lambda_{\msc{S}_e/\check{\msc{S}}_e}$, so \cite[Proposition~2.2]{CN} implies that $\rho\preceq \pi_r$ if and only if $\lambda_{\msc{S}_e/\check{\msc{S}}_e}$ factors through $C_e^*(\msc{S}_e)$, i.e., $\lambda_{\msc{S}_e/\check{\msc{S}}_e}\preceq\lambda_{\msc{S}_e}^e$. 

If $\lambda_{\msc{S}_e}^e=\lambda_{\msc{S}_e}$, then (ii) implies $\lambda_{\msc{S}_e/\check{\msc{S}}_e}\preceq\lambda_{\msc{S}_e}$, which in turns implies that $\check{\msc{S}}_e$ is amenable.

If $I$ is strongly 0-E-unitary, then $\cG_\sigma$ is a partial transformation groupoid (see \cite[Lemma~5.5.22.]{CELY}, for instance), so $||\cdot||_e$ and $||\cdot||_r$ coincide by \cite[Corollary~4.15]{CN}, i.e., $\lambda_{\msc{S}_e}^e=\lambda_{\msc{S}_e}$.
\eproof
\setlength{\parindent}{0cm} \setlength{\parskip}{0.5cm}

\bremark
In general, it is not clear when $||\cdot||_e = ||\cdot||_r$ on $\Cz\msc{S}_e$, i.e., when $\lambda_{\msc{S}_e}^e=\lambda_{\msc{S}_e}$. 
\eremark

\bcor
\label{cor:Summary}
Assume that $\sigma\colon S\acts G$ has a globalization $\tilde{\sigma}\colon\msc{S}\acts\msc{G}$ which satisfies \eqref{eqn:(JF)}, $\sigma$ is exact and one of the conditions in Theorem~\ref{thm:minimal} is satisfied. Then $\spkl{S}$ is amenable if and only if $\fA_{\sigma}$ is simple if and only if the map $\fA_\sigma\to C_r^*(\cG_\sigma) = C_{\es}^*(\cG_\sigma), \, \Lambda_\phi\mapsto v_\phi$ from Proposition~\ref{prop:essential} is an isomorphism.
\ecor

We will conclude this section by explaining how $\fA_\sigma$ is, in many cases, an exotic groupoid C*-algebra.

\blemma
\label{lem:rhoisomimpliesGamenable}
If $I$ is strongly 0-E-unitary and $\rho\colon C^*(\cG_\sigma)\to\fA_\sigma$ is an isomorphism, then $G$ is amenable.
\elemma
\setlength{\parindent}{0cm} \setlength{\parskip}{0cm}

\bproof
The map $G\ltimes\bd\to\cG_\sigma$ given by $(g,\chi)\mapsto [\bbt_g,\chi]$ identifies the transformation groupoid $G\ltimes\bd$ with the clopen subgroupoid $\{[\bbt_g,\chi]\in \cG_\sigma : g\in G,\chi\in\bd\}$ of $\cG_\sigma$. We thus have a canonical *-homomorphism $C(\bd)\rtimes G\cong C^*(G\ltimes\bd)\to C^*(\cG_\sigma)$. In fact, this map is injective (see, e.g., \cite[Exercise~3.3.6]{Renault2}). It is easy to see that this embedding sends the canonical unitary $u_g$ in $C(\bd)\rtimes G$ corresponding to $g\in G$ to $v_{\bbt_g}$.
Moreover, the canonical *-homomorphism $C^*(G) \to C(\bd)\rtimes G$ is injective because there is a $G$-invariant probability measure on $\bd$. Now the composition of these canonical embeddings with $\rho$ coincides with the left regular representation $\lambda_G$ of $C^*(G)$. Thus, if $\rho$ is an isomorphism, then $\lambda_G$ must also be an isomorphism, in which case $G$ is amenable.
\eproof
\setlength{\parindent}{0cm} \setlength{\parskip}{0.5cm}

\bcor
\label{cor:exotic}
Suppose $\sigma\colon S\acts G$ is exact and $I$ is strongly 0-E-unitary. If both $G$ and $\check{\msc{S}}_e$ are non-amenable, then $\fA_\sigma$ is an exotic groupoid C*-algebra, in the sense that it sits properly between the full and reduced C*-algebra of $\cG_{\sigma}$.
\ecor
\setlength{\parindent}{0cm} \setlength{\parskip}{0cm}

\bproof
Our assumptions imply that $\cG_{\sigma}$ is Hausdorff. Hence Proposition~\ref{prop:essential} produces the canonical projection map $\fA_\sigma\twoheadrightarrow C_r^*(\cG_\sigma)$. It is not injective if $\check{\msc{S}}_e$ is non-amenable by Proposition~\ref{prop:reducedtofA}. If $G$ is non-amenable, then $\rho$ is not injective by Lemma~\ref{lem:rhoisomimpliesGamenable}. Our claim follows.
\eproof

Example~\ref{ex:exotic} contains a concrete example class where the hypotheses in Corollary~\ref{cor:exotic} are satisfied.
\setlength{\parindent}{0cm} \setlength{\parskip}{0.5cm}

Note that, by Proposition~\ref{prop:idpure}, $I$ is strongly 0-E-unitary if $\sigma\colon S\acts G$ has a globalization $\tilde{\sigma}\colon\msc{S}\acts\msc{G}$ which satisfies \eqref{eqn:(JF)}.

\section{Comparisons with the boundary quotient of the semigroup C*-algebra}
\label{sec:semigpC*}
It is natural to compare our C*-algebras $C_\es^*(\cG_\sigma)$, $C_r^*(\cG_\sigma)$, and $\fA_\sigma$ to the boundary quotient of the semigroup C*-algebra $C_\lambda^*(P)$, where $P=G\rtimes S$. For background on semigroup C*-algebras and their boundary quotients, see \cite[\S~2]{Li:GarI} and \cite[Chapter~5]{CELY}.

\subsection{Comparison of groupoids} 
Let $I_l$ be the left inverse hull of $P$, $E$ the semilattice of idempotents of $I_l$, and $I_l \ltimes \partial \widehat{E}$ the associated boundary groupoid, so that $C_r^*(I_l \ltimes \partial \widehat{E})$ is the boundary quotient of the semigroup C*-algebra $C_\lambda^*(P)$. 
Let $\cJ_S$ denote the semilattice of constructible right ideals of $S$. A straightforward computation shows
\[
(g_1,s_1)^{-1} (g_2,s_2) \dotsm (g_{n-1},s_{n-1})^{-1} (g_n,s_n) = (\sigma_{s_1}^{-1} \bbt_{g_1^{-1}} \bbt_{g_2} \sigma_{s_2} \dotsm \sigma_{s_{n-1}}^{-1} \bbt_{g_{n-1}^{-1}} \bbt_{g_n} \sigma_{s_n}) \times (s_1^{-1} s_2 \dotsm s_{n-1}^{-1} s_n)
\]
for all $g_1, \dotsc, g_n \in G$ and $s_1, \dotsc, s_n \in S$. Here we identify $G \rtimes S$ with its canonical copy inside $I_l$. 
If $S$ is left reversible (i.e., $sS \cap tS \neq \emptyset$ for all $s,t \in S$), then $\emptyset\notin\cJ_S$ by \cite[Lemma~5.6.43]{CELY}, and the projection onto the $G$-component defines an inverse semigroup homomorphism $I_l \to I, \, \Phi \ma \Phi_G$. It is straightforward to see that this map is surjective. By restricting the map $I_l\to I$, we obtain a surjective semilattice homomorphism $E \to \cE$ and hence a continuous embedding $\widehat{\cE} \into \widehat{E}, \, \chi \ma \ti{\chi}$. 

For the remainder of this subsection, we assume that $S$ is left reversible.

\blemma
The map $\widehat{\cE} \into \widehat{E}, \, \chi \ma \ti{\chi}$ restricts to a bijection $\widehat{\cE}_{\max} \cong \widehat{E}_{\max}$.
\elemma
\setlength{\parindent}{0cm} \setlength{\parskip}{0cm}

\bproof
Recall that $\chi \in \widehat{\cE}$ is maximal if and only if whenever $\chi(gC) = 0$ for some $gC \in \cE$, there exists $hD \in \cE$ with $\chi(hD) = 1$ and $gC\cap hD = \emptyset$. Moreover, every element of $E$ is of the form $gC \times X$ for some $gC \in \cE$ and $X\in\cJ_S$. Since $S$ is left reversible, we have $(gC \times X) \cap (hD \times Y) = \emptyset$ if and only if $gC\cap hD = \emptyset$.
\setlength{\parindent}{0cm} \setlength{\parskip}{0.5cm}

Now suppose that $\chi \in \widehat{\cE}_{\max}$. Then $\ti{\chi}(gC \times X) = 1$ if and only if $\chi(gC) = 1$. Assume that $\ti{\chi}(gC \times X) = 0$. Then $\chi(gC) = 0$. Hence there exists $hD \in \cE$ with $\chi(hD) = 1$ and $gC\cap hD= \emptyset$. It follows that, for every $Y$ such that $hD \times Y \in E$, we have $\ti{\chi}(hD \times Y) = 1$. At the same time, $(gC \times X) \cap (hD \times Y) = \emptyset$. This shows that $\ti{\chi}$ is maximal.

Now take $\omega \in \widehat{E}_{\max}$. Define $\chi \in \widehat{\cE}$ by $\chi(gC) = 1$ if $\omega(gC \times X) = 1$ for some $X\in\cJ_S$. Then $\chi \in \widehat{\cE}_{\max}$. We claim that $\ti{\chi} = \omega$. Indeed, $\chi(gC) = 1$ if and only if $\omega(gC \times X) = 1$ for some $X$ if and only if $\omega(gC \times Z) = 1$ for all $Z$ such that $gC \times Z \in E$. The last equivalence follows from maximality of $\omega$.
\eproof
\setlength{\parindent}{0cm} \setlength{\parskip}{0.5cm}

\bcor
The map $\widehat{\cE} \into \widehat{E}, \, \chi \ma \ti{\chi}$ restricts to a homeomorphism $\partial \widehat{\cE} \cong \partial \widehat{E}$.
\ecor

Now we see that there is a canonical surjection $I_l\ltimes \partial \widehat{E}\twoheadrightarrow\cG_\sigma$ given by $[\Phi,\tilde{\chi}]\mapsto [\Phi_G,\chi]$.

\blemma
If the left inverse hull $I_l(S)$ of $S$ is E-unitary, then the surjection $I_l \onto I$, $\Phi\mapsto\Phi_G$, is an isomorphism if and only if the surjection $I_l\ltimes\partial \widehat{E}\onto \cG_\sigma$, $[\Phi,\tilde{\chi}]\mapsto [\Phi_G,\chi]$, is an isomorphism.
\elemma
\setlength{\parindent}{0cm} \setlength{\parskip}{0cm}

\bproof
\an{$\Rarr$} is clear. For \an{$\Larr$}, take $\Phi$ and assume that $\Phi_G \in \cE$, say $\Phi_G = \id_{gC}$ and $\chi(gC) = 1$, so that $[\Phi_G,\chi] \in (\cG_\sigma)^{(0)}$. Thus if $I_l\ltimes\partial \widehat{E} \onto \cG_\sigma$ is an isomorphism, we deduce $[\Phi,\ti{\chi}] = \ti{\chi}$ and thus $\Phi \vert_{gC \times X} = \id_{gC\times X}$ for some $X\in\cJ_S$. Therefore, $\Phi = \Phi_G \times \Phi_S$ with $\Phi_G \in \cE$ and $\Phi_S \vert_{X} = \id_X$. As $I_l(S)$ is E-unitary, the latter implies that $\Phi_S \in E(S)$, and hence $\Phi \in E$, as desired.
\eproof
\setlength{\parindent}{0cm} \setlength{\parskip}{0.5cm}

\subsection{Comparison of C*-algebras}

Let us now compare the C*-algebras. For $\Phi\in I_l$, let $w_\Phi$ denote the corresponding partial isometry in $\cC(I_l\ltimes\partial\widehat{E})$. 

\bprop
\label{prop:Srev+homs}
The following are equivalent:
\setlength{\parindent}{0cm} \setlength{\parskip}{0cm}

\begin{enumerate}[\upshape(i)]
	\item $S$ is left reversible;
	\item there exists a *-homomorphism $\vartheta\colon C^*(I_l \ltimes \partial \widehat{E}) \to \fA_\sigma$ such that $\vartheta(w_\Phi)= \Lambda_{\Phi_G}$ for all $\Phi\in I_l$; 
	\item there exists a *-homomorphism $\theta\colon C^*(I_l \ltimes \partial \widehat{E}) \to C^*_r(\cG_\sigma)$ such that $\theta(w_\Phi)=v_{\Phi_G}$ for all $\Phi\in I_l$.
\end{enumerate}
\eprop
\bproof
For all $s, t \in S$, we have $\sigma_sG\cap\sigma_tG\neq\emptyset$, so that $\Lambda_{\sigma_s} \Lambda_{\sigma_s}^* \Lambda_{\sigma_t} \Lambda_{\sigma_t}^*\neq 0$ and $v_{\sigma_s} v_{\sigma_s}^* v_{\sigma_t} v_{\sigma_t}^* \neq 0$. Thus, left reversibility of $S$ is necessary for existence in both cases. 
\setlength{\parindent}{0cm} \setlength{\parskip}{0.5cm}

Now assume $S$ is left reversible. To prove (ii), it suffices to show that the representation $I_l\to \fA_\sigma$, $\Phi\mapsto \Lambda_{\Phi_G}$ is a tight representation of $I_l$ in the sense of \cite{Exel08}; since the restriction to $E$ is unital, it suffices by \cite[Corollary~4.3]{Exel21} to prove that this representation is cover-to-join in the sense of \cite[\S~3]{Exel21}. Let $gC\times X\in E^\times$ and suppose $\mfc\subseteq E$ is a finite cover of $gC\times X$. Then for every $hD\times Y\subseteq gC\times X$, there exists $kB\times Z\in\mfc$ such that $(kB\times Z)\cap (hD\times Y)\neq\emptyset$. Let $\mfc_G:=\{kB : kB\times Z\in\mfc\text{ for some }Z\in\cJ_S\}$. It is easy to see that $\mfc_G$ is a (finite) cover of $gC$. We have 
\[
\bigvee_{hD\times Y\in \mfc}\Lambda_{(\id_{hD\times Y})_G}=\bigvee_{hD\in \mfc_G}\Lambda_{\id_{hD}}=\Lambda_{\id_{gC}},
\]
where the last equality uses that $\phi\mapsto \Lambda_\phi$ is a cover-to-join representation of $I$ in $\fA_\sigma$ (see the proof of Proposition~\ref{prop:Lambdatight}).

The proof of (iii) is essentially the same, using that $I_l\to C_r^*(\cG_\sigma)$, $\Phi\mapsto v_{\Phi_G}$ is a tight representation.
\eproof
\setlength{\parindent}{0cm} \setlength{\parskip}{0.5cm}

\bremark
A special case of the equivalence of (i)$\Leftrightarrow$(ii) in Proposition~\ref{prop:Srev+homs} was observed in \cite[Proposition~4.3]{BS} using different methods.
\eremark

For the remainder of this section, we assume that $S$ is left reversible. Let $\vartheta\colon C^*(I_l \ltimes \partial \widehat{E})\to\fA_\sigma$ be the *-homomorphism from part (ii) of Proposition~\ref{prop:Srev+homs}.

We now compare $\fA_\sigma$ and $C_r^*(I_l \ltimes \partial \widehat{E})$. Let $\pi_l\colon C^*(I_l \ltimes \partial \widehat{E})\to C_r^*(I_l \ltimes \partial \widehat{E})$ be the canonical projection map, and put $\check{\msc{T}}_e:=\{[\Phi,\ti{\chi}_e] : \Phi_G\in I^e \}\subseteq \msc{T}_e \defeq (I_l \ltimes \partial \widehat{E})_{\ti{\chi}_e}^{\ti{\chi}_e}$. Note that $\check{\msc{T}}_e=\msc{T}_e$ if $\sigma\colon S\acts G$ is exact.

\bprop
\label{prop:varthetaisinduced}
The representation $\vartheta$ is unitarily equivalent to $\Ind \lambda_{\msc{T}_e/\check{\msc{T}}_e}$. 
\eprop
\setlength{\parindent}{0cm} \setlength{\parskip}{0cm}

\bproof
The projection $I_l\ltimes\partial\widehat{E}\twoheadrightarrow \cG_\sigma$ descends to a bijection $(I_l \ltimes \partial \widehat{E})_{\ti{\chi}_e}/\check{\msc{T}}_e\cong (\cG_\sigma)_{\chi_e}/\check{\msc{S}}_e$. Composing this with the bijection $(\cG_\sigma)_{\chi_e}/\check{\msc{S}}_e\cong G$ from the proof of Proposition~\ref{prop:rhoisinduced}, we get that the map $G\to (I_l \ltimes \partial \widehat{E})_{\ti{\chi}_e}/\check{\msc{T}}_e$ given by $g\mapsto [g,\ti{\chi}_e]\check{\msc{T}}_e$ is a bijection. 
Similarly to the proof of Proposition~\ref{prop:rhoisinduced}, it now follows from Proposition~\ref{prop:indqr} that the unitary $\ell^2(G)\cong \ell^2((I_l \ltimes \partial \widehat{E})_{\ti{\chi}_e}/\check{\msc{T}}_e)$ induced by the above bijection intertwines $\vartheta$ and $\Ind\lambda_{\msc{T}_e/\check{\msc{T}}_e}$.
\eproof
\setlength{\parindent}{0cm} \setlength{\parskip}{0.5cm}

Let $C_e^*(\msc{T}_e)$ denote the completion of the complex group algebra $\Cz\msc{T}_e$ with respect to the norm $||\cdot||_e$ defined in \cite[Definition~2.1]{CN}, and denote by $\lambda_{\msc{T}_e}^e$ the canonical surjection $C^*(\msc{T}_e)\to C_e^*(\msc{T}_e)$. The following is analogous to Proposition~\ref{prop:reducedtofA}, using the observation that if $P$ embeds into a group, then $I_l\ltimes\partial\widehat{E}$ is a partial transformation groupoid by \cite[\S~5.7]{CELY}, so that \cite[Corollary~4.15]{CN} implies $\lambda_{\msc{T}_e}^e=\lambda_{\msc{T}_e}$.

\bcor
\label{cor:Theta<Pil}
Consider the following statements:
\setlength{\parindent}{0cm} \setlength{\parskip}{0cm}

\begin{enumerate}
    \item[\upshape(i)] $\check{\msc{T}}_e$ is amenable,
    \item[\upshape(ii)] $\lambda_{\msc{T}_e/\check{\msc{T}}_e}\preceq \lambda_{\msc{T}_e}^e$,
    \item[\upshape(iii)] $\vartheta \preceq \pi_l$.
\end{enumerate}
We always have (i) $\Rarr$ (ii) $\LRarr$ (iii). If $\lambda_{\msc{T}_e}^e=\lambda_{\msc{T}_e}$, then (ii) $\Rarr$ (i). In particular, if $P$ embeds into a group, then (i) $\LRarr$ (ii) $\LRarr$ (iii).
\ecor
\setlength{\parindent}{0cm} \setlength{\parskip}{0.5cm}

We lastly compare $C^*_{\es}(\cG_\sigma)$ and $C_{\es}^*(I_l \ltimes \partial \widehat{E})$. For this, we also assume that $S$ and $G$ are countable, so that we can use results from \cite{NS}.

The surjection $I_l\ltimes\partial \widehat{E} \onto \cG_\sigma$ induces a projection $\msc{T}_e \onto \msc{S}_e$. Let $\msc{N} \subseteq \msc{T}_e$ be the kernel of this map. Explicitly, $\msc{N}=\{[\Phi,\ti{\chi}_e]\in \msc{T}_e : \text{ there exists } C\in \cC\text{ with }\Phi_G\vert_C=\id_C\}$.

\bprop
\label{prop:ess2}
We have $C_{\es}^*(I_l\ltimes\partial\widehat{E})=(\Ind\lambda_{\msc{T}_e})(C^*(I_l\ltimes\partial\widehat{E}))$.
\eprop
\setlength{\parindent}{0cm} \setlength{\parskip}{0cm}

\bproof
Similarly to the proof of Proposition~\ref{prop:Indchi} (using that $S$ and $G$ are countable), it suffices to show that the point $\tilde{\chi}_e$ is not dangerous. Since $\chi_e$ is not dangerous in $\cG_\sigma$, it is enough to show that for any given $[\Phi,\ti{\chi}_e]\in\msc{N}\setminus\{\ti{\chi}_e\}$, there is no net in $\partial\widehat{E}$ that converges to $\ti{\chi}_e$ and $[\Phi,\ti{\chi}_e]$. So suppose $[\Phi,\ti{\chi}_e]\in\msc{N}\setminus\{\ti{\chi}_e\}$. Then, $[\Phi,\partial\widehat{E}(\dom(\Phi))]\cap \partial\widehat{E}=[\Phi,TF_\Phi]$, where $TF_\Phi:=\bigcup_{\varepsilon\in E,\varepsilon\subseteq\fix(\Phi)}\partial\widehat{E}(\varepsilon)$. Thus, it suffices to show that $\fix(\Phi)$ does not contain any member of $E^\times$. Suppose $hD\times X\in E^\times$ with $hD\times X\subseteq\fix(\Phi)$. 
Then $hD\subseteq \fix(\Phi_G)$ and $X\subseteq \fix(\Phi_S)$. Since $[\Phi,\ti{\chi}_e]$ lies in $\msc{N}$, there exists $C\in\cC$ such that $\Phi_G\vert_C=\id_C$. In particular, this implies that $\Phi_G\in I^e$, so that $\dom(\Phi_G)$ is a subgroup and $\Phi_G$ is a homomorphism on its domain (see Remark~\ref{rmk:partialgphom}). It follows that $\Phi_G(hd)=\Phi_G(h)\Phi_G(d)$ for all $d\in D$. Since $\Phi_G(h)=h$, this implies $D\subseteq\fix(\Phi_G)$. Now we have $D\times X\subseteq\fix(\Phi)$, and $D\times X=(h,1)^{-1}(hD \times X)\in E^\times$. Since $\ti{\chi}_e(D\times X)=1$, this implies $[\Phi,\ti{\chi}_e]=\ti{\chi}_e$, which is a contradiction. 
\eproof
\setlength{\parindent}{0cm} \setlength{\parskip}{0.5cm}

Let $\bar{\theta}\colon C^*(I_l \ltimes \partial\widehat{E})\to C^*_{\es}(\cG_\sigma)$ be the composition of $\theta$ from Proposition~\ref{prop:Srev+homs} with the quotient map $C^*_r(\cG_\sigma)\to C^*_{\es}(\cG_\sigma)$.

\bprop
\label{prop:piisinduced}
The representation $\bar{\theta}$ is unitarily equivalent to $\Ind \lambda_{\msc{T}_e/\msc{N}}$.
\eprop
\setlength{\parindent}{0cm} \setlength{\parskip}{0cm}

\bproof
The projection $I_l\ltimes \partial\widehat{E}\twoheadrightarrow\cG_\sigma$ induces a bijection $(I_l\ltimes\partial\widehat{E})_{\ti{\chi}_e}/\msc{N}\cong (\cG_\sigma)_{\chi_e}$.
By Proposition~\ref{prop:indqr}, the representation $\Ind \lambda_{\msc{T}_e/\msc{N}}$ is unitarily equivalent to the canonical representation of $C^*(I_l\ltimes\partial\widehat{E})$ on $\ell^2((I_l\ltimes\partial\widehat{E})_{\ti{\chi}_e}/\msc{N})$, and by Proposition~\ref{prop:Indchi}, $C^*_{\es}(\cG_\sigma)=(\Ind\lambda_{\msc{S}_e})(C^*(\cG_\sigma))\subseteq \cB(\ell^2((\cG_\sigma)_{\chi_e}))$.
\setlength{\parindent}{0.5cm} \setlength{\parskip}{0cm}

It remains to check that the unitary $\ell^2((I_l\ltimes\partial\widehat{E})_{\ti{\chi}_e}/\msc{N})\cong \ell^2((\cG_\sigma)_{\chi_e})$ associated with the above bijection implements a unitary equivalence between $\bar{\theta}$ and $\Ind \lambda_{\msc{T}_e/\msc{N}}$. This is similar to the proof of Proposition~\ref{prop:rhoisinduced}.
\eproof
\setlength{\parindent}{0cm} \setlength{\parskip}{0.5cm}

The following is analogous to Proposition~\ref{prop:reducedtofA} and Corollary~\ref{cor:Theta<Pil}
\bcor
\label{cor:ess2}
Consider the following statements:
\setlength{\parindent}{0cm} \setlength{\parskip}{0cm}

\begin{enumerate}
    \item[\upshape(i)] $\msc{N}$ is amenable,
    \item[\upshape(ii)] $\lambda_{\msc{T}_e/\msc{N}}\preceq\lambda_{\msc{T}_e}^e$,
    \item[\upshape(iii)] $\bar{\theta} \preceq \Ind\lambda_{\msc{T}_e}$, i.e., there is a *-homomorphism $\theta_{\es}: \: C_{\es}^*(I_l\ltimes\partial\widehat{E})\to C_{\es}^*(\cG_\sigma)$ sending $w_\Phi$ to $v_{\Phi_G}$.
\end{enumerate}
We always have (i) $\Rarr$ (ii) and (i) $\Rarr$ (iii). If $C_r^*(I_l\ltimes\partial\widehat{E})=C_{\es}^*(I_l\ltimes\partial\widehat{E})$ (e.g., if $I_l\ltimes\partial\widehat{E}$ is Hausdorff), then (ii) $\LRarr$ (iii). If $\lambda_{\msc{T}_e}^e=\lambda_{\msc{T}_e}$ (e.g., if $P$ embeds into a group), then (ii) $\Rarr$ (i).
\ecor 
\setlength{\parindent}{0cm} \setlength{\parskip}{0.5cm}

\bremark
Hausdorffness of $I_l\rtimes\partial\widehat{E}$ is characterized in \cite[\S~4]{Li:GarI}. 
\eremark

Recall that we are assuming that $S$ is left reversible.
\blemma
\label{lem:piNotInj}
Suppose $\theta_{\es}$ in Corollary~\ref{cor:ess2}~(iii) exists (for instance if $\msc{N}$ is amenable) and that $\Phi\in \gp{S}\subseteq\msc{T}_e$ is such that $\Phi \vert_{C \times X} \neq \id_{C\times X}$ for all $\emptyset \neq C \times X \in E$ with $C \in \cC$, but $\Phi_G = \id$. Then $\theta_{\es}(w_{\Phi} - w_{\Phi^{-1} \Phi}) = 0$ but $w_{\Phi} - w_{\Phi^{-1} \Phi} \neq 0$ in $C_{\es}^*(I_l\ltimes\partial\widehat{E})$.
\elemma
\setlength{\parindent}{0cm} \setlength{\parskip}{0cm}

\bproof
First note that $(\Ind \lambda_{\msc{T}_e})(w_{\Phi})(\delta_{[\id,\ti{\chi}_e]}) = \delta_{[\Phi,\ti{\chi}_e]}$ while $(\Ind \lambda_{\msc{T}_e})(w_{\Phi^{-1} \Phi})(\delta_{[\id,\ti{\chi}_e]}) = \delta_{[\id,\ti{\chi}_e]}$. Now $[\Phi,\ti{\chi}_e] = \ti{\chi}_e$ if and only if $\Phi \vert_{C \times X} = \id_{C\times X}$ for some $\emptyset \neq C \times X \in E$ with $C \in \cC$.
At the same time, $\Phi_G = \id$ implies $\Phi_G = \Phi_G^{-1} \Phi_G$. Hence
\[
\theta_{\es}(w_{\Phi}) = (\Ind \lambda_{\msc{S}_e})(v_{\Phi_G}) = (\Ind \lambda_{\msc{S}_e})(v_{\Phi_G^{-1} \Phi_G}) = (\Ind \lambda_{\msc{S}_e})(v_{\Phi_G^{-1}}) (\Ind \lambda_{\msc{S}_e})(v_{\Phi_G}) = \theta_{\es}(w_{\Phi^{-1} \Phi}).\qedhere
\]
\eproof
\setlength{\parindent}{0cm} \setlength{\parskip}{0.5cm}

\bremark
If $I_l$ is 0-E-unitary, then our condition above for non-injectivity of $\theta_{\es}$ is satisfied whenever there exists $\Phi \notin E$ with $\Phi_G = \id$ (i.e., $\Phi_G \in \cE$), i.e., $I_l \onto I$ is not injective. In other words, if $I_l$ is 0-E-unitary, and if $\theta_{\es}$ is an isomorphism, then the map $I_l \onto I$ from above must be an isomorphism. In particular, again if $I_l$ is 0-E-unitary, then $\theta_{\es}$ is not injective whenever $I$ is not 0-E-unitary.
\eremark

\bprop
\label{prop:isomiffN=1}
Suppose $\theta_{\es}$ in Corollary~\ref{cor:ess2}~(iii) exists (for instance if $\msc{N}$ is amenable). Then $\theta_{\es}$ is an isomorphism if and only if $\msc{N}$ is trivial.
\eprop
\bproof
If $\msc{N}$ is trivial, then $\Ind\lambda_{{\msc{T}_e}/\msc{N}}=\Ind \lambda_{\msc{T}_e}$, and the result follows from Propositions~\ref{prop:ess2} and \ref{prop:piisinduced}. If $\msc{N}$ is non-trivial, then $\theta_{\es}$ is not injective because for any non-trivial $[\Phi,\ti{\chi}_e] \in \msc{N}$, $\Phi$ satisfies the condition in Lemma~\ref{lem:piNotInj}.
\eproof

\section{Examples}

\subsection{Algebraic actions with the finite index property}
\label{sec:FI}

\bdefin
We say that $\sigma\colon S\acts G$ has the \emph{finite index property} if 
\begin{equation}
	\label{eqn:(FI)}\tag{FI}
	\# (G/\sigma_sG) <\infty\quad \text{for all } s\in S.
\end{equation}
\edefin

\bprop
\label{prop:(FI)}
If $\sigma\colon S\acts G$ satisfies \eqref{eqn:(FI)}, then every member of $\cC$ is a finite index subgroup of $G$.
\eprop
\setlength{\parindent}{0cm} \setlength{\parskip}{0cm}

\bproof
We proceed by induction. The induction start is provided by \eqref{eqn:(FI)}. For the induction step, suppose $C\in\cC$ with $\# (G/C) <\infty$. Now let $s,t\in S$. 
Since $\sigma_sG/\sigma_sC=\sigma_s(G/C)$, we see that $[\sigma_sG:\sigma_sC]$ is finite by the induction hypothesis. Hence, since $[G:\sigma_sG]$ is also finite by the induction hypothesis, $[G:\sigma_sC]$ is finite (see, e.g., \cite[Chapter~I,~Theorem~4.5]{Hunger}).
Moreover, we have $\sigma_t(G/\sigma_t^{-1}C)=\sigma_tG/((\sigma_tG)\cap C)$, and the latter is finite because we have an embedding $\sigma_tG/((\sigma_tG)\cap C)\hookrightarrow G/C$.
\eproof
\setlength{\parindent}{0cm} \setlength{\parskip}{0.5cm}

We immediately obtain the following:
\bcor
If $\sigma\colon S\acts G$ satisfies \eqref{eqn:(FI)}, then $\bar{G}:=\varprojlim_{C\in\cC}G/C$ is compact. Moreover, every character in $\bd$ is maximal, and $\bar{G}$ coincides with $\partial\E$ (cf. Lemma~\ref{lem:filtersOnEvsOnC}). 
\ecor

\bcor
\label{cor:(FI)implies(M)}
If $\sigma\colon S\acts G$ satisfies \eqref{eqn:(FI)}, then $\cG_\sigma$ is minimal.
\ecor
\setlength{\parindent}{0cm} \setlength{\parskip}{0cm}

\bproof
It is clear that \eqref{eqn:(FI)} implies \eqref{eqn:(M4)} from Theorem~\ref{thm:minimal}.
\eproof
\setlength{\parindent}{0cm} \setlength{\parskip}{0.5cm}

\bprop
\label{prop:(FI)to(JF)}
Assume $G$ has the property that $g^m=h^m$ implies $h=g$ for all $g,h\in G$ and all $m\in\Zz_{>0}$ and that $\sigma\colon S\acts A$ satisfies \eqref{eqn:(FI)}. Then $I$ is 0-E-unitary. If $\sigma\colon S\acts G$ admits a globalization $\tilde{\sigma}\colon\msc{S}\acts\msc{G}$, then it satisfies \eqref{eqn:(JF)}.
\eprop
\setlength{\parindent}{0cm} \setlength{\parskip}{0cm}

\bproof
By Corollary~\ref{cor:0Eunitary}, it suffices to prove that $I^e$ is E-unitary. Suppose $\varphi\in I^e$ is such that $\varphi\vert_C=\id_C$ for some $C\in\cC$. By Proposition~\ref{prop:(FI)}, for every $g \in G$ there exists $m\in\Zz_{>0}$ such that $g^m \in C$. Now we have $\varphi(g)^m=\varphi(g^m)=g^m$, which implies $\varphi(g)=g$ for all $g\in G$ by our assumption. The proof that \eqref{eqn:(JF)} is satisfied is similar. 
\eproof

Note that if $G$ is Abelian and torsion free, then $g^m=h^m$ implies $h=g$ for all $g,h\in G$ and all $m\in\Zz_{>0}$. Moreover, if $I$ is 0-E-unitary, then $\cG_\sigma$ is Hausdorff by \cite[Corollary~3.17]{ExelPardo}.
\setlength{\parindent}{0cm} \setlength{\parskip}{0.5cm}

\bex[Algebraic actions on tori and solenoids]
\label{ex:solenoids}
Let $G$ be a torsion-free Abelian group of finite rank $r\in\Zz_{>0}$, so that we can view $G$ as a subgroup of $\mathscr{G} \defeq \Qz^r$. Note that these assumptions on $G$ are equivalent to the dual group $\widehat{G}$ being a solenoid. Given an algebraic action $\sigma\colon S\acts G$, every $\sigma_s$ extends naturally to an automorphism $\ti{\sigma}_s$ of $\mathscr{G}$, so that we obtain a natural globalization by considering the action of the subgroup $\mathscr{S}$ of $\Aut(\mathscr{G})$ generated by $\ti{\sigma}_s$, $s \in S$. Moreover, \eqref{eqn:(JF)} is satisfied by Proposition~\ref{prop:(FI)to(JF)}, so that $\cG_\sigma \cong (\msc{G}\rtimes\mathscr{S}) \ltimes \bd$. In particular, our groupoid is Hausdorff.  
By \cite[Exercise~92.5]{Fuchs2}, \eqref{eqn:(FI)} is satisfied, so $\cG_\sigma$ is minimal by Corollary~\ref{cor:(FI)implies(M)} and hence purely infinite by Theorem~\ref{thm:PI}.
\eex

\bremark
Example~\ref{ex:solenoids} provides many algebraic actions which have globalizations even though the acting semigroup is not left Ore. For instance, it is not difficult to find faithful non-automorphic (even exact) actions of free monoids on tori. Such actions are necessarily very far from respecting the order in the sense of \cite[Definition~8.1]{BLS}.
\eremark

\bremark
For a non-automorphic algebraic action $\sigma: \: S \acts \Zz^r$, property ID for the dual action in the sense of \cite{Ber,Muc} implies exactness for $\sigma$.
\eremark

\bex[Actions on infinite rank groups that satisfy \eqref{eqn:(FI)}]
\label{ex:ssaction}
Fix a left cancellative semigroup $S$. Suppose $J$ is a non-empty set and that for each $j\in J$ we have an algebraic action $\sigma^j\colon S\acts G_j$ satisfying \eqref{eqn:(FI)}. Then the diagonal action $\delta\colon S\acts \prod_{i\in J} G_j$ given by $(\delta_s(g_j))_j=\sigma_j(g_j)$ satisfies \eqref{eqn:(FI)} if and only if for every $s\in S$, $\sigma^j_s\in\End(G_j)$ is invertible for all but finitely many $j$. For instance, if $\cO_K$ is the ring of integers in a number field $K$, then the diagonal action $\cO_K^\times\acts \prod_{\p\in\cP_K}\cO_{K,[\p]}$ satisfies \eqref{eqn:(FI)}, where $\cP_K$ is the set of non-zero prime ideals of $\cO_K$ and $\cO_{K,[\p]}$ is the localization of $\cO_K$ at the $\p\in\cP_K$. 
\eex

\bex[Algebraic actions from self-similar actions of groups]
\label{ex:ssactions}
We briefly explain how to obtain examples of algebraic actions from self-similar group actions. We refer the reader to \cite{Nek} and \cite{LW} for background on self-similar actions.
Let $d\in\Zz_{>1}$ and $X:=\{0,...,d-1\}$. Let $X^*$ denote the free monoid on $X$, and for each $n\in\Nz$, let $X^n\subseteq X^*$ be the set of words of length $n$.
Suppose $G\acts X^*$ is a faithful self-similar action of a non-trivial group $G$ on $X^*$ as in \cite[\S~3]{LW} (cf. \cite{Nek}). 
For each $\mu\in X^*$, let $G_\mu\subseteq G$ be the stabilizer subgroup of $\mu$, and let $\phi_\mu\colon G_\mu\to G$ be the homomorphism $\phi_\mu(g):=g\vert_\mu$, where $g\vert_\mu$ is the section of $g$ at $\mu$. 
Assume that $\phi_x$ is an isomorphism for all $x\in X$. By \cite[Lemma~3.10]{LW}, this is equivalent to assuming $\phi_\mu$ is an isomorphism for all $\mu\in X^*$. For each $\mu\in X^*$, put $\sigma_\mu:=\phi_\mu^{-1}\colon G\to G_\mu$. For $\mu,\nu\in X^*$, we have $\sigma_\mu\circ\sigma_\nu=\sigma_{\mu\nu}$, so that $S=\{\sigma_\mu : \mu\in X^*\}$ is a submonoid of $\End(G)$.
Since $G\acts X^*$ is faithful, $\bigcap_{\mu\in X^*}G_\mu=\{g\in G : g(\mu)=\mu\text{ for all } \mu\in X^*\}=\{e\}$, so $S\acts G$ is exact. Since $[G:G_\mu]\leq d^{|\mu|}$ for all $\mu\in X^*$, $S\acts G$ satisfies \eqref{eqn:(FI)}. Thus, our groupoid $\cG_\sigma$ is topologically free, minimal, and purely infinite.
\eex

\bex[Ring C*-algebras of non-commutative rings]

Let $\cR$ be a unital (not necessary commutative) ring with $1\neq 0$, and let $\cR^\times$ be the set of left regular elements of $\cR$, i.e., $\cR^\times:=\{a\in \cR : ax=ay \text{ implies } x=y\text{ for all }x,y\in\cR\}$. Then $\cR^\times$ acts on the additive group of $\cR$ by injective endomorphisms. Since $\cR$ is unital, the algebraic action $\cR^\times\acts\cR$ is faithful. The concrete C*-algebra associated with $\cR^\times\acts\cR$ is called the \emph{reduced ring C*-algebra of $\cR$} (see \cite{Li:Ring}) and is denoted by $\fA_r[\cR]$. 
Assume that the additive group of $\cR$ is torsion-free and of finite rank. Examples of such rings include integral group rings of finite groups and $R^n$ or $M_n(R)$, where $R$ is an order in a central simple algebra over an algebraic number field.
Assume $\cR^\times\acts\cR$ is exact (this occurs if and only if there is no group embedding of $\Qz$ into the additive group of $\cR$, e.g., if the additive group of $\cR$ is isomorphic to $\Zz^d$ for some $d\in\Zz_{>0}$). It is straightforward to check that $\cR^\times\acts\cR$ satisfies \eqref{eqn:(FI)} and that $(\Qz\otimes_\Zz\cR)^*\acts\Qz\otimes_\Zz\cR$ is a globalization for $\cR^\times\acts\cR$ that satisfies \eqref{eqn:(JF)}.
By Corollary~\ref{cor:(FI)implies(M)}, Corollary~\ref{cor:Summary} and Theorem~\ref{thm:Amenability} imply that the following are equivalent:
\setlength{\parindent}{0cm} \setlength{\parskip}{0cm}

\begin{enumerate}[\upshape(i)]
\item $(\Qz\otimes_\Zz\cR)^*$ is amenable;
\item $\fA_r[\cR]$ is nuclear;
\item $\fA_r[\cR]$ is simple.
\end{enumerate}
If the above equivalent conditions are satisfied, then $\fA_r[\cR]$ is a UCT Kirchberg algebra by Corollary~\ref{cor:PI}. 
\eex
\setlength{\parindent}{0cm} \setlength{\parskip}{0.5cm}

\subsection{Algebraic actions by left reversible monoids}
\label{sec:leftreversible}

Recall that the monoid $S$ is said to be left reversible if $sS\cap tS\neq\emptyset$ for all $s,t\in S$. Right reversibility is defined analogously. We shall now demonstrate that our conditions from \S~\ref{sec:properties} are especially easy to check for actions by left reversible monoids. 

\subsubsection{General results for actions satisfying (PC)}
We shall call $S$-constructible subgroups of the form $\sigma_sG$ \emph{principal constructible subgroups}. This terminology comes from the ring-theoretic examples where the principal constructible subgroups are principal ideals of the ring. 
Consider the following condition on $\sigma\colon S\acts G$:
\begin{equation}
	\label{eqn:(PC)}\tag{PC}
	\text{For every }C\in\cC,\text{ there exists }s\in S \text{ such that }\sigma_sG\subseteq C.
\end{equation}
This condition means that the family of principal constructible subgroups is co-final in $\cC$.

\bprop
\label{prop:exactness4Sreversible}
If $S$ is left reversible, then $\sigma\colon S\acts G$ satisfies \eqref{eqn:(PC)}. 
\eprop
\setlength{\parindent}{0cm} \setlength{\parskip}{0cm}

\bproof
Let $C =\sigma_{s_1}^{-1} \sigma_{t_1} \dotsm \sigma_{s_n}^{-1} \sigma_{t_n}G\in\cC$. As $S$ is left reversible, there exists $s \in S$ such that $sS \subseteq s_1^{-1} t_1 \dotsm s_n^{-1} t_nS$ (see \cite[Lemma~5.6.43.]{CELY}). Hence $\sigma_sG \subseteq C$.
\eproof
\setlength{\parindent}{0cm} \setlength{\parskip}{0.5cm}

\bremark
\label{rmk:exactness4Sreversible}
If $\sigma\colon S\acts G$ satisfies \eqref{eqn:(PC)}, then $\bigcap_{C\in \cC} C=\bigcap_{s\in S}\sigma_sG$, so that $\sigma\colon S\acts G$ is exact if and only if $\bigcap_{s\in S}\sigma_sG=\{e\}$.

If $G$ is Abelian, then $\bigcap_{s\in S}\sigma_sG=\{e\}$ if and only if $\bigcup_{s\in S}\fix(\hat{\sigma}_s)$ is dense in $\widehat{G}$.
\eremark

\blemma
\label{lem:PCimpliesM}
If $\sigma\colon S\acts G$ satisfies \eqref{eqn:(PC)} (e.g., if $S$ is left reversible), then $\sigma\colon S\acts G$ satisfies \eqref{eqn:(M2)}.
\elemma
\setlength{\parindent}{0cm} \setlength{\parskip}{0cm}

\bproof
Let $C\in\cC$. Since $\sigma\colon S\acts G$ satisfies \eqref{eqn:(PC)}, there exists $s\in S$ such that $\sigma_sG\subseteq C$. Now $G=\sigma_s^{-1}C$, so \eqref{eqn:(M2)} holds.
\eproof
\setlength{\parindent}{0cm} \setlength{\parskip}{0.5cm}

\bcor
\label{cor:(JF)4Sreversible}
Assume $G$ is Abelian, and $S$ is cancellative and right reversible (i.e., left Ore). If $\sigma\colon S\acts G$ is faithful and satisfies \eqref{eqn:(PC)}, then $\sigma\colon \gp{S}\acts S^{-1}G$ satisfies \eqref{eqn:(JF)}.
\ecor
\setlength{\parindent}{0cm} \setlength{\parskip}{0cm}

\bproof
For convenience, let us write $G$ additively. Since $S$ is left Ore, by (iii) in Example~\ref{ex:InvSgp-Ore}, it suffices to show that $C\subseteq\ker(\sigma_s-\sigma_t)\implies s=t$ for all $C\in\cC$ and $s\in S$. Suppose we have $C\in\cC$ and $s,t\in S$ with $C\subseteq \ker(\sigma_s-\sigma_t)$. By assumption, $S\acts G$ satisfies \eqref{eqn:(PC)}, so there exists $r\in S$ such that $\sigma_rG\subseteq C$, so we have $\sigma_s\sigma_r(g)=\sigma_t\sigma_r(g)$ for all $g\in G$. By faithfulness, it follows that $sr = tr$, and hence $s=t$ by right cancellation.
\eproof
\setlength{\parindent}{0cm} \setlength{\parskip}{0.5cm}

\blemma
Assume that $\sigma\colon S\acts G$ has a globalization $\tilde{\sigma}\colon\mathscr{S}\acts\mathscr{G}$ and that \eqref{eqn:(PC)} holds. Then \eqref{eqn:(JF)} is satisfied if and only if $\ti{\sigma}_g \vert_G = \id$ implies $g=1$ for all $g\in \mathscr{G}\rtimes\mathscr{S}$. 
\elemma
\setlength{\parindent}{0cm} \setlength{\parskip}{0cm}

\bproof
Assume $\ti{\sigma}_g \vert_G = \id$ implies $g=1$ for all $g\in \mathscr{G}\rtimes\mathscr{S}$. Suppose $\tilde{\sigma}_g\vert_C=\id_C$ for some $C\in\cC$. Then by \eqref{eqn:(PC)}, we can find $s\in S$ with $\sigma_sG\subseteq C$. Now we have that $\ti{\sigma}_g \circ \ti{\sigma}_s$ and $\ti{\sigma}_s$ agree on $G$. Therefore our assumption implies $gs = s$ and thus $g = 1$.
\eproof
\setlength{\parindent}{0cm} \setlength{\parskip}{0.5cm}

By Corollary~\ref{cor:PI} and Corollary~\ref{cor:Summary}, we have:
\btheo
\label{thm:Sreversible}
Assume $\sigma\colon S\acts G$ is a non-automorphic algebraic action with $G$ Abelian and $S$ left and right reversible and cancellative. Then the groupoid $\cG_\sigma\cong (\gp{S}\rtimes S^{-1}G)\ltimes\bd$ is Hausdorff, minimal, and purely infinite. Moreover, if additionally $\sigma\colon S\acts G$ is exact, then $\cG_\sigma$ is topologically free, and
\setlength{\parindent}{0cm} \setlength{\parskip}{0cm}

\begin{enumerate}
\item[\upshape(i)] $C_r^*(\cG_\sigma)$ is simple and purely infinite;
\item[\upshape(ii)] The map $\fA_\sigma\to C_r^*(\cG_\sigma)$, $\Lambda_\phi\mapsto v_\phi$ from Proposition~\ref{prop:essential} is an isomorphism if and only if $\gp{S}$ is amenable. In particular, $\fA_\sigma$ is simple if and only if $\gp{S}$ is amenable.
\end{enumerate}
\setlength{\parindent}{0cm} \setlength{\parskip}{0.5cm}

\etheo

\subsubsection{Algebraic actions from commutative algebra}
Inspired by the rich class of algebraic actions of groups arising from considerations in commutative algebra, see, e.g., \cite{Sch}, we now turn to examples of algebraic actions of semigroups arising from modules over commutative rings.

Let $R$ be an infinite, commutative, unital ring with $0\neq 1$, and let $M$ be a non-zero $R$-module. For $a\in R$, let $\alpha_a^M\in\End_\Zz(M)$ be the associated endomorphism of $M$. We often omit the superscript and simply write $\alpha_a$ when the module in question is clear from context.
Let $\Reg(M):=\{a\in R : \alpha_a \text{ is injective}\}$
be the commutative monoid of $M$-regular elements. For the $R$-module $R$, $R^\times:=\Reg(R)$ is the multiplicative monoid of non-zerodivisors in $R$.
Let $S$ be a non-trivial submonoid of $R^\times\cap\Reg(M)$; since $R^\times$ is cancellative, so is $S$. Consider the action $\alpha\colon S\acts M$ given by $a\mapsto \alpha_a$.
For $M=R$ and $S=R^\times$, the concrete C*-algebra $\fA_\alpha$ associated with $\alpha \colon R^\times\acts R$ is the reduced ring C*-algebra of $R$ in the sense of \cite[Definition~7]{Li:Ring}, which is denoted by $\fA_r[R]$.

\bremark
In this setting, the canonical globalization has a particularly nice form. Let $S^{-1}R$ and $S^{-1}M$ denote the localizations of $R$ and $M$, respectively, at $S$. The canonical map $R\to S^{-1}R$ is injective since $S\subseteq R^\times$, and the canonical map $M\to S^{-1}M$ is injective since $S\subseteq\Reg(M)$. Let $\gp{S}$ denote the subgroup of $(S^{-1}R)^*$ generated by $S$. The canonical algebraic action $\gp{S}\acts S^{-1}M$ is a globalization of $S\acts M$. By Corollary~\ref{cor:(JF)4Sreversible}, $\gp{S}\acts S^{-1}M$ satisfies \eqref{eqn:(JF)} if and only if $S\acts M$ is faithful.
\setlength{\parindent}{0.5cm} \setlength{\parskip}{0cm}

For $M=R$ and $S=R^\times$, the ring $Q(R):=(R^\times)^{-1}R$ is the total quotient ring of $R$.
\eremark

If $\alpha\colon S\acts M$ is faithful, then since $S$ is Abelian, we are in the setting of Theorem~\ref{thm:Sreversible}, and $\fA_\alpha$ will be a UCT Kirchberg algebra whenever $\alpha\colon S\acts M$ is non-automorphic and exact (nuclearity comes from Theorem~\ref{thm:Amenability}). Thus, we set out to establish criteria for these conditions to be satisfied. We start with the case $M=R$.
 
\bremark
If $M=R$, then $S\acts R$ is always faithful, and $S\acts R$ is automorphic if and only if $S\subseteq R^*$.
\eremark

\bprop
\label{prop:IDexact}
Let $R$ be a Noetherian integral domain and $S\subseteq R^\times$ a submonoid. Then $S\acts R$ is exact if and only if $S$ contains a non-unit.
\eprop
\setlength{\parindent}{0cm} \setlength{\parskip}{0cm}

\bproof
If $S\subseteq R^*$, then the action $S\acts R$ is by automorphisms and is thus not exact.
Suppose there exists $a\in S\setminus R^*$ such that $(a)\subsetneq R$. Since $R$ is Noetherian by assumption, \cite[Corollary~10.18]{AM} says that $\bigcap_{n=1}^\infty(a)^n=\{0\}$. Since each $(a)^n=(a^n)=a^nR$ lies in $\cC_{S\acts R}$ for every $n$, exactness follows.
\eproof
\setlength{\parindent}{0cm} \setlength{\parskip}{0.5cm}

Applying Theorem~\ref{thm:Amenability} and Theorem~\ref{thm:Sreversible} with the above observations gives:

\bcor
\label{cor:IDKirchberg}
\begin{enumerate}[\upshape(i)]
 \item If $\bigcap_{a\in R^\times}aR=\{0\}$, then the reduced ring C*-algebra $\fA_r[R]$ associated with $R^\times\acts R$ is a UCT Kirchberg algebra.
\item If $R$ is a Noetherian integral domain, and $S\subseteq R^\times$ is a submonoid containing a non-unit, then the C*-algebra $\fA_\alpha$ associated with $\alpha\colon S\acts R$ is a UCT Kirchberg algebra.
\end{enumerate}
\ecor
\bremark
\label{rmk:IDKirchberg}
Corollary~\ref{cor:IDKirchberg}(i) applies to a larger class of rings than those treated by the results in \cite[\S~5.3]{Li:Ring}; indeed, in order to apply the results of \cite{Li:Ring} to the reduced ring C*-algebra $\fA_r[R]$, one needs $\bigcap_{a\in R^\times}aR=\{0\}$ and also condition (**) from \cite[\S~5.3]{Li:Ring}. 
\setlength{\parindent}{0.5cm} \setlength{\parskip}{0cm}

If $R$ is an integral domain that is not a field, then $\bigcap_{a\in R^\times}aR=\{0\}$ (see the proof of \cite[Corollary~9]{Li:Ring}), so Corollary~\ref{cor:IDKirchberg}(i) applies to all integral domains that are not fields. This special case is also covered by \cite[Corollaries~8~\&~9]{Li:Ring} or \cite[Corollary~8.4]{LaSe}.

Let $R$ be the ring of integers in an algebraic number field. The boundary quotients $\partial C_\lambda^*(R\rtimes R_{\m,\Gamma})$ associated with the action of a congruence monoid $R_{\m,\Gamma}$ on $R$ are UCT Kirchberg algebras by \cite[\S~8]{Bru1} and \cite[Theorem~3.1]{BruLi}. Corollary~\ref{cor:IDKirchberg}(ii) generalizes and explains this. 
\eremark
\setlength{\parindent}{0cm} \setlength{\parskip}{0.5cm}

Let us now turn to more general modules. Let $\Ann(M):=\{a\in R : a.x=0\text{ for all }x\in M\}\unlhd R$ denote the annihilator ideal of $M$. It is easy to see that $S\acts M$ is faithful if and only if $(S-S)\cap \Ann(M)=\{0\}$, where $S-S:=\{a-b \in R:a,b\in S\}$.

\bex[Prime actions]
\label{ex:primeactions}
Let $\p$ be a prime ideal of $R$, and consider $M=R/\p$ as an $R$-module. Then $\Reg_R(R/\p)=R\setminus\p$, so $S\acts R/\p$ acts by injective endomorphisms if and only if $S\subseteq R\setminus\p$. Such actions are called \emph{prime actions}. Since $\Ann(R/\p)=\p$, the action $S\acts R/\p$ is faithful if and only if $(S-S)\cap\p=\{0\}$. The action $S\acts R/\p$ is non-automorphic if and only if $S$ contains an element $a$ such that $a+\p\in (R/\p)^\times$ is a non-unit. In particular, this means that $\p$ cannot be a maximal ideal.
If $R/\p$ is Noetherian (e.g., if $R$ is Noetherian), then by Proposition~\ref{prop:IDexact}, $S\acts R/\p$ is exact if and only if there exists $a\in S$ such that $a+\p\in (R/\p)^\times$ is a non-unit.
\eex

A prime $\p$ of $R$ is said to be \emph{associated with $M$} if $\p=\ann_R(x)$ for some $x\in M$. Let $\Asc(M)$ denote the (possibly empty) set of primes associated with $M$. If $\p$ is a prime of $R$, then $\p\in\Asc(R)$ if and only if there is an embedding of $R$-modules $R/\p\hookrightarrow M$. Note that $\Ann(M)\subseteq\p$ for every $\p\in\Asc(M)$. 
If $R$ is Noetherian, then $\bigcup_{\p\in\Asc(M)}\p=R\setminus\Reg(M)$ (see, e.g., \cite[Theorem~6.1, p.38]{Mats2}). If $R$ is Noetherian and $M$ is finitely generated, then $\Asc(M)$ is finite and non-empty (see, e.g., \cite[Theorem~3.1]{Eis}).

\blemma
\label{lem:Noetherfaithful}
Assume $\Asc(M)$ is non-empty and that $S\subseteq R\setminus\bigcup_{\p\in\Asc(M)}\p$. Then $S\acts M$ is faithful if and only if there exists a prime $\p\in\Asc(M)$ such that the canonical prime action $S\acts R/\p$ is faithful.
\elemma
\setlength{\parindent}{0cm} \setlength{\parskip}{0cm}

\bproof
If $S\acts M$ is not faithful, then there exist $a,b\in S$ such that $a\neq b$ and $a.x=b.x$ for every $x\in M$, i.e., $a-b\in\Ann(M)$. Since $\Ann(M)\subseteq \p$ for every $\p\in\Asc(M)$, we see that $S\acts R/\p$ is not faithful for every $\p\in\Asc(M)$. For every $\p\in\Asc(M)$, we have an embedding of $R$-modules $R/\p\hookrightarrow M$. Hence, if $S\acts R/\p$ is faithful for some $\p\in\Asc(M)$, then $S\acts M$ must be faithful.
\eproof
\setlength{\parindent}{0cm} \setlength{\parskip}{0.5cm}

\bprop
\label{prop:Noetherexact}
Let $R$ be a Noetherian ring, $M$ a finitely generated $R$-module, and $S\subseteq R\setminus\bigcup_{\p\in\Asc(M)}\p$ a submonoid. The following are equivalent: 
\setlength{\parindent}{0cm} \setlength{\parskip}{0cm}

\begin{enumerate}[\upshape(i)]
	\item there exists $a\in S$ such that $a^\Nz\acts M$ is exact;
	\item $S\acts M$ is exact;
	\item $S\acts R/\p$ is exact for every $\p\in\Asc(M)$.
\end{enumerate}
\eprop
\setlength{\parindent}{0cm} \setlength{\parskip}{0cm}

\bproof
(i)$\Rightarrow$(ii) is obvious. (ii)$\Rightarrow$(iii): Assume $S\acts M$ is exact. Let $\p\in\Asc(M)$, so that there is an $R$-module embedding $\iota\colon R/\p\hookrightarrow M$. By Proposition~\ref{prop:exactness4Sreversible}, $S\acts M$ and every $S\acts R/\p$ satisfy \eqref{eqn:(PC)}. Since $\bigcap_{s\in S}\alpha^{\iota(R/\p)}_s(\iota(R/\p))\subseteq \bigcap_{s\in S}\alpha^M_s(M)$, exactness of $S\acts R/\p$ follows from exactness of $S\acts M$ (cf. Remark~\ref{rmk:exactness4Sreversible}).
\setlength{\parindent}{0.5cm} \setlength{\parskip}{0cm}

(iii)$\Rightarrow$(i): Assume $S\acts R/\p$ is exact for every $\p\in\Asc(M)$. By Example~\ref{ex:primeactions}, this is equivalent to the following statement: For each $\p\in\Asc(M)$, there exists $a_\p\in S$ such that $(a_\p)+\p\subsetneq R$. Let $a:=\prod_{\p\in\Asc(M)}a_\p\in S$. Since $R$ is Noetherian and $M$ finitely generated, \cite[Theorem~10.17]{AM} implies
\begin{equation}
	\label{eqn:Krull}
\bigcap_{n=0}^\infty (a)^n.M=\{x\in M : \text{there exists } r\in R \text{ such that }(1-ar).x=0\}.
\end{equation}
Suppose $x\in \bigcap_{n=0}^\infty a^n.M$ is non-zero. Since $a^n.M\subseteq (a)^n.M$, we see from \eqref{eqn:Krull} that there exists $r\in R$ such that $(1-ar).x=0$, i.e., $1-ar$ is not $M$-regular and thus lies in $\bigcup_{\p\in\Asc(M)}\p=R\setminus\Reg(M)$. Hence, $1-ar\in\p$ for some $\p\in\Asc(M)$, i.e., $a$ is invertible modulo $\p$, so that $(a)+\p=R$. But we know that $(a)+\p\subseteq (a_\p)+\p\subsetneq R$, so this is a contradiction. Hence, $x=0$.
\eproof
\setlength{\parindent}{0cm} \setlength{\parskip}{0.5cm}

\begin{example}[A result of Krzy{\.z}ewski]
\label{ex:Kryz}
Fix $n\geq 1$. Let $A\in \M_n(\Zz)$ be a matrix with non-zero determinant, and let $\sigma_A$ be the associated endomorphism of $\Zz^n$. Let $\chi_A(u)$ be the characteristic polynomial of $A$. The main result in \cite{Krz} says that $\sigma_A\colon\Nz\acts \Zz^n$ is exact if and only if $\chi_A(u)$ not divisible by any polynomial with constant term $\pm 1$ (i.e., $\chi_A(u)$ has no unimodular factors). We demonstrate that this characterization follows from Proposition~\ref{prop:Noetherexact}.

View $\Zz^n$ as a $\Zz[u]$-module via $f(u).\bm{n}=f(A)\bm{n}$, and similarly view $\Qz^n$ as a $\Qz[u]$-module. Let $m_A(u)$ be the minimal polynomial of $A$, i.e., $m_A(u)$ is the unique monic generator of the ideal $\Ann_{\Qz[u]}(\Qz^n)\unlhd \Qz[u]$. Since $A$ has integer entries, $m_A(u)\in\Zz[u]$. Let $m_A(u)=\prod_{i=1}^rf_i(u)^{k_i}$ be the factorization of $m_A(u)$ into powers of irreducible elements, which are defined up to multiplication by $\pm 1$ (here, we are using that $\Zz[u]$ is a unique factorization domain with unit group $\{\pm 1\}$). 
\setlength{\parindent}{0.5cm} \setlength{\parskip}{0cm}

By Proposition~\ref{prop:IDexact} applied to $S = u^{\Nz} \cong \Nz$, the canonical $\Nz$-action $\Nz\acts \Zz[u]/(f_i(u))$ generated by $g(u)+(f_i(u))\mapsto ug(u)+(f_i(u))$ is exact if and only if $u+(f_i(u))$ is a non-unit, which is equivalent to $f_i(0)\notin\{\pm 1\}$. Thus, in order to deduce  Krzy{\.z}ewski's characterization of exactness from Proposition~\ref{prop:Noetherexact} it suffices to prove the following lemma.

\blemma
In the situation of Example~\ref{ex:Kryz}, we have $\Asc_{\Zz[u]}(\Zz^n)=\{(f_i(u)) : 1\leq i\leq r\}$.
\elemma
\setlength{\parindent}{0cm} \setlength{\parskip}{0cm}

\bproof
If $\p\in \Asc_{\Zz[u]}(\Zz^n)$, then $\Zz[u]/\p$ embeds as a $\Zz[u]$-submodule of $\Zz^n$; in particular, $\Zz[u]/\p$ is torsion-free as an additive group, so that $\p\cap\Zz=(0)$. Hence, taking localizations with respect to $\Zz^\times$ and applying \cite[Theorem~3.1(c)]{Eis} gives us
\begin{equation}
	\label{eqn:AscQ[u]}
\Asc_{\Qz[u]}(\Qz^n)=\{\p_\Qz : \p\in\Asc_{\Zz[u]}(\Zz^n)\},
\end{equation}
where $\p_\Qz$ denotes the prime ideal of $\Qz[u]$ generated by $\p$. Let $\q\in\Asc_{\Qz[u]}(\Qz^n)$. Since $\Qz[u]$ is a principal ideal domain, we can write $\q=(p(u))_\Qz$ for some irreducible polynomial $p(u)\in\Zz[u]$. Here, we write $(p(u))_\Qz$ for the ideal of $\Qz[u]$ generated by $p(u)$. Since $(p(u))_\Qz\supseteq\Ann_{\Qz[u]}(\Qz^n)=(m_A(u))_\Qz$, we have $p(u)\mid m_A(u)$ in $\Qz[u]$, so that $(p(u))_\Qz=(f_i(u))_\Qz$ for some $1\leq i\leq r$. Thus, we have $\Asc_{\Qz[u]}(\Qz^n)\subseteq \{(f_i(u))_\Qz : 1\leq i\leq r\}$.
\setlength{\parindent}{0cm} \setlength{\parskip}{0.5cm}

Fix $1\leq i\leq r$, and let $\bm{x}\in (m_A(u)/f_i(u)).\Qz^n\setminus\{\bm{0}\}$ (such an $\bm{x}$ exists by definition of $m_A$). Then $f_i(u).\bm{x}=\bm{0}$, so $(f_i(u))_\Qz\subseteq\Ann_{\Qz[u]}(\bm{x})$. Since $1.\bm{x}\neq \bm{0}$, $\Ann_{\Qz[u]}(\bm{x})$ is a proper ideal of $\Qz[u]$; since $(f_i(u))_\Qz$ is a prime ideal and $\Qz[u]$ has Krull dimension $1$, we must have $(f_i(u))_\Qz=\Ann_{\Qz[u]}(\bm{x})$. Hence, $(f_i(u))_\Qz\in\Asc_{\Qz[u]}(\Qz^n)$.

Therefore, using \eqref{eqn:AscQ[u]}, we have $\{\p_\Qz : \p\in\Asc_{\Zz[u]}(\Zz^n)\}=\{(f_i(u))_\Qz : 1\leq i\leq r\}$. It remains to observe that if $f(u)\in\Zz[u]$ is a monic polynomial, then $(f(u))_\Qz\cap\Zz[u]=f\Zz[u]$ by Gauss's lemma.
\eproof
\end{example}

\bex[Algebraic $\Nz^d$-actions]
\label{ex:AlgNdAct}
Fix $d\in\Zz_{>0}$, and let $R_d^+:=\Zz[u_1,...,u_d]$ be the polynomial ring with integer coefficients in the $d$ commuting variables $u_1,...,u_d$. 
For $\bm{n}=(n_1,...,n_d)\in\Nz^d$, we let $u^{\bm{n}}:=u_1^{n_1}\cdots u_d^{n_d}$. Given $f\in R_d^+$, we can write $f=\sum_{\bm{n}\in\Nz^d}f_{\bm{n}}u^{\bm{n}}$, where $f_{\bm{n}}\in\Zz$ is zero for all but finitely many $\bm{n}$.
Given any algebraic $\Nz^d$-action $\Nz^d\acts M$, where $M$ is an Abelian group, $M$ naturally becomes a module over $R_d^+$ via $f.x:=\sum_{\bm{n}}f_{\bm{n}}u^{\bm{n}}.x$. Note that $R_d^+$ is Noetherian (see, e.g., \cite[Theorem~7.5]{AM}).
\eex

\bprop
\label{prop:Ndactions}
Let $M$ be a finitely generated module over $\Zz[u_1,...,u_d]$, and assume that $u_i\notin\p$ for every $1\leq i\leq d$ and every $\p\in\Asc(M)$, so that we get an algebraic $\Nz^d$-action $\Nz^d\acts M$ as in Example~\ref{ex:AlgNdAct}. Then
\setlength{\parindent}{0cm} \setlength{\parskip}{0cm}

\begin{enumerate}[\upshape(i)]
	\item $\Nz^d\acts M$ is faithful if and only if there exists $\p\in\Asc(M)$ such that 
	\begin{equation}
		\{u^{\bm{n}}-u^{\bm{m}} : \bm{n},\bm{m}\in\Nz^d\}\cap\p=\{0\};
	\end{equation}
\item Let $\pi_j: \: R_d^+ \onto \Zz[u_1, \dotsc, u_{j-1}, u_{j+1}, \dotsc, u_d]$ be the canonical projection. $\Nz^d\acts M$ is exact if and only if for all $\p\in\Asc(M)$ there exists $1 \leq j \leq d$ such that $1 \notin \pi_j(\p)$. In particular, $\Nz^d\acts M$ is exact if we can find $f_1,...,f_m\in \Zz[u_1,...,u_d]$ and $\bm{z}\in\Zz^d$ with $\bm{z}_j=0$ for some $1\leq j\leq d$ such that $\p=(f_1,...,f_m)$ and $\gcd(f_1(\bm{z}),...,f_m(\bm{z}))\neq 1$. 
\end{enumerate}
\eprop
\setlength{\parindent}{0cm} \setlength{\parskip}{0cm}

\bproof
(i): By Lemma~\ref{lem:Noetherfaithful}, $\Nz^d\acts M$ is faithful if and only if there exists a prime $\p\in\Asc(M)$ such that $\Nz^d\acts R_d^+/\p$ is faithful, which happens if and only if $\{u^{\bm{n}}-u^{\bm{m}} : \bm{n},\bm{m}\in\Nz^d\}\cap\p=\{0\}$. 
\setlength{\parindent}{0.5cm} \setlength{\parskip}{0cm}

(ii): By Example~\ref{ex:primeactions} and Proposition~\ref{prop:Noetherexact}, $\Nz^d\acts M$ is exact if and only if the image of the canonical map $\mon{u_1,...,u_d}\to R_d^+/\p$ contains a non-unit for every $\p\in\Asc(M)$. The latter holds if and only if there exists $1 \leq j \leq d$ such that $u^{\bm{e}_j} + \p$ is not a unit in $R_d^+/\p$, where $\gekl{\bm{e}_j}$ denotes the canonical generators of $\Nz^d$. Now $u^{\bm{e}_j} + \p$ is not a unit in $R_d^+/\p$ if and only if $1 \notin (u^{\bm{e}_j}) + \p$ if and only if $1 \notin \pi_j(\p)$.
\eproof
\setlength{\parindent}{0cm} \setlength{\parskip}{0.5cm}

\subsection{Shifts over semigroups}
\label{sec:shifts}

Throughout this section, let $S$ be a left cancellative monoid with $S\neq S^*$ and $\Sigma$ any non-trivial group. As before, we denote the identity of $\Sigma$ by $e$. Let us recall the definition of shifts over semigroups from Example~\ref{ex:genshift}.

\bdefin
\label{def:genshift}
The \emph{full $S$-shift over $\Sigma$} is the algebraic $S$-action
\[
\sigma\colon S\acts \textstyle{\bigoplus_S}\Sigma,\quad \sigma_s(x)_t:=\begin{cases}
x_{s^{-1}t} & \text{ if } t\in sS,\\
e & \text{ if } t\notin sS.
\end{cases}
\]
\edefin

For $t\in S$ and $x\in\Sigma$, we let $x\eps_t\in \textstyle{\bigoplus_S}\Sigma$ be the element defined by $(x\eps_t)_t=x$ and $(x\eps_t)_s=e$ for $s\neq t$. Then $\sigma_s(x\eps_t)=x\eps_{st}$, and we can write every element of $a=(a_t)_t\in\textstyle{\bigoplus_S}\Sigma$ as $a=\prod_{t\in S}a_t\eps_t$. Thus, we see that the algebraic $S$-action $S\acts \textstyle{\bigoplus_S}\Sigma$ is faithful because the left translation action $S\acts S$ is faithful.

For $a\in \textstyle{\bigoplus_S}\Sigma$, let $\supp(a):=\{s\in S : a_s\neq e\}$. Given $X\subseteq S$, we identify $\textstyle{\bigoplus_X}\Sigma$ with the subgroup $\{a\in \textstyle{\bigoplus_S}\Sigma : \supp(a)\subseteq X\}$.
\bprop
\label{prop:Cforshift}
We have $\cC_{S\acts\textstyle{\bigoplus_S}\Sigma}=\{\textstyle{\bigoplus_X}\Sigma: X\in\cJ_S\}$. Moreover, the map $\cJ_S\to \cC_{S\acts\textstyle{\bigoplus_S}\Sigma}$ given by $X\mapsto \textstyle{\bigoplus_X}\Sigma$ is an isomorphism of semilattices.
\eprop
\setlength{\parindent}{0cm} \setlength{\parskip}{0cm}

Here, as in \S~\ref{sec:semigpC*}, $\cJ_S$ denotes the semilattice of constructible right ideals of $S$.
\bproof
Given $a=\prod_{t\in S}a_t\eps_t\in\textstyle{\bigoplus_S}\Sigma$, $X\subseteq S$, and $s\in S$, we have $\supp(a)\subseteq X$ if and only if $\supp(\sigma_s(a))\subseteq sX$, and we have $\supp(a)\subseteq s^{-1}X$ if and only if $\supp(\sigma_s(a))\subseteq X$. 
Hence, $\sigma_s(\textstyle{\bigoplus_X}\Sigma)=\textstyle{\bigoplus_{sX}}\Sigma$ and $\sigma_s^{-1}(\textstyle{\bigoplus_X}\Sigma)=\textstyle{\bigoplus_{s^{-1}X}}\Sigma$, which is enough for the first claim.
\setlength{\parindent}{0.5cm} \setlength{\parskip}{0cm}

For the second claim, note that $\textstyle{\bigoplus_X}\Sigma\cap\textstyle{\bigoplus_Y}\Sigma=\textstyle{\bigoplus_{X\cap Y}}\Sigma$ for all $X,Y\in\cJ_S$.
\eproof
\setlength{\parindent}{0cm} \setlength{\parskip}{0.5cm}

We immediately obtain the following:
\bcor
\label{cor:genshifts}
\begin{enumerate}[\upshape(i)]
    \item $S\acts\textstyle{\bigoplus_S}\Sigma$ is exact;
    \item $S\acts\textstyle{\bigoplus_S}\Sigma$ satisfies \eqref{eqn:(M2)} from Theorem~\ref{thm:minimal} if and only if $S$ is left reversible;
    \item $S\acts\textstyle{\bigoplus_S}\Sigma$ satisfies \eqref{eqn:(FI)} from \S~\ref{sec:FI} if and only if $\#\Sigma<\infty$ and $\#(S\setminus sS)<\infty$ for all $s\in S$.
\end{enumerate}
\ecor
\setlength{\parindent}{0cm} \setlength{\parskip}{0cm}

\bproof
(i): Suppose we have $t\in\bigcap_{s\in S}sS$. Then there exists $t'\in S$ such that $t=t^2 t'$, which by left cancellation implies $1=tt'$. Since $S$ is left cancellative, it follows that $t\in S^*$. But $\bigcap_{s\in S}sS$ is a proper right ideal of $S$ because $S$ contains a non-invertible element, so we have a contradiction. It follows that $\bigcap_{X\in\cJ_{S}^\times}X\subseteq\bigcap_{s\in S}sS=\emptyset$. Thus, $S\acts\bigoplus_S\Sigma$ is exact by Proposition~\ref{prop:Cforshift}.
\setlength{\parindent}{0cm} \setlength{\parskip}{0.5cm}

(ii): If $S$ is left reversible, then \eqref{eqn:(M2)} holds by Lemma~\ref{lem:PCimpliesM}.
Conversely, if $S$ is not left reversible, $\emptyset\in \cJ_S$ by \cite[Lemma~5.6.43.]{CELY}, so that $\{e\}$ is constructible by Proposition~\ref{prop:Cforshift}. This implies that \eqref{eqn:(M2)} does not hold.

(iii): For $s\in S$, $\textstyle{\bigoplus_S}\Sigma / \textstyle{\bigoplus_{sS}}\Sigma\cong\textstyle{\bigoplus_{S\setminus sS}}\Sigma$, which is finite if and only if $\#\Sigma<\infty$ and $\#(S\setminus sS)<\infty$.
\eproof
\setlength{\parindent}{0cm} \setlength{\parskip}{0.5cm}

\blemma
Assume $S$ is left reversible.
\setlength{\parindent}{0cm} \setlength{\parskip}{0cm}

\begin{enumerate}[\upshape(i)]
    \item If $S$ is right cancellative, then $S\acts\textstyle{\bigoplus_S}\Sigma$ satisfies \eqref{eqn:(H)} from \S~\ref{ss:Hd}.
    \item If $S$ embeds in a group $\msc{S}$, then the globalization $\ti{\sigma}\colon\msc{S}\acts\textstyle{\bigoplus_\msc{S}}\Sigma$ satisfies \eqref{eqn:(JF)} from \S~\ref{ss:PartialGlobal}.
\end{enumerate}
\elemma
\setlength{\parindent}{0cm} \setlength{\parskip}{0cm}

\bproof
(i): Suppose we are given a constructible subgroup $\textstyle{\bigoplus_X}\Sigma$ and $s\in S$ with $\textstyle{\bigoplus_X}\Sigma\subseteq\fix(\sigma_s)$. Then $x\eps_{st}=\ti{\sigma}_s(x\eps_t)=x\eps_t$ for all $x\in \Sigma$ and all $t\in X$, so that $st=t$ for all $t\in X$. Since $S$ is left reversible, $X\neq \emptyset$. Thus, since $S$ is right cancellative, we have $s=1$. It follows that \eqref{eqn:(H)} is satisfied. The proof of (ii) is similar.  
\eproof
\setlength{\parindent}{0cm} \setlength{\parskip}{0.5cm}

Using Theorems~\ref{thm:Amenability} and Corollaries~\ref{cor:PI}, \ref{cor:Summary}, the following is now an immediate consequence.
\bcor
\label{cor:FullShiftKirchberg}
Let $S$ be a countable, left reversible submonoid of a group $\msc{S}$ such that $S\neq S^*$, and let $\Sigma$ be a countable non-trivial group. Assume $\gp{S}=\msc{S}$ and that $\msc{S}$ and $\Sigma$ are both amenable. Then the concrete C*-algebra $\fA_\sigma$ associated with $\sigma: \: S\acts\textstyle{\bigoplus_S}\Sigma$ is a UCT Kirchberg algebra.
\ecor

\bex
\label{ex:exotic}
Consider $\sigma: \: S\acts\textstyle{\bigoplus_S}\Sigma$, where $S$ is a countable, left reversible submonoid of a group $\msc{S}$ such that $S\neq S^*$, and $\Sigma$ is a countable group. Assume $\gp{S}=\msc{S}$ and that $\msc{S}$ and $\Sigma$ are both non-amenable. Then $\fA_\sigma$ is an exotic groupoid C*-algebra by Corollary~\ref{cor:exotic}.
\eex

\bremark[Dual subshifts]
Now let us consider algebraic actions of the form $\sigma \colon S \acts G \subseteq \bigoplus_S \Sigma$ given by restricting the full shift $S \acts \bigoplus_S \Sigma$ (which we will also denote by $\sigma$). The general recipe to construct invariant subgroups $G \subseteq \bigoplus_S \Sigma$ is to start with an arbitrary subgroup $G_0 \subseteq \bigoplus_S \Sigma$ and consider the smallest subgroup of $\bigoplus_S \Sigma$ generated by $\bigcup_{s \in S} \sigma_s(G_0)$,
$$
 G \defeq \Big\langle \bigcup_{s \in S} \sigma_s(G_0) \Big\rangle \subseteq \bigoplus_S \Sigma.
$$
We call such algebraic actions dual subshifts because we are forming subshifts on the dual side (if $\Sigma$ is Abelian, then what we call full shifts are really dual actions of classical full shifts). Analogous arguments as for Corollary~\ref{cor:genshifts}~(i), (ii) show that all these dual subshifts are exact and satisfy \eqref{eqn:(M2)} if and only if $S$ is left reversible.

Let us now assume that our semigroup $S$ is a subsemigroup of a group $\mathscr{S}$. Then we could take an arbitrary subgroup $G_0$ of $\bigoplus_{\mathscr{S}} \Sigma$, and set $G \defeq (\bigoplus_S \Sigma) \cap \ti{\sigma}_{\mathscr{S}}(G_0)$, where $\ti{\sigma}\colon\msc{S}\acts \bigoplus_\msc{S} \Sigma$ is the canonical globalization. In this case, we automatically have an enveloping action (the restriction of $\ti{\sigma}$ to $\ti{\sigma}_{\msc{S}}(G_0) = \bigcup_{g \in \msc{S}} \ti{\sigma}_g(G_0)$). Hausdorffness of our groupoid $\cG_\sigma$ is then related to the zero-divisor conjecture for group rings: If $\msc{S}$ is torsion-free and $\Sigma$ is a subgroup of the additive group of a field $\Kz$, and the zero-divisor conjecture is true for $\Kz[\msc{S}] = \bigoplus_{\msc{S}} \Kz$, then an equation of the form $(1-s)a = 0$ in $\bigoplus_\msc{S} \Sigma \subseteq \Kz[\msc{S}]$ for $1 \neq s \in \msc{S}$ implies $a = 0$. Therefore, this would then imply condition \eqref{eqn:(JF)} for the enveloping action and thus that $\cG_\sigma$ is Hausdorff.
\eremark

\bremark[Subgroup shifts]
Let us now consider algebraic actions of the form $\sigma \colon S \acts G \defeq (\bigoplus_S \Sigma)/ \msc{I}$ induced by the full shift $S \acts \bigoplus_S \Sigma$ (also denoted by $\sigma$), where $\msc{I}$ is an invariant subgroup of $\bigoplus_S \Sigma$. To ensure injectivity of $\sigma_s$ for all $s \in S$, we need $\sigma_s^{-1}(\msc{I}) = \msc{I}$ for all $s \in S$. Here is a recipe to construct such $\msc{I}$: Assume that our semigroup $S$ is a subsemigroup of a group $\mathscr{S}$. Take a left ideal $\msc{J}$ of $\bigoplus_\msc{S} \Sigma$, and set $\msc{I} \defeq (\bigoplus_S \Sigma) \cap \msc{J}$. In this case, we automatically have an enveloping action (the natural action $\msc{S} \acts (\bigoplus_\msc{S} \Sigma) /\msc{J}$). We call such algebraic actions subgroup shifts because in case $\Sigma$ is Abelian, the dual actions will be subgroup shifts in the classical sense.

In this setting, exactness is satisfied if for all finite subsets $F \subseteq S$ there exists a finite subset $F' \subseteq S$ such that for all $a \in \bigoplus_S \Sigma$ with $\supp(a) \subseteq F$, $a \vert_{F'} \in \pi_{F'}(\msc{I})$ implies that $a \in \msc{I}$, where $\pi_{F'}$ is the canonical projection from $\bigoplus_S \Sigma$ onto $\bigoplus_{F'} \Sigma$, and that for all finite subsets $F' \subseteq S$, there exists $t \in S$ with $F' \cap (tS) = \emptyset$. Indeed, suppose $a \in \bigoplus_S \Sigma$ satisfies $a \in \bigcap_{s \in S} \sigma_s(\bigoplus_S \Sigma)$ mod $\msc{I}$. Let $\supp(a) \subseteq F$ for a finite subset $F \cap S$ and let $F'$ be a finite subset for $F$ as above. Choose $t \in S$  with $F' \cap (tS) = \emptyset$. Since $a \in \sigma_t(\bigoplus_S \Sigma)$ mod $\msc{I}$, we can write $a = \sigma_t(b) + c$ for some $b \in \bigoplus_S \Sigma$ and $c \in \msc{I}$. Now $F' \cap (tS) = \emptyset$ implies that $\supp(\sigma_t(b)) \cap F' = \emptyset$, so that $(\sigma_t(b)) \vert_{F'} = e$. Hence $a \vert_{F'} \equiv c \vert_{F'}$. It follows by assumption that $a \in \msc{I}$. Hence exactness holds.
\eremark

\subsubsection{K-theory formulas and classification results}
\label{sss:KTheoryClass}

Let us now discuss K-theory. We first derive a general K-theory formula for inverse semigroup C*-algebras in our setting, i.e., for the (reduced) Toeplitz-type C*-algebra attached to $\sigma\colon S\acts G$ given by the reduced C*-algebra $C_\lambda^*(I)$ of the inverse semigroup $I$. For such C*-algebras, there are powerful tools for computing K-theory. Given $C\in\cC$, we consider the group $I_C:=\{\phi\in I: \dom(\phi)=\im(\phi)=C\}$. Applying \cite[Theorem~1.1]{Li:Ktheory}, we obtain the following K-theory formula:

\btheo
\label{thm:Ktheory}
Assume that $S$ and $G$ are countable and that $I$ admits an idempotent pure partial homomorphism to a group that satisfies the Baum--Connes conjecture with coefficients. Then,
\begin{equation}
    \label{eqn:ktheory}
    K_*(C_\lambda^*(I))\cong \bigoplus_{[C]\in I^e\backslash \cC}K_*(C_\lambda^*(I_C)),
\end{equation}
where $I^e\backslash \cC$ is the set of orbits for the action of $I^e$ on $\cC$.
\etheo

The following condition will allow us to pass from $C_\lambda^*(I)$ to $C_r^*(\cG_\sigma)$.
\bdefin
We say that $\sigma\colon S\acts G$ has the \emph{infinite index property} if 
\begin{equation}
	\label{eqn:(II)}\tag{II}
	\# (C/D) =\infty\quad \text{for all } C,D\in\cC \text{ with } D\lneq C.
\end{equation}
\edefin

\bcor
If $\sigma\colon S\acts G$ satisfies \eqref{eqn:(II)}, then $\bd=\widehat{\cE}$ and $C_r^*(\cG_\sigma)=C_\lambda^*(I)$. In this case, if the conditions in Theorem~\ref{thm:Ktheory} are satisfied, then 
$K_*(C_r^*(\cG_\sigma))\cong \bigoplus_{[C]\in I^e\backslash \cC}K_*(C_\lambda^*(I_C))$.
\ecor
\setlength{\parindent}{0cm} \setlength{\parskip}{0cm}

\bproof
If $\sigma\colon S\acts G$ satisfies \eqref{eqn:(II)}, then every character on $\cE$ is tight by Lemma~\ref{lem:CharTightChar}. Our claims follow immediately.
\eproof
\setlength{\parindent}{0cm} \setlength{\parskip}{0.5cm}

\bex
Assume $S$ embeds into a countable group $\msc{S}$.
Let $S\acts \bigoplus_S \Sigma$ be the full $S$-shift over a non-trivial group $\Sigma$ such that $\left(\bigoplus_{\msc{S}}\Sigma\right)\rtimes\msc{S}$ satisfies the Baum--Connes conjecture with coefficients. If either $\#\Sigma=\infty$ or $\#(X\setminus Y)=\infty$ for all $X,Y\in\cJ_S\reg$ with $Y\subsetneq X$, then $S\acts \bigoplus_S \Sigma$ satisfies \eqref{eqn:(II)} by Proposition~\ref{prop:Cforshift}, so that Theorem~\ref{thm:Ktheory} gives us
\[
K_*(C_r^*(\cG_\sigma))\cong\bigoplus_{[X]\in \msc{S}\backslash\cJ_S^\times}K_*(C_\lambda^*\left(\left(\textstyle{\bigoplus_X}\Sigma\right)\rtimes \msc{S}_X\right)),
\]
where $\msc{S}_X=\{\gamma\in \msc{S}: \gamma X=X\}$. If $\Sigma$ is Abelian with $\#\Sigma<\infty$, then $K_*(C_\lambda^*\left(\left(\textstyle{\bigoplus_{s\in X}}\Sigma\right)\rtimes \msc{S}_X\right))$ can be explicitly computed using \cite[Theorem~1.1]{Li:Lamp}.

If, in addition, $S$ is right LCM (i.e., $\cJ_S^\times = \menge{sS}{s \in S}$) and $S^* = \gekl{1}$, then our K-theory formula simplifies to
\begin{equation}
\label{e:KThLCMNoUnits}
K_*(C_r^*(\cG_\sigma))\cong K_*(C_\lambda^*\left(\textstyle{\bigoplus_S}\Sigma\right)),
\end{equation}
\eex

Note that if $S$ is right LCM, then the condition that $\#(X\setminus Y)=\infty$ for all $X,Y\in\cJ_S\reg$ with $Y\subsetneq X$ is equivalent to $\#(S \setminus sS) = \infty$ for all $s \in S \setminus S^*$. Let us present two example classes where the latter condition holds. 
\bex
Assume that $S$ is right LCM and $s \in S \setminus S^*$ satisfies $\#(S \setminus sS) < \infty$, say $S \setminus sS = \gekl{r_j}$. Further suppose that $S$ is right Noetherian, i.e., we cannot find an infinite chain of the form $\dotso \supsetneq Ss_3 \supsetneq Ss_2 \supsetneq Ss_1$ for $s_1, s_2, s_3, \dotsc \in S$. Then every element of $S$ is of the form $s^n r_j$ for some non-negative integer $n$ and $r_j \in S \setminus sS$. Indeed, since $S$ is right Noetherian, for every $x \in S$ there exists a maximal non-negative integer $n$ such that $x = s^n y$ with $y \notin sS$, so that $y \in \gekl{r_j}$. In particular, if $S$ is Abelian and cancellative, then this would imply that the enveloping group $\msc{S} = S^{-1} S$ is virtually Abelian. 
\eex
\bex
Now suppose that $S$ is given by generators and relations, i.e., $S = \spkl{\Gamma \vert R}^+$ as in \cite[\S~2.1.1]{LOS}, with the same standing assumptions as in \cite[\S~2.1.2]{LOS}. Assume that for every generator $\gamma \in \Gamma$ there exists an infinite word $w$ in $\Gamma$ not starting with $\gamma$ such that no relator appears as a finite subword of $w$ (in particular, this implies that $\# \Gamma > 1$ and no generator is a relator). If $w_l$ denotes the finite subword of $w$ consisting of the first $l$ letters, then our assumptions imply that $w_l \notin \gamma S$ as well as $w_l \neq w_m$ whenever $l \neq m$, so that $\#(S \setminus \gamma S) = \infty$. It then follows that $\#(S \setminus s S) = \infty$ for all $s \in S \setminus S^*$.
\eex

Let us now present isomorphism results for two classes of full shifts.
\bcor
Assume that $S_i$, $i=1,2$, are two non-trivial, countable, left reversible monoids which are cancellative, right LCM, satisfy $S_i^* = \gekl{1}$ as well as $\#(S_i \setminus s_i S_i) = \infty$ for all $1 \neq s_i \in S_i$, and that their enveloping groups $\msc{S}_i$ are amenable. Let $\Sigma_i$, $i=1,2$, be any two non-trivial finite Abelian groups. Consider the full shifts $\sigma_i: \: S_i \acts \bigoplus_{S_i} \Sigma_i$. Then we have $\fA_{\sigma_1} \cong \fA_{\sigma_2}$.
\ecor
\setlength{\parindent}{0cm} \setlength{\parskip}{0cm}

\bproof
This follows from Corollary~\ref{cor:FullShiftKirchberg}, \eqref{e:KThLCMNoUnits} and the Kirchberg-Phillips classification theorem \cite{KP,Phi}, together with the observation that for any infinite, countable monoid $S$ and any non-trivial, finite, Abelian group $\Sigma$, $C^*_{\lambda}(\bigoplus_S \Sigma)$ is isomorphic to the C*-algebra of continuous functions on the Cantor space, hence independent of $S$ and $\Sigma$.
\eproof
\setlength{\parindent}{0cm} \setlength{\parskip}{0.5cm}

\bcor
Assume that $S_i$, $i=1,2$, are two non-trivial, countable, left reversible monoids which are cancellative, right LCM, satisfy $S_i^* = \gekl{1}$, and that their enveloping groups $\msc{S}_i$ are amenable. Let $\Sigma$ be an arbitrary infinite, amenable group. Consider the full shifts $\sigma_i: \: S_i \acts \bigoplus_{S_i} \Sigma$. Then we have $\fA_{\sigma_1} \cong \fA_{\sigma_2}$.
\ecor
\setlength{\parindent}{0cm} \setlength{\parskip}{0cm}

\bproof
Apply Corollary~\ref{cor:FullShiftKirchberg}, \eqref{e:KThLCMNoUnits} and the Kirchberg-Phillips classification theorem \cite{KP,Phi}.
\eproof
\setlength{\parindent}{0cm} \setlength{\parskip}{0.5cm}

\subsubsection{Non-simple examples}
\label{sss:nonsimple}

Let us now use Theorem~\ref{thm:ClosedInvariantSubsp} to describe a class of algebraic actions whose C*-algebras have exactly one proper, non-zero ideal. Throughout this section, let us assume that the trivial subgroup is constructible, i.e., $\{e\}\in\cC$. This implies that the dense invariant subset $\fU:=\{\chi_k : k\in G\}$ of $\bd$ is open, so that $\fZ:=\bd\setminus\fU$ is a non-empty, proper, closed invariant subset of $\bd$.
Consider the following condition on $\sigma\colon S\acts G$:
\begin{equation}
	\label{eqn:(WPC)}\tag{WPC}
	\text{For every }C\in\cC \text{ with }C\neq\{e\},\text{ there exists }s\in S \text{ such that }\sigma_sG\subseteq C.
\end{equation}

\bprop
\label{lem:oneinvsubset}
Assume $\sigma\colon S\acts G$ satisfies \eqref{eqn:(WPC)} and that $\{e\}\in\cC$. Then
\setlength{\parindent}{0cm} \setlength{\parskip}{0cm}

\begin{enumerate}[\upshape(i)]
    \item $\fZ$ is the only non-empty, proper, closed invariant subset of $\bd$;
    \item the reduction groupoid $I\ltimes\fZ$ is purely infinite.
\end{enumerate} 
\eprop
\setlength{\parindent}{0cm} \setlength{\parskip}{0cm}

\bproof
(i): By Theorem~\ref{thm:ClosedInvariantSubsp}, this is equivalent to $\bm{\fF}_\cC\setminus \{\cC\}$ being the only non-empty, proper, $I^e$-invariant, $\subseteq$-closed subset of $\bm{\fF}_\cC$. 
Suppose $\emptyset\neq\bm{\fF}\subsetneq\bm{\fF}_\cC$ is an $I^e$-invariant and $\subseteq$-closed subset. Since $\bm{\fF}$ is $\subseteq$-closed and proper, $\bm{\fF}\subseteq\bm{\fF}_\cC\setminus \{\cC\}$. Let $C\in \cC\setminus\{e\}$, and choose any $\fF\in \bm{\fF}$. By \eqref{eqn:(WPC)}, there exists $s\in S$ such that $\sigma_sG\subseteq C$. Since $G\in\fF$, we have $\sigma_sG\in\sigma_s.\fF$, so that $C\in\sigma_s.\fF$ because $\sigma_s.\fF$ is a filter. By $I^e$-invariance of $\bm{\fF}$, we have $\sigma_s.\fF\in \bm{\fF}$, so $C\in \bigcup_{\fF\in\bm{\fF}}\fF$. Thus, $\cC\setminus\{\{e\}\}\subseteq \bigcup_{\fF\in\bm{\fF}}\fF$. Since $\bm{\fF}$ is $\subseteq$-closed, it follows that $\bm{\fF}$ contains every member of $\bm{\fF}_\cC\setminus \{\cC\}$.
\setlength{\parindent}{0cm} \setlength{\parskip}{0.5cm}

(ii): The basic (compact) open subsets of $\fZ$ are of the form $\fZ\cap \bd(kB;\{k_iB_i\})$, where $kB\in\cE^\times$ and $\{k_iB_i\}\subseteq\cE^\times$ is a finite subset. Moreover, each compact open subset of $\fZ$ can be written as a finite disjoint union of such sets (cf. \cite[Lemma~4.1]{Li:IMRN}), so to prove $I\ltimes\fZ$ is purely infinite it suffices to prove that every non-empty compact open subset of the form $\fZ\cap \bd(kB;\{k_iB_i\})$ is properly infinite.
\setlength{\parindent}{0.5cm} \setlength{\parskip}{0cm}

If $kB=\{k\}$, then $ \bd(kB)=\{\chi_k\}$, so $\fZ\cap \bd(kB;\{k_iB_i\})=\emptyset$ in this case. If $k_iB_i=\{k_i\}$, then $\bd(k_jB_j)=\{\chi_{k_j}\}$, so that $\fZ\cap \bd(kB;\{k_iB_i\})=\fZ\cap \bd(kB;\{k_iB_i : i\neq j \})$. Thus, we may even assume that $\{e\}\notin\{B\}\cup\{B_i\}_i$. In this case, the proof of Theorem~\ref{thm:PI} goes through and gives that $\bd(kB;\{k_iB_i\})$ is properly infinite in $\cG_\sigma$, and thus $\fZ\cap \bd(kB;\{k_iB_i\})$ is properly infinite in $I\ltimes\fZ$. 
\eproof
\setlength{\parindent}{0cm} \setlength{\parskip}{0.5cm}

\blemma
If $\sigma\colon S\acts G$ satisfies \eqref{eqn:(WPC)} and $\bigcap_{s\in S}\sigma_sG=\{e\}$, then $I\ltimes\fZ$ is topologically free. 
\elemma
\setlength{\parindent}{0cm} \setlength{\parskip}{0cm}

\bproof
By Theorem~\ref{thm:TopFree}, it suffices to show that for all $D \in \bigcup_{\cC\neq \fF \in \bm{\fF}_\cC} \fF$, we have $\bigcap_{\cC\neq\fF \in \bm{\fF}_\cC,\, D \in \fF}\cap\fF = \gekl{e}$. Given $D \in \bigcup_{\cC\neq \fF \in \bm{\fF}_\cC} \fF$, \eqref{eqn:(WPC)} implies that there exists $t \in S$ with $\sigma_t G \subseteq D$. It follows that for all $s \in S$, $\sigma_{ts}G \subseteq \sigma_t G \subseteq D$. Thus $D \in \sigma_{ts}.\fF$ because $\sigma_{ts} G \in \sigma_{ts}.\fF$ for all $\fF \in \bm{\fF}_\cC$. It follows that, for all $\fF \in \bm{\fF}_\cC$, $\bigcap_{\cC\neq\fF \in \bm{\fF}_\cC,\, D \in \fF}\cap\fF\subseteq \bigcap_{s \in S} \sigma_{ts} G = \sigma_t(\bigcap_{s \in S} \sigma_s G) = \sigma_t(\gekl{e}) = \gekl{e}$.
\eproof
\setlength{\parindent}{0cm} \setlength{\parskip}{0.5cm}

\bcor
\label{cor:nonsimple}
If $\{e\}\in\cC$, $\sigma\colon S\acts G$ satisfies \eqref{eqn:(WPC)} and $\bigcap_{s\in S}\sigma_sG=\{e\}$, then the canonical map $C_\es^*(\cG_\sigma)\to \fA_\sigma$ from Corollary~\ref{cor:rho<Indchi} is an isomorphism, under which $C_r^*(I\ltimes\fU) = C_{\es}^*(I\ltimes\fU)$ is identified with $\cK(\ell^2(G))$, and $C^*_\es(I\ltimes\fZ)$ is simple and purely infinite. If moreover $S$ and $G$ are countable and $\cG_\sigma$ is Hausdorff and inner exact, then $\cK(\ell^2(G))$ is the unique non-zero, proper ideal of $\fA_{\sigma}$ and $\fA_{\sigma} / \cK(\ell^2(G)) \cong C_r^*(\cG_\sigma)/C_r^*(I\ltimes\fU) \cong C^*_r(I\ltimes\fZ)$ is simple and purely infinite.
\ecor

\bremark
\label{rmk:WPC}
Assume $\sigma\colon S\acts A$ satisfies the following condition: 
\begin{equation}
	\label{lem:enoughforWPC}
	\sigma_{s_1}^{-1}\sigma_{t_1}\cdots\sigma_{s_m}^{-1}\sigma_{t_m}G\neq\{e\}\implies s_1^{-1}t_1\cdots s_m^{-1}t_mS\neq\emptyset
\end{equation}
for all $s_i,t_i\in S$ and $m\in\Zz_{>0}$. Arguing as in Proposition~\ref{prop:exactness4Sreversible}, we see $\sigma\colon S\acts A$ satisfies \eqref{eqn:(WPC)}.
\eremark

\bex
\label{ex:FullShiftNotRev}
Let $\sigma\colon S\acts\bigoplus_S\Sigma$ be the full $S$-shift over $\Sigma$ as in Definition~\ref{def:genshift}.  If $S$ is not left reversible, then $\{e\}$ is constructible. Using Remark~\ref{rmk:WPC}, it is not hard to see that $\sigma\colon S\acts\bigoplus_S\Sigma$ satisfies \eqref{eqn:(WPC)}. We have seen in the proof of Corollary~\ref{cor:genshifts} that $\bigcap_{s\in S}\sigma_s(\bigoplus_S\Sigma)=\gekl{e}$. Hence Corollary~\ref{cor:nonsimple} applies.
\eex

\bremark
Let $\sigma_i \colon S_i \acts G_i$, $i=1,2$, be two algebraic actions. Form the product action $\sigma_1 \times \sigma_2 \colon S_1 \times S_2 \acts G_1 \times G_2$. Then it is easy to see that $\fA_{\sigma_1 \times \sigma_2} \cong \fA_{\sigma_1} \otimes_{\min} \fA_{\sigma_2}$. Hence, by forming products of algebraic actions as in Example~\ref{ex:FullShiftNotRev}, we obtain algebraic actions whose C*-algebras have more complicated ideal structures.
\eremark


\begin{thebibliography}{99}
	


\bibitem{AR} C. \textsc{Anantharaman-Delaroche} and J. \textsc{Renault}, \emph{Amenable groupoids}, Groupoids in analysis, geometry, and physics (Boulder, CO, 1999), 35--46, Contemp. Math., 282, Amer. Math. Soc., Providence, RI, 2001.

	
\bibitem{AM} M. F. \textsc{Atiyah} and I. G. \textsc{Macdonald}, \emph{Introduction to commutative algebra}, Addison-Wesley Publishing Co., Reading, Mass.-London-Don Mills, Ont. 1969.

	
\bibitem{Ber} D. \textsc{Berend}, \emph{Multi-invariant sets on tori}, Trans. Amer. Math. Soc. \emph{280} (1983), no. 2, 509--532.


\bibitem{BoLi} C. \textsc{B\"{o}nicke} and K. \textsc{Li}, \emph{Ideal structure and pure infiniteness of ample groupoid C*-algebras}, Ergodic Theory Dynam. Systems \emph{40} (2020), no. 1, 34--63. 


\bibitem{BS} N. \textsc{Brownlowe} and N. \textsc{Stammeier}, \emph{The boundary quotient for algebraic dynamical systems}, J. Math. Anal. Appl. \emph{438} (2016), no. 2, 772--789.

\bibitem{BLS} N. \textsc{Brownlowe}, N. S. \textsc{Larsen}, and N. \textsc{Stammeier}, \emph{On C*-algebras associated to right LCM semigroups}, Trans. Amer. Math. Soc. \emph{369} (2017), no. 1, 31--68.

\bibitem{BLS2} N. \textsc{Brownlowe}, N. S. \textsc{Larsen}, and N. \textsc{Stammeier}, \emph{C*-algebras of algebraic dynamical systems and right LCM semigroups}, Indiana Univ. Math. J. \emph{67} (2018), no. 6, 2453--2486.

\bibitem{Bru1} C. \textsc{Bruce}, \emph{C*-algebras from actions of congruence monoids on rings of algebraic integers}, Trans. Amer. Math. Soc. \emph{373} (2020), no. 1, 699--726. 

\bibitem{BruLi} C. \textsc{Bruce} and X. \textsc{Li}, \emph{On K-theoretic invariants of semigroup C*-algebras from actions of congruence monoids}, Amer. J. Math. \emph{145} (2023), no. 1, 251--285.

\bibitem{BruLi:rigidity} C. \textsc{Bruce} and X. \textsc{Li}, \emph{Algebraic actions II. Groupoid rigidity}, preprint, \href{https://arxiv.org/abs/2301.04459}{arXiv:2301.04459}.

\bibitem{BFS} A. \textsc{Buss}, D. \textsc{Ferraro}, and C. F. \textsc{Sehnem}, \emph{Nuclearity for partial crossed products by exact discrete groups}, J. Operator Theory \emph{88} (2022), no. 1, 83--115.

\bibitem{CN} J. \textsc{Christensen} and S. \textsc{Neshveyev}, \emph{(Non)exotic Completions of the Group Algebras of Isotropy Groups}, Int. Math. Res. Not. IMRN 2022, no. 19, 15155--15186.

\bibitem{CP67} A. H. \textsc{Clifford} and G. B. \textsc{Preston}, \emph{The algebraic theory of semigroups}, Vol. II. Mathematical Surveys, No. 7 American Mathematical Society, Providence, R.I. 1967.

\bibitem{Connes} A. \textsc{Connes}, \emph{A survey of foliations and operator algebras}, Operator algebras and applications,
Part I (Kingston, Ont., 1980). Vol. 38. Proc. Sympos. Pure Math. Amer. Math. Soc.,
1982, pp.


\bibitem{Cuntz} J. \textsc{Cuntz}, \emph{C*-algebras associated with the $ax+b$-semigroup over $\Nz$}, $K$-theory and noncommutative geometry, 201--215, EMS Ser. Congr. Rep., Eur. Math. Soc., Z\"urich, 2008.

\bibitem{CELY} J. \textsc{Cuntz}, S. \textsc{Echterhoff}, X. \textsc{Li}, and G. \textsc{Yu}, \emph{K-theory for group C*-algebras and semigroup C*-algebras}, Oberwolfach Seminars, 47, Birkh{\"a}user/Springer, Cham, 2017.


\bibitem{CuLi} J. \textsc{Cuntz} and X. \textsc{Li}, \emph{The regular C*-algebra of an integral domain}, Quanta of Maths, Clay Math. Proc., Vol. 11, Amer. Math. Soc., 2010, pp. 149--170.

\bibitem{CV} J. \textsc{Cuntz} and A. \textsc{Vershik}, \emph{C*-algebras associated with endomorphisms and polymorphisms of compact Abelian groups}, Comm. Math. Phys. \emph{321} (2013), no. 1, 157--179.


\bibitem{Eis} D. \textsc{Eisenbud}, \emph{Commutative algebra. With a view toward algebraic geometry}, Graduate Texts in Mathematics, 150. Springer-Verlag, New York, 1995.

\bibitem{ExelPAMS} R. \textsc{Exel}, \emph{Partial actions of groups and actions of inverse semigroups}, Proc. Amer. Math. Soc. \emph{126} (1998), no. 12, 3481--3494.

\bibitem{Exel08} R. \textsc{Exel}, \emph{Inverse semigroups and combinatorial C*-algebras}, Bull. Braz. Math. Soc. (N.S.) \emph{39} (2008), no. 2, 191--313.

\bibitem{Exel09} R. \textsc{Exel}, \emph{Tight representations of semilattices and inverse semigroups}, Semigroup Forum \emph{79} (2009), no. 1, 159--182. 

\bibitem{Exel21} R. \textsc{Exel}, \emph{Tight and cover-to-join representations of semilattices and inverse semigroups}, Operator theory, functional analysis and applications, 183--192, Oper. Theory Adv. Appl., 282, Birkh\"{a}user/Springer, Cham, 2021.

\bibitem{ExelPardo} R. \textsc{Exel} and E. \textsc{Pardo}, \emph{The tight groupoid of an inverse semigroup}, Semigroup Forum \emph{92} (2016), no. 1, 274--303.

\bibitem{ExelPitts} R. \textsc{Exel} and D. R. \textsc{Pitts}, \emph{Characterizing Groupoid C*-algebras of Non-Hausdorff {\'E}tale Groupoids}, Lecture Notes in Mathematics, 2306, Springer, Cham, 2022.

\bibitem{Fuchs2} L. \textsc{Fuchs}, \emph{Infinite Abelian groups. Vol. II}, Pure and Applied Mathematics. Vol. 36-II. Academic Press, New York-London, 1973.

\bibitem{GT} E. \textsc{Gardella} and O. \textsc{Tanner}, \emph{Generalisations of Thompson's group V arising from purely infinite groupoids}, preprint, \href{https://arxiv.org/abs/2302.04078}{arXiv:2302.04078}.


\bibitem{Hirsh} I. \textsc{Hirshberg}, \emph{On C*-algebras associated to certain endomorphisms of discrete groups}, New York J. Math. \emph{8} (2002), 99--109.

\bibitem{Hunger} T. W. \textsc{Hungerford}, \emph{Algebra}, Reprint of the 1974 original. Graduate Texts in Mathematics, 73. Springer-Verlag, New York-Berlin, 1980.

\bibitem{KellLaw} J. \textsc{Kellendonk} and M. V. \textsc{Lawson}, \emph{Partial actions of groups}, Internat. J. Algebra Comput. 14 (2004), no. 1, 87--114. 

\bibitem{KL} D. \textsc{Kerr} and H. \textsc{Li}, \emph{Ergodic theory. Independence and dichotomies}, Springer Monographs in Mathematics. Springer, Cham, 2016.


\bibitem{KKLRU} M. \textsc{Kennedy}, S. J. \textsc{Kim}, X. \textsc{Li}, S. \textsc{Raum}, and D. \textsc{Ursu}, \emph{The ideal intersection property for essential groupoid C*-algebras}, preprint, arXiv:2107.03980. 

\bibitem{Ker} D. \textsc{Kerr}, \emph{Dimension, comparison, and almost finiteness}, J. Eur. Math. Soc. \emph{22} (2020), no. 11, 3697--3745. 

\bibitem{KS} M. \textsc{Khoshkam} and G. \textsc{Skandalis}, \emph{Regular representation of groupoid C*-algebras and applications to inverse semigroups}, J. Reine Angew. Math. \emph{546} (2002), 47--72. 

\bibitem{KP} E. \textsc{Kirchberg} and N.C. \textsc{Phillips}, \emph{Embedding of exact C*-algebras in the Cuntz algebra $\cO_2$}, J. Reine Angew. Math. \emph{525} (2000), 17--53.

\bibitem{Krz} K. \textsc{Krzy{\.z}ewski}, \emph{On exact toral endomorphisms}, Monatsh. Math. \emph{116} (1993), no. 1, 39--47.

\bibitem{KM} B. K. \textsc{Kwa\'sniewski} and R. \textsc{Meyer}, \emph{Essential crossed products for inverse semigroup actions: simplicity and pure infiniteness}, Doc. Math. \emph{26} (2021), 271--335.

\bibitem{LaSe} M. \textsc{Laca} and C. F. \textsc{Sehnem}, \emph{Toeplitz algebras of semigroups}, Trans. Amer. Math. Soc. \emph{375} (2022), no. 10, 7443--7507.

\bibitem{Law}M. V. \textsc{Lawson}, \emph{Inverse semigroups. The theory of partial symmetries}, World Scientific Publishing Co., Inc., River Edge, NJ, 1998.

\bibitem{LW} M. V. \textsc{Lawson} and A. R. \textsc{Wallis}, \emph{A correspondence between a class of monoids and self-similar group actions II}, Internat. J. Algebra Comput. \emph{25} (2015), no. 4, 633--668.


\bibitem{Li:Ring} X. \textsc{Li}, \emph{Ring C*-algebras}, Math. Ann. \emph{348} (2010), no. 4, 859--898.

\bibitem{Li:JFA} X. \textsc{Li}, \emph{Semigroup C*-algebras and amenability of semigroups}, J. Funct. Anal. 262 (2012), no. 10, 4302--4340.



\bibitem{Li:IMRN} X. \textsc{Li}, \emph{Partial transformation groupoids attached to graphs and semigroups}, Int. Math. Res. Not. IMRN 2017, no. 17, 5233--5259.

\bibitem{Li:Lamp} X. \textsc{Li}, \emph{$K$-theory for generalized Lamplighter groups}, Proc. Amer. Math. Soc. \emph{147} (2019), no. 10, 4371--4378. 

\bibitem{Li:Ktheory}  X. \textsc{Li}, \emph{K-theory for semigroup C*-algebras and partial crossed products},  Comm. Math. Phys. \emph{390} (2022), no. 1, 1--32.

\bibitem{Li:GarI} X. \textsc{Li}, \emph{Left regular representations of Garside categories I. C*-algebras and groupoids}, Glasg. Math. J. \emph{65} (2023), no. S1, S53--S86.

\bibitem{Li:TFG} X. \textsc{Li}, \emph{Ample groupoids, topological full groups, algebraic K-theory spectra and infinite loop spaces}, preprint,  \href{https://arxiv.org/abs/2209.08087}{arXiv:2209.08087}.


\bibitem{LOS} X. \textsc{Li}, T. \textsc{Omland} and J. \textsc{Spielberg}, \emph{C*-algebras of right LCM one-relator monoids and Artin-Tits monoids of finite type}, Comm. Math. Phys. \emph{381} (2021), no. 3, 1263--1308. 

\bibitem{LindSchmidt} D. \textsc{Lind} and K. \textsc{Schmidt}, \emph{A survey of algebraic actions of the discrete Heisenberg group}. (Russian); translated from Uspekhi Mat. Nauk \emph{70} (2015), no. 4(424), 77--142 Russian Math. Surveys \emph{70} (2015), no. 4, 657--714 

\bibitem{Ma} X. \textsc{Ma}, \emph{Purely infinite locally compact Hausdorff \'{e}tale groupoids and their C*-algebras},  Int. Math. Res. Not. IMRN 2022, no. 11, 8420--8471.

\bibitem{MW} X. \textsc{Ma} and J. \textsc{Wu}, \emph{Almost elementariness and fiberwise amenability for \'{e}tale groupoids}, preprint, \href{https://arxiv.org/abs/2011.01182}{arXiv:2011.01182}.

\bibitem{Mats2} H. \textsc{Matsumura}, \emph{Commutative ring theory}. Second edition. Cambridge Studies in Advanced Mathematics, 8. Cambridge University Press, Cambridge, 1989.

\bibitem{Mat12} H. Matui, \emph{Homology and topological full groups of \'{e}tale groupoids on totally disconnected spaces}, Proc. Lond. Math. Soc. (3) 104 (2012), no. 1, 27--56.

\bibitem{Mat15} H. \textsc{Matui}, \emph{Topological full groups of one-sided shifts of finite type}, J. Reine Angew. Math. \emph{705} (2015), 35--84.

\bibitem{Muc} R. \textsc{Muchnik}, \emph{Semigroup actions on $\Tz^n$}, Geom. Dedicata \emph{110} (2005), 1--47.


\bibitem{Nek} V. \textsc{Nekrashevych}, \emph{Self-similar groups}. Mathematical Surveys and Monographs, 117. American Mathematical Society, Providence, RI, 2005.


\bibitem{NS} S. \textsc{Neshveyev} and G. \textsc{Schwartz}, \emph{Non-Hausdorff etale groupoids and C*-algebras of left cancellative monoids}, M\"{u}nster J. Math. \emph{16} (2023), no. 1, 147--175.


\bibitem{Neu} B. H. \textsc{Neumann}, \emph{Groups covered by permutable subsets}, J. London Math. Soc. \emph{29} (1954), 236--248.

\bibitem{CEPSS} L. O. \textsc{Clark}, R. \textsc{Exel}, E. \textsc{Pardo}, A. \textsc{Sims}, and C. \textsc{Starling}, \emph{Simplicity of algebras associated to non-Hausdorff groupoids}, Trans. Amer. Math. Soc. \emph{372} (2019), no. 5, 3669--3712.

\bibitem{Pat} A.L.T. \textsc{Paterson}, \emph{Groupoids, inverse semigroups, and their operator algebras}, Progress in Mathematics, 170. Birkh{\"a}user Boston, Inc., Boston, MA, 1999.

\bibitem{Phi} N. C. \textsc{Phillips}, \emph{A classification theorem for nuclear purely infinite simple C*-algebras}, Doc. Math. \emph{5} (2000), 49--114.

\bibitem{RS} T. \textsc{Rainone} and A. \textsc{Sims}, \emph{A dichotomy for groupoid C*-algebras}, Ergodic Theory Dynam. Systems \emph{40} (2020), no. 2, 521--563.

\bibitem{Renault} J. \textsc{Renault}, \emph{A groupoid approach to C*-algebras}, Lecture Notes in Mathematics, 793. Springer, Berlin, 1980.


\bibitem{Renault2} J. \textsc{Renault}, \emph{C*-algebras and dynamical systems}, Publicac\~{o}es Matem\'{a}ticas do IMPA. [IMPA Mathematical Publications] 27$\sp {\rm o}$ Col\'{o}quio Brasileiro de Matem\'{a}tica. [27th Brazilian Mathematics Colloquium] Instituto Nacional de Matem\'{a}tica Pura e Aplicada (IMPA), Rio de Janeiro, 2009.

\bibitem{Rohlin} V. A. \textsc{Rohlin}, \emph{Exact endomorphisms of a Lebesgue space}, (Russian) Izv. Akad. Nauk SSSR Ser. Mat. \emph{25} (1961) 499--530. 

\bibitem{Sch} K. \textsc{Schmidt}, \emph{Dynamical systems of algebraic origin}, Progress in Mathematics, 128, Birkh{\"a}user Verlag, Basel, 1995.


\bibitem{SSW} A. \textsc{Sims}, G. \textsc{Szab\'{o}}, and D. \textsc{Williams}, \emph{Operator Algebras and Dynamics: Groupoids, Crossed Products, and Rokhlin
Dimension}, Advanced Courses in Mathematics. CRM Barcelona, Birkh\"{a}user, 2020.


\bibitem{Sta} N. \textsc{Stammeier}, \emph{On C*-algebras of irreversible algebraic dynamical systems}, J. Funct. Anal. \emph{269} (2015), no. 4, 1136--1179.

\bibitem{Sta2} N. \textsc{Stammeier}, \emph{Topological freeness for $\ast$-commuting covering maps}, Houston J. Math. \emph{45} (2019), no. 2, 525--551.

\bibitem{Vieira1} F. \textsc{Vieira}, \emph{C*-algebras associated with endomorphisms of groups}, J. Operator Theory \emph{79} (2018), no. 1, 3--31.



\end{thebibliography}
\end{document}